\documentclass[a4paper,12pt,amsart,frenchb]{article}

\usepackage{amsmath,amsbsy,amsfonts,amssymb}
\usepackage[french]{babel}
\newcommand{\ass}[2]{\vskip0.3cm\noindent
{\bf {#1}}. { \sl {#2}}\vskip0.3cm\noindent
}
  
\begin{document}
 
   \title{Stabilisation de la formule des traces tordue IV:    transfert spectral archim\'edien}
\author{J.-L. Waldspurger}
\date{6 mars 2014}
\maketitle

\bigskip

{\bf Introduction}

Comme la r\'ef\'erence [I], dont nous reprenons les notations, cet article  contient des r\'esultats pr\'eparatoires \`a la stabilisation de la formule des traces tordue. Il concerne  exclusivement les groupes r\'eels. La premi\`ere section \'enonce un th\'eor\`eme de Paley-Wiener pour les fonctions $C^{\infty}$ \`a support compact. Ce th\'eor\`eme est d\^u \`a Renard mais on en modifie quelque peu la formulation pour y faire appara\^{\i}tre les repr\'esentations elliptiques qui, depuis le travail qu'Arthur leur a consacr\'e, sont devenus les blocs de base de ce type d'analyse harmonique. Dans les sections 2 et 3, on prouve les analogues dans le cas tordu, et sur le corps de base ${\mathbb R}$, des r\'esultats contenus dans l'article [A] d'Arthur. A savoir les deux r\'esultats suivants, exprim\'es ici de fa\c{c}on lapidaire. On consid\`ere un triplet $(G,\tilde{G},{\bf a})$, o\`u $G$ est un groupe r\'eductif connexe d\'efini sur ${\mathbb R}$, $\tilde{G}$ est un espace tordu sur $G$ et ${\bf a}$ est un \'el\'ement de $H^1(W_{{\mathbb R}};Z(\hat{G}))$
Supposons d'abord $G$ quasi-d\'eploy\'e, $\tilde{G}$  \`a torsion int\'erieure et ${\bf a}=1$. Alors une combinaison lin\'eaire finie de caract\`eres de repr\'esentations elliptiques de $\tilde{G}({\mathbb R})$ qui est stable sur les \'el\'ements r\'eguliers elliptiques est stable partout.  Dans le cas g\'en\'eral, soit ${\bf G}'$ une donn\'ee endoscopique elliptique de $(G,\tilde{G},{\bf a})$. Consid\'erons une combinaison lin\'eaire finie $\Sigma$ de caract\`eres de repr\'esentations elliptiques de $\tilde{G}'({\mathbb R})$ et une combinaison lin\'eaire finie $\Pi$ de caract\`eres de $\omega$-repr\'esentations elliptiques de $\tilde{G}({\mathbb R})$. Supposons que $\Sigma$ est stable et que la restriction de $\Pi$ aux \'el\'ements r\'eguliers elliptiques de $\tilde{G}({\mathbb R})$ est \'egale \`a la m\^eme restriction du transfert de $\Sigma$. Alors $\Pi$ est le transfert de $\Sigma$. On a n\'eglig\'e ici comme dans la suite de cette introduction le fait qu'en g\'en\'eral, il faut remplacer $\tilde{G}'({\mathbb R})$ par un espace auxiliaire $\tilde{G}'_{1}({\mathbb R})$. 
Une premi\`ere cons\'equence de ces r\'esultats est une version "stable" du th\'eor\`eme de Paley-Wiener (th\'eor\`eme 2.3(ii)). Une deuxi\`eme est  la d\'efinition du transfert spectral (corollaire 3.3). Une troisi\`eme cons\'equence est l'existence du transfert g\'eom\'etrique $K$-fini: si $f\in C_{c}^{\infty}(\tilde{G}({\mathbb R}))$ est $K$-finie, il existe une fonction $f'\in C_{c}^{\infty}(\tilde{G}'({\mathbb R}))$  qui est $K'$-finie et qui est un transfert de $f$, cf. corollaire 3.4.  Tout cela est certainement cons\'equence des r\'esultats beaucoup plus fins obtenus par Mezo dans son article  r\'ecent [M]. Nos preuves sont tr\`es diff\'erentes. Elles s'appuient sur le th\'eor\`eme de Renard repris dans la premi\`ere section, sur le r\'esultat de Shelstad affirmant l'existence du transfert entre fonctions $C^{\infty}$ \`a support compact et sur le r\'esultat suivant: une combinaison lin\'eaire finie de caract\`eres de repr\'esentations temp\'er\'ees de $\tilde{G}({\mathbb R})$ est supertemp\'er\'ee si et seulement si toutes les repr\'esentations qui interviennent sont elliptiques. Dans le cas non tordu,  ce r\'esultat est d\^u \`a Harish-Chandra. Il vaut aussi d'apr\`es Herb sur un corps de base non-archim\'edien. Il a \'et\'e r\'ecemment g\'en\'eralis\'e au cas tordu par Moeglin dans [Moe], que le corps de base soit r\'eel ou non-archim\'edien.

 \bigskip
\section{Le th\'eor\`eme de Paley-Wiener pour les fonctions $C^{\infty}$ \`a support compact}
\bigskip

\subsection{La situation}
Dans cet article, le corps de base est ${\mathbb R}$.     On consid\`ere un triplet $(G,\tilde{G},{\bf a})$, o\`u $G$ est un groupe r\'eductif connexe d\'efini sur ${\mathbb R}$, $\tilde{G}$ est un espace tordu sur $G$ et ${\bf a}$ est un \'el\'ement de $H^1(W_{{\mathbb R}};Z(\hat{G}))$, cf. [I] 1.1 et 1.5. Le terme ${\bf a}$ d\'etermine un caract\`ere $\omega$   de $G({\mathbb R})$. On suppose

$\bullet$ $\tilde{G}({\mathbb R})\not=\emptyset$;

$\bullet$ l'automorphisme $\theta$ de $Z(G)$ est d'ordre fini.

On fixe un espace de Levi minimal $\tilde{M}_{0}$ de $\tilde{G}$ et un sous-groupe compact maximal $K$ de $G({\mathbb R})$. On suppose que les alg\`ebres de Lie $\mathfrak{k}$ de $K$ et $\mathfrak{a}_{M_{0}}$ de $A_{M_{0}}$ sont orthogonales pour la forme de Killing. Introduisons $\underline{la}$ paire de Borel $(B^*,T^*)$ de $G$. On pose
$\mathfrak{h}_{{\mathbb R}}=X_{*}(T^*)\otimes{\mathbb R}$
et $\mathfrak{h}=\mathfrak{h}_{{\mathbb R}}\otimes {\mathbb C}$. Soit $\tilde{S}$ un sous-tore tordu maximal de $G$ d\'efini sur ${\mathbb R}$.
Le groupe $\Gamma_{{\mathbb R}}\simeq \{\pm 1\}$ agit sur $X_{*}(S)$ et sur $X_{*,{\mathbb Q}}(S)=X_{*}(S)\otimes_{{\mathbb Z}}{\mathbb Q}$.  De cette action  se d\'eduit une d\'ecomposition
$$X_{*,{\mathbb Q}}(S)=X_{*,{\mathbb Q}}(S) ^+\oplus X_{*,{\mathbb Q}}(S)^-$$
o\`u $\Gamma_{{\mathbb R}}$ agit trivialement sur le premier sous-espace et par le caract\`ere non trivial sur le second. On en d\'eduit une d\'ecomposition de
 l'alg\`ebre de Lie $\mathfrak{s}({\mathbb R})=(X_{*}(S)\otimes_{{\mathbb Z}}{\mathbb C})^{\Gamma_{{\mathbb R}}}$  en somme directe 
 $$\mathfrak{s}({\mathbb R})= X_{*,{\mathbb Q}}(S)^+\otimes_{{\mathbb Q}}{\mathbb R}\oplus X_{*,{\mathbb Q}}(S)^-\otimes_{{\mathbb Q}}i{\mathbb R}.$$
 Modulo le choix d'un groupe de Borel contenant $S$ et stable par $\tilde{S}$, on peut identifier $\mathfrak{s}$ \`a $\mathfrak{h}$. Le premier facteur ci-dessus s'identifie \`a un sous-espace de $\mathfrak{h}_{{\mathbb R}}$ et le second s'identifie \`a un sous-espace de $i\mathfrak{h}_{{\mathbb R}}$.   Si on prend pour $\tilde{S}$ un sous-tore tordu maximal de $\tilde{M}_{0}$, le premier facteur n'est autre que $\mathfrak{a}_{\tilde{M}_{0}}({\mathbb R})$, qui s'identifie ainsi \`a un sous-espace de $\mathfrak{h}_{{\mathbb R}}$. 
 
   Le groupe $\Gamma_{{\mathbb R}}$ agit sur $T^*$ donc aussi sur $\mathfrak{h}_{{\mathbb R}}$. On fixe sur $\mathfrak{h}_{\mathbb R}$ une forme quadratique d\'efinie positive invariante par l'action du groupe de Weyl de $G$ relatif \`a $T^*$, par celle de $\Gamma_{{\mathbb R}}$  et par l'automorphisme $\theta$ de $T^*$.  Par dualit\'e, on en d\'eduit une telle forme sur le dual $\mathfrak{h}_{\mathbb R}^*$.  Il se d\'eduit aussi de $\theta$ un automorphisme dual de $\mathfrak{h}_{\mathbb R}^*$ que l'on note encore $\theta$. Remarquons que $\mathfrak{a}_{M_{0}}({\mathbb R})$ s'identifie \`a ${\cal A}_{M_{0}}$. Pour tout Levi $M$ de $G$ contenant $M_{0}$, ${\cal A}_{M}$ est un sous-espace de ${\cal A}_{M_{0}}$ et on munit cet espace de la restriction de la forme quadratique. Plus g\'en\'eralement, pour tout Levi $M$, on peut choisir $g\in G({\mathbb R})$ tel que $g^{-1}Mg$ contienne $M_{0}$ et on munit ${\cal A}_{M}$ de la forme quadratique sur ${\cal A}_{g^{-1}Mg}$ transport\'ee par l'isomorphisme $ad_{g}$. Cela ne d\'epend pas du choix de $g$. On munit tout sous-espace de ${\cal A}_{M}$  de la mesure de Haar associ\'ee \`a la restriction \`a ce sous-espace de cette forme quadratique. 

On \'etend toutes ces formes quadratiques en des formes hermitiennes sur les complexifi\'es de nos espaces. Pour tout espace vectoriel $V$ sur ${\mathbb R}$, on note $V_{{\mathbb C}}=V\otimes_{{\mathbb R}}{\mathbb C}$ son complexifi\'e. On fait une exception pour l'espace $\mathfrak{h}_{{\mathbb R}}$ dont on note simplement $\mathfrak{h}$ le complexifi\'e. 

On note $\mathfrak{A}_{G}$, resp. $\mathfrak{A}_{\tilde{G}}$, la composante neutre pour la topologie r\'eelle de $A_{G}({\mathbb R})$, resp. $A_{\tilde{G}}({\mathbb R})$. On sait que la restriction \`a $\mathfrak{A}_{G}$ de l'application $H_{G}:G({\mathbb R})\to {\cal A}_{G}$   est un isomorphisme. On munit $\mathfrak{A}_{G}$, resp $\mathfrak{A}_{\tilde{G}}$, de la mesure de Haar d\'eduite par cet isomorphisme de celle fix\'ee sur ${\cal A}_{G}$, resp. ${\cal A}_{\tilde{G}}$. 

  L'homomorphisme
 $H_{\tilde{G}}:G({\mathbb R})\to {\cal A}_{\tilde{G}}$ est le compos\'e de $H_{G}$ et de la projection ${\cal A}_{G}\to {\cal A}_{\tilde{G}}$ sur les invariants par $\theta$. On fixe arbitrairement une application encore not\'ee $H_{\tilde{G}}:\tilde{G}({\mathbb R})\to {\cal A}_{\tilde{G}}$  telle que 
  $H_{\tilde{G}}(g\gamma)=H_{\tilde{G}}(g)+H_{\tilde{G}}(\gamma)$ pour tous $g\in G({\mathbb R})$ et $\gamma\in \tilde{G}({\mathbb R})$.
Notons $\tilde{G}({\mathbb R})^1$ l'ensemble des $\gamma\in \tilde{G}({\mathbb R})$ tels que $H_{\tilde{G}}(\gamma)=0$. On a des isomorphismes inverses l'un de l'autre
$$\begin{array}{ccccccc}\tilde{G}({\mathbb R}) &\to&{\cal A}_{\tilde{G}}\times \tilde{G}({\mathbb R})^1&\qquad&{\cal A}_{\tilde{G}}\times \tilde{G}({\mathbb R})^1&\to& \tilde{G}({\mathbb R})\\ \gamma&\mapsto &(H_{\tilde{G}}(\gamma),exp(-H_{\tilde{G}}(\gamma))\gamma)&\qquad &(H,\gamma^1)&\mapsto& e^H\gamma^1\\ \end{array}$$

Posons $W(\tilde{M}_{0})=Norm_{G({\mathbb R})}(\tilde{M}_{0})/M_{0}({\mathbb R})$. On a

(1) il existe une unique application $H_{\tilde{M}_{0}}:\tilde{M}_{0}({\mathbb R})\to {\cal A}_{\tilde{M}_{0}}$ telle que

(i) $H_{\tilde{M}_{0}}(m\gamma)=H_{\tilde{M}_{0}}(m)+H_{\tilde{M}_{0}}(\gamma)$ pour tous $m\in M_{0}({\mathbb R})$ et $\gamma\in \tilde{M}_{0}({\mathbb R})$;

(ii) $H_{\tilde{M}_{0}}(ad_{n}(\gamma))=w(H_{\tilde{M}_{0}}(\gamma))$, pour tous $\gamma\in \tilde{M}_{0}({\mathbb R})$, $n\in Norm_{G({\mathbb R})}(\tilde{M}_{0})$, o\`u $w$ est l'image de $n$ dans $W(\tilde{M}_{0})$;

(iii) la compos\'ee de $H_{\tilde{M}_{0}}$ et de la projection ${\cal A}_{\tilde{M}_{0}}\to {\cal A}_{\tilde{G}}$ est la restriction de $H_{\tilde{G}}$ \`a $\tilde{M}_{0}({\mathbb R})$. 

Preuve.  Fixons $\gamma_{0}\in \tilde{M}_{0}({\mathbb R})$. La fonction $H'_{\tilde{M}_{0}}$ d\'efinie par $H'_{\tilde{M}_{0}}(m\gamma_{0})=H_{\tilde{M}_{0}}(m)+H_{\tilde{G}}(\gamma_{0})$ pour tout $m\in M_{0}({\mathbb R})$ v\'erifie (i) et (iii). Posons 
$$H_{\tilde{M}_{0}}(\gamma)=\vert W(\tilde{M}_{0})\vert^{-1} \sum_{w\in W(\tilde{M}_{0})}w^{-1}H'_{\tilde{M}_{0}}(ad_{n_{w}}(\gamma)),$$
o\`u $n_{w}$ est un rel\`evement de $w$ dans $Norm_{G({\mathbb R})}(\tilde{M}_{0})$. Cette fonction r\'epond \`a la question. Si $H''_{\tilde{M}_{0}}$ v\'erifie les m\^emes conditions, la condition (i) entra\^{\i}ne qu'il existe $H\in {\cal A}_{\tilde{M}_{0}}$ tel que $H''_{\tilde{M}_{0}}(\gamma)=H+H_{\tilde{M}_{0}}(\gamma)$ pour tout $\gamma\in \tilde{M}_{0}({\mathbb R})$. La condition (ii) entra\^{\i}ne que $H$ est invariant par $W(\tilde{M}_{0})$, donc appartient \`a ${\cal A}_{\tilde{G}}$. La condition (iii) entra\^{\i}ne alors que $H=0$. $\square$

On d\'efinit $H_{\tilde{M}_{0}}$ par la condition (1). Plus g\'en\'eralement, pour tout $\tilde{M}\in {\cal L}(\tilde{M}_{0})$, on d\'efinit $H_{\tilde{M}}:\tilde{M}({\mathbb R})\to {\cal A}_{\tilde{M}}$ comme l'unique application telle que 

- $H_{\tilde{M}}(m\gamma)=H_{\tilde{M}}(m)+H_{\tilde{M}}(\gamma)$ pour tous $m\in M({\mathbb R})$ et $\gamma\in \tilde{M}({\mathbb R})$;

- la restriction de $H_{\tilde{M}}$ \`a $\tilde{M}_{0}({\mathbb R})$ est la compos\'ee de $H_{\tilde{M}_{0}}$ et de la projection ${\cal A}_{\tilde{M}_{0}}\to {\cal A}_{\tilde{M}}$.

\bigskip

La th\'eorie est vide si $\omega$ n'est pas trivial sur $Z(G)^{\theta}({\mathbb R})$. Nous ne supposons toutefois pas que $\omega$ est  trivial sur ce groupe car l'inconv\'enient de cette hypoth\`ese est qu'elle ne se conserve pas si l'on remplace $\tilde{G}$ par un espace de  Levi $\tilde{M}$.  Par contre, nous supposerons que $\omega$ est trivial sur la composante neutre de $Z(G)^{\theta}({\mathbb R})$ pour la topologie r\'eelle. 
Cette hypoth\`ese se conserve si l'on remplace $\tilde{G}$ par un espace de Levi $\tilde{M}$.

  \bigskip
  
  \subsection{Rappels sur les $\omega$-repr\'esentations}
  Rappelons qu'une $\omega$-repr\'esentation (admissible) de $\tilde{G}({\mathbb R})$ est un couple $(\pi,\tilde{\pi})$, o\`u $\pi$ est une repr\'esentation admissible de $G({\mathbb R})$ dans un espace complexe $V$ et $\tilde{\pi}$ est une application de $\tilde{G}({\mathbb R})$ dans le groupe des automorphismes de $V$ qui v\'erifie la condition $\tilde{\pi}(g\gamma g')=\pi(g)\tilde{\pi}(\gamma)\pi(g')\omega(g')$ pour tous $g,g'\in G({\mathbb R})$ et $\gamma\in \tilde{G}({\mathbb R})$.  On supposera toujours $\pi$ de longueur finie.  En pratique, on notera simplement $\tilde{\pi}$ la $\omega$-repr\'esentation, la repr\'esentation $\pi$ de $G({\mathbb R})$ \'etant sous-entendue.
  
  A une telle $\omega$-repr\'esentation est associ\'e son caract\`ere, qui est une forme lin\'eaire sur $C_{c}^{\infty}(\tilde{G}({\mathbb R}))\otimes Mes(G({\mathbb R}))$, que l'on note
  $${\bf f}\mapsto I^{\tilde{G}}(\tilde{\pi},{\bf f}).$$
  Pr\'ecis\'ement, pour une fonction $f\in C_{c}^{\infty}(\tilde{G}({\mathbb R}))$ et une mesure de Haar $dg$ sur $G({\mathbb R})$, on d\'efinit l'op\'erateur $\tilde{\pi}(f\otimes dg)=\int_{\tilde{G}({\mathbb R})}\tilde{\pi}(\gamma)f(\gamma)\,d\gamma$ (la mesure $d\gamma$ \'etant naturellement associ\'ee \`a $dg$); alors 
$$ I^{\tilde{G}}(\tilde{\pi},f\otimes dg)=trace(\tilde{\pi}(f\otimes dg)).$$ 
Cette distribution est continue quand on munit $C_{c}^{\infty}(\tilde{G}({\mathbb R}))$ de sa topologie usuelle. Elle se factorise en une forme lin\'eaire continue sur $I(\tilde{G}({\mathbb R}),\omega)\otimes Mes(G({\mathbb R}))$. Ici, $I(\tilde{G}({\mathbb R}),\omega)$ d\'esigne l'espace des $\omega$-int\'egrales orbitales sur $\tilde{G}({\mathbb R})$. Il est muni d'apr\`es Bouaziz d'une topologie (cf. [I] 5.2) pour laquelle l'application naturelle $C_{c}^{\infty}(\tilde{G}({\mathbb R}))\to I(\tilde{G}({\mathbb R}),\omega)$ est continue et ouverte ([R] th\'eor\`eme 9.4). On note $D_{spec}(\tilde{G}({\mathbb R}),\omega)$ l'espace engendr\'e par les distributions $I^{\tilde{G}}(\tilde{\pi},.)$, quand  $\tilde{\pi}$ d\'ecrit les $\omega$-repr\'esentations de longueur finie, tensoris\'e par $Mes(G({\mathbb R}))$. Ainsi, ces distributions appartiennent \`a $D_{spec}(\tilde{G}({\mathbb R}),\omega)\otimes Mes(G({\mathbb R}))^*$. On note $D_{temp}(\tilde{G}({\mathbb R}),\omega)$ le sous-espace de $D_{spec}(\tilde{G}({\mathbb R}),\omega)$ qui, apr\`es tensorisation par $Mes(G({\mathbb R}))^*$, est engendr\'e par les caract\`eres de repr\'esentations temp\'er\'ees.

  Pour $\lambda\in {\cal A}_{\tilde{G},{\mathbb C}}^*$ et pour une $\omega$-repr\'esentation $(\pi,\tilde{\pi})$, on d\'efinit $(\pi_{\lambda},\tilde{\pi}_{\lambda})$ par $\pi_{\lambda}(g)=e^{<\lambda,H_{\tilde{G}}(g)>}\pi(g)$ et $\tilde{\pi}_{\lambda}(\gamma)=e^{<\lambda,H_{\tilde{G}}(\gamma)>}\tilde{\pi}(\gamma)$. L'action ainsi obtenue de ${\cal A}_{\tilde{G},{\mathbb C}}^*$ sur l'ensemble des $\omega$-repr\'esentations (\`a isomorphisme pr\`es) est libre.

Notons $\mathfrak{Z}(G)$ le centre de l'alg\`ebre enveloppante de l'alg\`ebre de Lie de $G({\mathbb R})$.   On consid\`ere qu'il  agit sur les espaces de fonctions sur $G({\mathbb R})$ via l'action par translations \`a gauche de $G({\mathbb R})$, et qu'il agit sur les espaces de distributions par dualit\'e, c'est-\`a-dire par  la formule $(ZD)(f)=D(Zf)$ pour une distribution $D$, une fonction $f$ et un \'el\'ement $ Z\in \mathfrak{Z}(G)$. Comme on sait, $\mathfrak{Z}(G)$ est isomorphe \`a $Sym(\mathfrak{h})^W$, ou encore \`a l'alg\`ebre des polyn\^omes invariants par $W$ sur $\mathfrak{h}^*$. Notons $\mathfrak{h}_{Z}$ la partie centrale de $\mathfrak{h}$, c'est-\`a-dire $\mathfrak{h}_{Z}=X_{*}(Z(G)^0)\otimes_{{\mathbb Z}}{\mathbb C}$.  On a dualement un sous-espace $\mathfrak{h}_{Z}^*\subset \mathfrak{h}^*$. Au caract\`ere $\omega$ est associ\'e un caract\`ere infinit\'esimal de $\mathfrak{Z}(G)$, qui est l'\'evaluation en un point $\mu(\omega)\in \mathfrak{h}_{Z}^*$.  L'hypoth\`ese que $\omega$ est trivial sur la composante neutre de $Z(G;{\mathbb R})^{\theta}$ pour la topologie r\'eelle implique que $\mu(\omega)$ appartient au sous-espace $(1-\theta)(\mathfrak{h}_{Z}^*)$. Il existe un unique point $\tilde{\mu}(\omega)$ dans ce sous-espace tel que $ \mu(\omega)=(\theta^{-1}-1)(\tilde{\mu}(\omega))$. Consid\'erons l'espace affine $\tilde{\mu}(\omega)+\mathfrak{h}^{\theta,*}$.   Il est invariant par l'action de $W^{\theta}$. Notons $Pol(\tilde{\mu}(\omega)+\mathfrak{h}^{\theta,*})$ l'alg\`ebre des polyn\^omes sur cet espace affine. Par restriction, tout \'el\'ement de $\mathfrak{Z}(G)$ d\'efinit  un \'el\'ement de $Pol(\tilde{\mu}(\omega)+\mathfrak{h}^{\theta,*})^{W^{\theta}}$. On a

(1) l'application $\mathfrak{Z}(G)\to Pol(\tilde{\mu}(\omega)+\mathfrak{h}^{\theta,*})^{W^{\theta}}$ est un homomorphisme d'alg\`ebres surjectif.

Preuve. Parce que $\tilde{\mu}(\omega)$ est central, il existe un unique automorphisme $\iota'$ de $\mathfrak{Z}(G)$ qui, \`a un \'el\'ement $X\in \mathfrak{h}\subset \mathfrak{Z}(G)$, associe l'\'el\'ement $X+<X,\tilde{\mu}(\omega)>$. Pour des \'el\'ements $Z\in \mathfrak{Z}(G)$ et  $\lambda\in \mathfrak{h}^{\theta,*}$, \'evaluer $Z$ en $\tilde{\mu}+\lambda$ revient \`a \'evaluer $\iota'(Z)$ en $\lambda$. Cela ram\`ene l'assertion au cas o\`u $\tilde{\mu}(\omega)=0$. C'est alors le th\'eor\`eme 5 de [DM]. $\square$

Soit $(\pi,\tilde{\pi})$ une $\omega$-repr\'esentation, supposons $\pi$ irr\'eductible. A $\pi$ est associ\'e son caract\`ere infinit\'esimal, qui est param\'etr\'e par une orbite dans $\mathfrak{h}^*$ pour l'action de $W$. On note $\mu(\pi)$ ou  $\mu(\tilde{\pi})$ cette orbite. Parce que $\pi$ se prolonge en une $\omega$-repr\'esentation, on a l'\'egalit\'e $\mu(\pi)+\mu(\omega)=\theta^{-1}(\mu(\pi))$. L'ensemble $\mu(\pi)-\tilde{\mu}(\omega)$ est alors une $W$-orbite qui est invariante par $\theta$. On v\'erifie que l'intersection d'une telle orbite avec $\mathfrak{h}^{\theta,*}$  est non vide et est une seule orbite sous l'action de $W^{\theta}$. Autrement dit,  l'ensemble
$$(\tilde{\mu}(\omega)+\mathfrak{h}^{\theta,*}) \cap \mu(\tilde{\pi}).$$
 est une unique orbite sous l'action de $W^{\theta}$ dans l'espace affine  $\tilde{\mu}(\omega)+\mathfrak{h}^{\theta,*}$.

Si $(\pi,\tilde{\pi})$ est une $\omega$-repr\'esentation de longueur finie et si toutes les composantes irr\'eductibles de $\pi$ ont un m\^eme param\`etre $\mu$, on pose $\mu(\tilde{\pi})=\mu$.  Plus g\'en\'eralement, on dira qu'un \'el\'ement de $D_{spec}(\tilde{G}({\mathbb R}),\omega)$ est de param\`etre $\mu$ si, modulo le choix d'une mesure de Haar, c'est une combinaison lin\'eaire de caract\`eres de repr\'esentations irr\'eductibles $\tilde{\pi}$ dont le param\`etre est $\mu$.  On note $D_{spec,\mu}(\tilde{G}({\mathbb R}),\omega)$ le sous-espace des \'el\'ements de param\`etre $\mu$.
 
On a d\'efini en [W] 6.2 l'ensemble ${\cal E}_{ell}(\tilde{G},\omega)$. Il est form\'e de triplets $\boldsymbol{\tau}=(M,\sigma,\boldsymbol{\tilde{r}})$, o\`u $M$ est un Levi de $G$ contenant $M_{0}$, $\sigma$ est une repr\'esentation irr\'eductible de la s\'erie discr\`ete de $M({\mathbb R})$ et $\boldsymbol{\tilde{r}}$ est un \'el\'ement de l'ensemble ${\cal R}^{\tilde{G}}(\sigma)$ d\'efini en [W] 2.8. Ces \'el\'ements sont soumis \`a des conditions telles qu'\`a $\boldsymbol{\tau}$ est associ\'ee une repr\'esentation "elliptique" $\tilde{\pi}_{\boldsymbol{\tau}}$ de $\tilde{G}({\mathbb R})$. Deux triplets peuvent \^etre conjugu\'es par $G({\mathbb R})$ et donnent dans ce cas la m\^eme repr\'esentation de $\tilde{G}({\mathbb R})$.  D'autre part, le groupe ${\mathbb U}=\{z\in {\mathbb C}; \vert z\vert =1\}$ agit naturellement sur l'ensemble  ${\cal R}^{\tilde{G}}(\sigma)$, donc sur ${\cal E}_{ell}(\tilde{G},\omega)$ par $z(M,\sigma,\boldsymbol{\tilde{r}})=(M,\sigma,z\boldsymbol{\tilde{r}})$. On a $\tilde{\pi}_{z\boldsymbol{\tau}}=z\tilde{\pi}_{\boldsymbol{\tau}}$ (on rappelle que, dans notre situation tordue, les "repr\'esentations" peuvent \^etre multipli\'ees par un nombre complexe).    Pour un couple $(M,\sigma)$ comme ci-dessus et pour $\lambda\in {\cal A}_{\tilde{G},{\mathbb C}}^*$, on d\'efinit la repr\'esentation $\sigma_{\lambda}$. La restriction \`a $i{\cal A}_{\tilde{G}}^*$ de cette action de ${\cal A}_{\tilde{G},{\mathbb C}}^*$ s'\'etend en une action $\boldsymbol{\tau}\mapsto \boldsymbol{\tau}_{\lambda}$ de $i{\cal A}_{\tilde{G}}^*$ sur ${\cal E}_{ell}(\tilde{G},\omega)$. On a $\tilde{\pi}_{\boldsymbol{\tau}_{\lambda}}=(\tilde{\pi}_{\boldsymbol{\tau}})_{\lambda}$.

Pour $M$ et $\sigma$ comme ci-dessus, $\sigma$ poss\`ede un caract\`ere central $\chi_{\sigma}$ et un caract\`ere central infinit\'esimal qui est param\'etr\'e par une orbite $\mu(\sigma)$ du groupe de Weyl $W^M$ dans $\mathfrak{h}^*$, plus pr\'ecis\'ement dans $\mathfrak{h}_{\mathbb R}^{M,*}\oplus i\mathfrak{a}^*_{M}({\mathbb R})$.    Pour $\boldsymbol{\tau}=(M,\sigma,\boldsymbol{\tilde{r}})\in {\cal E}_{ell}(\tilde{G},\omega)$,    on note $\mu(\boldsymbol{\tau})$ la $W$-orbite engendr\'ee par $\mu(\sigma)$. On note ${\cal E}_{ell,0}(\tilde{G},\omega)$ le sous-ensemble des $\boldsymbol{\tau}=(M,\sigma,\boldsymbol{\tilde{r}})\in {\cal E}_{ell}(\tilde{G},\omega)$ tels que  $\chi_{\sigma}$ soit trivial sur $\mathfrak{A}_{\tilde{G}}$.  Tout \'el\'ement de ${\cal E}_{ell}(\tilde{G},\omega)$ s'\'ecrit de fa\c{c}on unique $\boldsymbol{\tau}_{\lambda}$ pour un couple $(\boldsymbol{\tau},\lambda)\in {\cal E}_{ell,0}(\tilde{G},\omega)\times i{\cal A}_{\tilde{G}}^*$. On note $D_{ell}(\tilde{G}({\mathbb R}),\omega)$, resp. $D_{ell,0}(\tilde{G}({\mathbb R}),\omega)$, le sous-espace de $D_{temp}(\tilde{G}({\mathbb R}),\omega)$ engendr\'e par les caract\`eres des repr\'esentations $\tilde{\pi}_{\boldsymbol{\tau}}$ pour $\boldsymbol{\tau}\in {\cal E}_{ell}(\tilde{G},\omega)$, resp. ${\cal E}_{ell,0}(\tilde{G},\omega)$. Pour toute $W$-orbite $\mu$ dans $\mathfrak{h}^*$, on d\'efinit les espaces $D_{ell,\mu}(\tilde{G}({\mathbb R}),\omega)$, resp. $D_{ell,0,\mu}(\tilde{G}({\mathbb R}),\omega)$, en se limitant aux $\boldsymbol{\tau}$ tels que $\mu(\boldsymbol{\tau})=\mu$. On a $D_{ell,0,\mu}(\tilde{G}({\mathbb R}),\omega)=D_{ell,\mu}(\tilde{G}({\mathbb R}),\omega)$ si la projection de $\mu$ dans $\mathfrak{a}_{\tilde{G}}^*({\mathbb C})$ est nulle. Sinon, $D_{ell,0,\mu}(\tilde{G}({\mathbb R}),\omega)=0$. 

On note $D_{ell,{\mathbb C}}(\tilde{G}({\mathbb R}),\omega)$ le sous-espace de $D_{spec}(\tilde{G}({\mathbb R}),\omega)$ engendr\'e par les caract\`eres de repr\'esentations $(\tilde{\pi}_{\boldsymbol{\tau}})_{\lambda}$ pour $\boldsymbol{\tau}\in {\cal E}_{ell}(\tilde{G},\omega)$ et $\lambda\in {\cal A}_{\tilde{G},{\mathbb C}}^*$. On a l'\'egalit\'e
$$(2) \qquad D_{spec}(\tilde{G}({\mathbb R}),\omega)=\left(\oplus_{\tilde{L}\in {\cal L}(\tilde{M}_{0})}Ind_{\tilde{L}}^{\tilde{G}}(D_{ell,{\mathbb C}}(\tilde{L}({\mathbb R}),\omega))\right)^{W(\tilde{M}_0)}.$$

\bigskip

\subsection{Espaces de Paley-Wiener}
On consid\`ere les donn\'ees suivantes:

-  $E$ est un ensemble;

- $D$ est un entier  positif ou nul;

-  $d:E\to {\mathbb R}_{\geq 0}$ est une fonction; 

- pour tout $e\in E$, $V_{e}$ est un espace vectoriel r\'eel de dimension inf\'erieure ou \'egale \`a $D$ muni d'une forme quadratique d\'efinie positive. 

L'espace dual $V_{e}^*$ est donc lui-aussi muni d'une telle forme. On prolonge ces formes en des formes hermitiennes sur $V_{e,{\mathbb C}}$ et $V_{e,{\mathbb C}}^*$. 

Pour tout r\'eel $r>0$, notons $PW^r$ l'espace des familles ${\bf f}=(f_{e})_{e\in E}$, o\`u, pour tout $e\in E$, $f_{e}$ est une fonction enti\`ere sur $V_{e,{\mathbb C}}^*$, qui v\'erifient la condition

(1) pour tout $N\in {\mathbb N}$, il existe $C_{N}>0$ tel que, pour tout $e\in E$ et tout $\lambda\in V_{e,{\mathbb C}}^*$, on ait l'in\'egalit\'e
$$\vert f_{e}(\lambda)\vert \leq C_{N}(1+d(e)+\vert \lambda\vert )^{-N}e^{r\vert Re(\lambda)\vert }.$$

On munit cet espace de la famille de semi-normes
$${\cal N}^r_{N}({\bf f})=sup_{e\in E, \lambda\in V_{e,{\mathbb C}}^*}(1+d(e)+\vert \lambda\vert )^{N}e^{-r\vert Re(\lambda)\vert }\vert f_{e}(\lambda)\vert $$
pour $N\in {\mathbb N}$. C'est un espace de Fr\'echet. Pour $r<r'$, l'injection $PW^r\to PW^{r'}$ est continue et on note $PW$ la limite inductive des $PW^r$, muni de la topologie limite inductive.   

Pour tout r\'eel $r>0$, notons $\underline{PW}^r$ l'espace des familles  ${\bf f}=(f_{e})_{e\in E}$, o\`u, pour tout $e\in E$, $f_{e}$ est une fonction enti\`ere sur $V_{e,{\mathbb C}}^*$, qui v\'erifient les conditions

(2) pour tout $N$ et tout $e\in E$, il existe $C_{N}(e)>0$ tel que, pour tout $\lambda\in V_{e,{\mathbb C}}^*$, on ait l'in\'egalit\'e
$$\vert f_{e}(\lambda)\vert \leq C_{N}(e)(1+\vert \lambda\vert )^{-N}e^{r\vert Re(\lambda)\vert };$$

(3) pour tout $N$, il existe $\underline{C}_{N}>0$ tel que, pour tout $e\in E$ et tout $\lambda\in iV_{e}^*$, on ait l'in\'egalit\'e
$$\vert f_{e}(\lambda)\vert \leq \underline{C}_{N}(1+d(e)+\vert \lambda\vert )^{-N}.$$

On munit $\underline{PW}^r$ de la famille de semi-normes
$$\underline{{\cal N}}_{N}({\bf f})=sup_{e\in E,\lambda\in iV_{e}^*}(1+d(e)+\vert \lambda\vert )^{N}\vert f_{e}(\lambda)\vert $$
pour $N\in {\mathbb N}$. De nouveau, on note $\underline{PW}$ la limite inductive des $\underline{PW}^r$ munie de la topologie limite inductive.

Il est clair que, pour tout $r$,  $PW^r$ est inclus dans $\underline{PW}^r$ et que cette injection est continue. D'o\`u une injection continue $PW\subset \underline{PW}$.

\ass{Lemme}{Cette application est bijective et c'est un hom\'eomorphisme.}

{\bf Remarque.} Ce lemme est \'el\'ementaire. On n'en donne une preuve que pour la commodit\'e du r\'edacteur.

Preuve. On peut d\'ecomposer $E$ en union finie disjointe de sous-ensembles sur lesquels la fonction $e\mapsto dim(V_{e})$ est constante. On voit qu'il suffit de d\'emontrer le lemme analogue obtenu en rempla\c{c}ant $E$ par un tel sous-ensemble. En oubliant cela, on peut supposer que l'espace $V_{e}$, muni de sa forme quadratique,  est constant et on l'identifie \`a un espace fixe $V$. Pour simplifier, on identifie $V$ \`a son dual $V^*$ \`a l'aide de la forme quadratique. On consid\`ere donc que $f_{e}$ est d\'efinie sur $V_{{\mathbb C}}$ pour tout $e$. 

 Fixons $\epsilon>0$. On va montrer que, pour tout $r>0$, $\underline{PW}^r$ est inclus dans $PW^{r+\epsilon}$ et que cette injection est continue. Pour cela, fixons $r>0$ et $N\in {\mathbb N}$. On va prouver plus pr\'ecis\'ement qu'il existe $c>0$ (d\'ependant de $N$, $r$ et $\epsilon$) tel que, pour tout ${\bf f}\in \underline{PW}^r$, on a l'in\'egalit\'e
$$(4) \qquad {\cal N}^{r+\epsilon}_{N}({\bf f})\leq c\underline{{\cal N}}_{4N+4D}({\bf f}).$$
Soit ${\bf f}=(f_{e})_{e\in E}\in \underline{PW}^r$.  Pour tout $e\in E$, on d\'efinit une fonction $\varphi_{e}$ sur $V$ par
$$(5) \qquad \varphi_{e}(x)=\int_{iV}f_{e}(\lambda)e^{-2\pi (x,\lambda)}\,d\lambda.$$
La condition (2)  entra\^{\i}ne que la fonction $\varphi_{e}$ est $C^{\infty}$ et, par un proc\'ed\'e usuel de d\'eplacement de contour,  que cette fonction est \`a support dans l'ensemble des $x\in V$ tels que $\vert x\vert \leq r/2\pi$. Par inversion de Fourier, on a
$$(6) \qquad f_{e}(\lambda)=\int_{V}\varphi_{e}(x)e^{2\pi(x,\lambda)}\,dx$$
pour $\lambda\in iV$ et cette \'egalit\'e persiste pour tout $\lambda\in V_{{\mathbb C}}$ par continuation holomorphe. On veut majorer l'expression
$$(1+d(e)+\vert \lambda\vert )^Ne^{-(r+\epsilon)\vert Re(\lambda)\vert }\vert f_{e}(\lambda)\vert $$
pour tout $e\in E$ et tout $\lambda\in V_{{\mathbb C}}$. On a tout d'abord
$$(1+d(e)+\vert \lambda\vert )^N\leq(1+d(e))^N(1+\vert\lambda\vert)^N$$
car $1+x+y\leq(1+x)(1+y)$ pour tous $x,y\geq0$. On a
 $$(1+\vert\lambda\vert)^N\leq 2^N(1+\vert \lambda\vert ^2)^N=2^N(1+\vert Im(\lambda)\vert ^2+\vert Re(\lambda)\vert ^2)^N$$
 car $1+x\leq 2(1+x^2)$ pour tout $x\geq0$. Introduisons des coordonn\'ees sur $V$ relatives \`a une base orthogonale. On a
 $$1+\vert Im(\lambda)\vert ^2+\vert Re(\lambda)\vert ^2=1+\sum_{j=1,...,D}(Im(\lambda_{j})^2+Re(\lambda_{j})^2)=1+\sum_{j=1,...,D}(-\lambda_{j}^2+2\lambda_{j}Re(\lambda_{j})).$$
 Ainsi $(1+\vert Im(\lambda)\vert ^2+\vert Re(\lambda)\vert ^2)^N$ s'exprime comme combinaison lin\'eaire finie de produits de mon\^omes de degr\'e au plus $2N$ en les $\lambda_{j}$ et de mon\^omes de degr\'e au plus $N$ en les $Re(\lambda_{j})$. On peut aussi bien fixer deux tels mon\^omes $P(\lambda)$ et $Q(Re(\lambda))$ et majorer
 $$(1+d(e))^N\vert P(\lambda)\vert \vert Q(Re(\lambda))\vert e^{-(r+\epsilon)\vert Re(\lambda)\vert }\vert f_{e}(\lambda)\vert .$$
Il existe des constantes $c_{1},c_{2}>0$ telles que 
$$\vert Q(Re(\lambda))\vert \leq c_{1}(1+\vert Re(\lambda)\vert)^N$$
et $(1+x)^Ne^{-\epsilon x}\leq c_{2}$ pour tout $x\geq0$. Ainsi l'expression pr\'ec\'edente est major\'ee par le produit d'une constante et de  
 $$(7) \qquad (1+d(e))^N\vert P(\lambda)\vert   e^{-r\vert Re(\lambda)\vert }\vert f_{e}(\lambda)\vert .$$
 Par les r\`egles usuelles de d\'erivation, on d\'eduit de (6) l'existence d'un op\'erateur diff\'erentiel \`a coefficients constants $\partial_{P}$ sur $V$, d'ordre au plus $2N$, tel que
 $$(8) \qquad P(\lambda)f_{e}(\lambda)= \int_{V}\partial_{P}\varphi_{e}(x)e^{2\pi(x,\lambda)}\,dx.$$
 Par les m\^emes r\`egles de d\'erivation, appliqu\'ees cette fois \`a l'expression (5), on a
 $$(9) \qquad \partial_{P}\varphi_{e}(x)=\int_{iV}P(\lambda)f_{e}(\lambda)e^{-2\pi (x,\lambda)}\,d\lambda.$$
 Il existe une constante $c_{3}>0$ telle que
 $$\vert P(\lambda)\vert \leq c_{3}(1+\vert \lambda\vert )^{2N}.$$
 Alors, pour $\lambda\in iV$, 
 $$\vert P(\lambda) f_{e}(\lambda)\vert \leq c_{3}\underline{{\cal N}}_{4N+4D}({\bf f})(1+\vert \lambda\vert )^{2N}(1+d(e)+\vert \lambda\vert )^{-4N-4D}.$$
 On utilise  que $(1+x+y)\geq (1+x)^{1/2}(1+y)^{1/2}$ pour tous $x,y\geq0$ et on obtient
  $$\vert P(\lambda) f_{e}(\lambda)\vert \leq c_{3}\underline{{\cal N}}_{4N+4}({\bf f})(1+\vert \lambda\vert )^{-2D}(1+d(e))^{-2N-2D}.$$
  La fonction $\lambda\mapsto (1+\vert \lambda\vert )^{-2D}$ est int\'egrable sur $iV$. Gr\^ace \`a (9), on en d\'eduit l'existence de $c_{4}>0$ tel que
  $$\partial_{P}\varphi_{e}(x)\leq c_{4}\underline{{\cal N}}_{4N+4D}({\bf f})(1+d(e))^{-2N-2D}$$
  pour tout $x\in V$. Gr\^ace \`a (8) et \`a la propri\'et\'e du support de $\varphi_{e}$, on obtient une constante $c_{5}>0$ telle que
  $$\vert P(\lambda)f_{e}(\lambda)\vert \leq c_{5}\underline{{\cal N}}_{4N+4D}({\bf f})(1+d(e))^{-2N-2D}e^{r\vert Re(\lambda)\vert }$$
  pour tout $\lambda\in V_{{\mathbb C}}$. Alors, l'expression (7) est born\'ee par le produit d'une constante et de $\underline{{\cal N}}_{4N+4D}({\bf f})$. Cela prouve la majoration (4) et le lemme. $\square$
  
  Ce lemme \'etant d\'emontr\'e, on n'aura plus besoin de distinguer $PW$ de $\underline{PW}$ et on ne conservera que la notation $PW$.

\bigskip

\subsection{Enonc\'e du th\'eor\`eme}
On d\'efinit une fonction $d:{\cal E}_{ell,0}(\tilde{G},\omega)\to {\mathbb R}_{\geq0}$ par $d(\boldsymbol{\tau})=\vert \mu(\boldsymbol{\tau})\vert $, o\`u l'on d\'esigne ainsi la norme pour la forme hermitienne fix\'ee sur $\mathfrak{h}$ d'un \'el\'ement quelconque de $\mu(\boldsymbol{\tau})$.   On d\'efinit l'espace $PW_{ell}^{\infty}(\tilde{G},\omega)$ comme celui des fonctions $\varphi:{\cal E}_{ell}(\tilde{G},\omega)\to {\mathbb C}$ qui v\'erifient les conditions suivantes:

(1) si deux \'el\'ements $\boldsymbol{\tau}$ et $\boldsymbol{\tau}'$ de ${\cal E}_{ell}(\tilde{G},\omega)$ sont conjugu\'es par $G({\mathbb R})$, alors $\varphi(\boldsymbol{\tau})=\varphi(\boldsymbol{\tau}')$;

(2) pour $\boldsymbol{\tau}\in {\cal E}_{ell}(\tilde{G},\omega)$ et $z\in {\mathbb U}$, on a $\varphi(z\boldsymbol{\tau})=z\varphi(\boldsymbol{\tau})$;

(3) pour $\boldsymbol{\tau}\in {\cal E}_{ell}(\tilde{G},\omega)$, la fonction $\lambda\mapsto \varphi(\boldsymbol{\tau}_{\lambda})$ sur $i{\cal A}_{\tilde{G}}^*$ s'\'etend en une fonction enti\`ere $\varphi_{\boldsymbol{\tau}}$ sur ${\cal A}_{\tilde{G},{\mathbb C}}^*$;

(4) fixons un ensemble de repr\'esentants $\underline{{\cal E}}_{ell,0}(\tilde{G},\omega)$ de repr\'esentants des classes d'\'equivalence dans ${\cal E}_{ell,0}(\tilde{G},\omega)$ pour l'\'equivalence engendr\'ee par la conjugaison par $G({\mathbb R})$ et par l'action de ${\mathbb U}$; alors la famille $(\varphi_{\boldsymbol{\tau}})_{\boldsymbol{\tau}\in \underline{{\cal E}}_{ell,0}(\tilde{G},\omega)}$ appartient \`a l'espace de Paley-Wiener d\'efini en 1.3 relatif \`a l'ensemble $\underline{{\cal E}}_{ell,0}(\tilde{G},\omega)$ muni de la fonction $d$. 

Le groupe $W(\tilde{M}_0) $ agit naturellement dans 
$$\oplus_{\tilde{L}\in {\cal L}(\tilde{M}_{0})}PW_{ell}^{\infty}(\tilde{L},\omega).$$
On note $PW^{\infty}(\tilde{G},\omega)$ le sous-espace des invariants. 

{\bf Remarque.} On peut se limiter aux $\tilde{L}$ tels que la restriction de $\omega$ \`a $Z(L;{\mathbb R})^{\theta}$ est triviale. Pour les autres, les espaces correspondants sont nuls pour la simple raison que les ensembles ${\cal E}_{ell}(\tilde{L},\omega)$ sont vides: il n'y a pas de $\omega$-repr\'esentations de $\tilde{L}({\mathbb R})$.

  \bigskip
Soit ${\bf f}\in C_{c}^{\infty}(\tilde{G}({\mathbb R}))\otimes Mes(G({\mathbb R}))$.  Pour $\tilde{L}\in {\cal L}(\tilde{M}_{0})$ et $\boldsymbol{\tau}\in {\cal E}_{ell}(\tilde{L},\omega)$, posons $\varphi_{{\bf f}}(\boldsymbol{\tau})=I^{\tilde{L}}(\tilde{\pi}_{\boldsymbol{\tau}},{\bf f}_{\tilde{L},\omega})=trace(\tilde{\pi}_{\boldsymbol{\tau}}({\bf f}_{\tilde{L},\omega}))$. On peut aussi dire que 
$\varphi_{{\bf f}}(\boldsymbol{\tau})=I^{\tilde{G}}(Ind_{\tilde{Q}}^{\tilde{G}}(\tilde{\pi}_{\boldsymbol{\tau}}),{\bf f})=trace(Ind_{\tilde{Q}}^{\tilde{G}}(\tilde{\pi}_{\boldsymbol{\tau}})({\bf f}))$, o\`u $\tilde{Q}$ est un espace parabolique quelconque de composante de Levi $\tilde{L}$. On a ainsi d\'efini une application lin\'eaire qui, \`a ${\bf f}\in C_{c}^{\infty}(\tilde{G}({\mathbb R}))\otimes Mes(G({\mathbb R}))$, associe une famille de fonctions sur $\sqcup_{\tilde{L}\in {\cal L}(\tilde{M}_{0})}{\cal E}_{ell}(\tilde{L},\omega)$. Elle se quotiente en une application lin\'eaire sur $I(\tilde{G}({\mathbb R}),\omega)\otimes Mes(G({\mathbb R}))$.

\ass{Th\'eor\`eme}{Cette application est un hom\'eomorphisme de $I(\tilde{G}({\mathbb R}),\omega)\otimes Mes(G({\mathbb R}))$ sur $PW^{\infty}(\tilde{G},\omega)$.}

Ce th\'eor\`eme est prouv\'e par Renard ([R] th\'eor\`eme 17.5). La formulation de Renard \'etant largement diff\'erente de la n\^otre, nous montrerons dans les  deux paragraphes suivants pourquoi l'\'enonc\'e de Renard implique l'\'enonc\'e ci-dessus.

Rappelons pour m\'emoire le r\'esultat de Delorme et Mezo. On note $C_{c}^{\infty}(\tilde{G}({\mathbb R}),K)$ l'espace des \'el\'ements de $C_{c}^{\infty}(\tilde{G}({\mathbb R}))$ qui sont $K$-finis \`a droite et \`a gauche. On note $I(\tilde{G}({\mathbb R}),\omega,K)$ son image dans $I(\tilde{G}({\mathbb R}),\omega)$.  Notons $PW_{ell}(\tilde{G},\omega)$ le sous-espace des familles $(\varphi_{\boldsymbol{\tau}})_{\boldsymbol{\tau}\in \underline{{\cal E}}_{ell,0}(\tilde{G},\omega)}\in PW^{\infty}_{ell}(\tilde{G},\omega)$ telles que $\varphi_{\boldsymbol{\tau}}=0$ pour presque tout $\boldsymbol{\tau}$. En supprimant les exposants $^{\infty}$, on d\'efinit comme ci-dessus l'espace $PW(\tilde{G},\omega)$, qui est un sous-espace de $PW^{\infty}(\tilde{G},\omega)$. Alors l'application du th\'eor\`eme se restreint en un isomorphisme de $I(\tilde{G}({\mathbb R}),\omega,K)$ sur $PW(\tilde{G},\omega)$. 

\bigskip

\subsection{La transition entre le th\'eor\`eme de Renard et le th\'eor\`eme 1.4}
On suppose $\omega=1$. On oublie les questions de mesures en fixant des mesures de Haar sur tous les groupes intervenant. Consid\'erons un sous-ensemble $\Omega$ de $\tilde{G}({\mathbb R})$ qui est r\'eunion de composantes connexes pour la topologie r\'eelle et qui est engendr\'e par conjugaison sous $G({\mathbb R})$ par une seule telle composante. Il peut n'exister aucun sous-tore tordu maximal elliptique $\tilde{T}$ de $\tilde{G}$ tel que $\tilde{T}({\mathbb R})\cap \Omega$ soit non vide. Supposons qu'il en existe un. Alors il n'en existe qu'un, \`a conjugaison pr\`es par $G({\mathbb R})$ ([R] lemme 12.12). Fixons un tel tore tordu $\tilde{T}$. Notons $\mathfrak{t}^{\theta,\tilde{G}}$ l'orthogonal de $\mathfrak{a}_{\tilde{G}}$ dans $\mathfrak{t}^{\theta}$. Renard introduit un certain sous-ensemble de $i\mathfrak{t}^{\theta,\tilde{G}}({\mathbb R})$, notons-le $H^*(\Omega)$. 

{\bf Remarque.} Plus pr\'ecis\'ement, Renard introduit un tel sous-ensemble sur lequel agit un certain groupe de Weyl. Les constructions de Renard pour deux \'el\'ements conjugu\'es par ce groupe sont essentiellement les m\^emes. Aussi, nous prendrons pour $H^*(\Omega)$ un ensemble de repr\'esentants des orbites dans l'ensemble de Renard pour l'action de ce groupe.

Pour $h^*\in H^*(\Omega)$, il d\'efinit une distribution $\Theta_{\Omega,h^*}$ sur $\Omega$ qui  v\'erifie de nombreuses propri\'et\'es. C'est une distribution propre pour l'action de $\mathfrak{Z}(G)$, donc il lui est associ\'ee une $W$-orbite $\mu(h^*)$ dans $\mathfrak{h}^*$.
D'apr\`es [R] paragraphe 18, c'est la restriction \`a $\Omega$ d'un \'el\'ement de $D_{temp}(\tilde{G}({\mathbb R}))$.  On peut donc tensoriser la distribution $\Theta_{\Omega,h^*}$ par un \'el\'ement $\lambda\in {\cal A}_{\tilde{G}}^*$: on note $\Theta_{\Omega,h^*,\lambda}$ la distribution obtenue. Notons $E_{ell}^{\tilde{G}}$ l'ensemble des paires $(\Omega,h^*)$, o\`u $\Omega$ parcourt les sous-ensembles de $\tilde{G}({\mathbb R})$ v\'erifiant les conditions ci-dessus   et $h^*$ parcourt $H^*(\Omega)$.  Remarquons  qu'il n'y  a qu'un nombre fini de $\Omega$ puisqu'il n'y a qu'un nombre fini de composantes connexes pour la topologie r\'eelle. On d\'efinit une fonction $d$ sur $E_{ell}^{\tilde{G}}$: pour  $(\Omega,h^*)\in E^{\tilde{G}}_{ell}$, $d(\Omega,h^*)=\vert \mu(h^*)\vert $. On introduit l'espace de Paley-Wiener associ\'e \`a cet ensemble $E_{ell}^{\tilde{G}}$ et \`a cette fonction $d$. Notons-le $PW_{ell}'(\tilde{G})$. De nouveau, le groupe $W(\tilde{M}_0)$ agit naturellement sur
$$\oplus_{\tilde{L}\in {\cal L}(\tilde{M}_{0})}PW_{ell}'(\tilde{L}).$$
On note $PW'(\tilde{G})$ le sous-espace des invariants. On pose $E=\cup_{\tilde{L}\in {\cal L}(\tilde{M}_{0})}E_{ell}^{\tilde{L}}$.

 Pour $f\in C_{c}^{\infty}(\tilde{G}({\mathbb R}))$ pour $\tilde{L}\in {\cal L}(\tilde{M}_{0})$ et pour $(\Omega,h^*)\in E^{\tilde{L}}$, on d\'efinit une fonction $\varphi_{\Omega,h^*}$ sur ${\cal A}_{\tilde{L},{\mathbb C}}^*$ par $\varphi_{\Omega,h^*}(\lambda)=\Theta_{\Omega,h^*,\lambda}(f_{\tilde{L}})$. Renard d\'emontre que l'application qui, \`a $f$, associe la famille de fonctions $(\varphi_{\Omega,h^*})_{(\Omega,h^*)\in E}$ se quotiente en un hom\'eomorphisme de $I(\tilde{G}({\mathbb R}))$ sur $PW'(\tilde{G})$.
 
 D'apr\`es [BT] corollaire 14.5, le groupe $\Pi$ des composantes connexes de $G({\mathbb R})$ pour la topologie r\'eelle est un groupe ab\'elien fini, en fait une puissance de ${\mathbb Z}/2{\mathbb Z}$. Il est muni d'un automorphisme $\theta$, \'egal \`a l'automorphisme d\'eduit de $ad_{\gamma}$ pour n'importe quel $\gamma\in \tilde{G}({\mathbb R})$. Notons $\Omega_{G}$ la r\'eunion des composantes appartenant au sous-groupe $(1-\theta)(\Pi)$. Alors tout ensemble $\Omega$ intervenant ci-dessus est une unique classe modulo $\Omega_{G}$, \`a droite ou \`a gauche. Notons $\Xi$ le groupe des caract\`eres de $G({\mathbb R})$ triviaux sur $\Omega_{G}$. Pour $\Omega$ intervenant ci-dessus et pour $\xi\in \Xi$, notons $\tilde{\xi}_{\Omega}$ l'unique fonction sur $\tilde{G}({\mathbb R})$ telle que $\tilde{\xi}_{\Omega}(g\gamma)=\xi(g)$ pour tous $g\in G({\mathbb R})$ et $\gamma\in \Omega$. La fonction caract\'eristique de $\Omega$ est \'egale \`a
 $$\vert \Xi\vert ^{-1}\sum_{\xi\in \Xi}\tilde{\xi}_{\Omega}.$$
 Soit $(\Omega,h^*)\in E_{ell}^{\tilde{G}}$. Comme on l'a dit, la distribution $\Theta_{\Omega,h^*}$ est la restriction \`a $\Omega$ d'un \'el\'ement de $D_{temp}(\tilde{G}({\mathbb R}))$, que l'on peut \'ecrire $\sum_{i\in I}c_{i}\tilde{\pi}_{i}$, o\`u $I$ est un ensemble fini d'indices et, pour tout $i\in I$, $c_{i}$ est un coefficient complexe et $\tilde{\pi}_{i}$ est une repr\'esentation irr\'eductible et temp\'er\'ee. On a alors l'\'egalit\'e
 $$\Theta_{\Omega,h^*}=\vert \Xi\vert ^{-1}\sum_{\xi\in \Xi}\sum_{i\in I}c_{i}\tilde{\xi}_{\Omega}\tilde{\pi}_{i},$$
 o\`u le produit $\tilde{\xi}_{\Omega}\tilde{\pi}_{i}$ se d\'efinit de fa\c{c}on \'evidente. Il est clair que $\tilde{\xi}_{\Omega}\tilde{\pi}_{i}$ est encore une repr\'esentation temp\'er\'ee et irr\'eductible.  Cela nous d\'ebarrasse du passage par la restriction \`a $\Omega$: la distribution $\Theta_{\Omega,h^*}$ s'identifie \`a un \'el\'ement de $D_{temp}(\tilde{G}({\mathbb R}))$. On la note d\'esormais  $f\mapsto I^{\tilde{G}}(\tilde{\pi}_{\Omega,h^*},f)$ pour un certain \'el\'ement $\tilde{\pi}_{\Omega,h^*}$ de $D_{temp}(\tilde{G}({\mathbb R}))$. 
 
 Soit $(\Omega,h^*)\in E^{\tilde{G}}_{ell}$ et soit comme ci-dessus $\tilde{T}$ un tore tordu maximal elliptique tel que $\tilde{T}({\mathbb R})$ coupe $\Omega$. Alors Renard calcule le caract\`ere de $\tilde{\pi}_{\Omega,h^*}$ sur $\tilde{T}({\mathbb R})$ ([R] 15.17).  Fixons $\gamma\in \tilde{T}({\mathbb R})\cap \Omega$. Tout \'el\'ement de cette intersection est conjugu\'e \`a un \'el\'ement de $exp(\mathfrak{t}^{\theta}({\mathbb R}))\gamma$. En un point $exp(X)\gamma$, avec $X\in \mathfrak{t}^{\theta}({\mathbb R})$, c'est le produit d'un module explicite avec  une combinaison lin\'eaire finie explicite de $e^{<w(h^*),X>}$, pour certains $w\in W$. Il en r\'esulte d'abord que 
 
 (1) $\tilde{\pi}_{\Omega,h^*}$  admet un caract\`ere central pour l'action de $\mathfrak{A}_{\tilde{G}}$  et que celui-ci est trivial. 
 
 Rappelons que l'on peut d\'efinir le produit scalaire elliptique $(\tilde{\pi}_{1},\tilde{\pi}_{2})_{ell}$ de deux repr\'esentations $\tilde{\pi}_{1}$ et $\tilde{\pi}_{2}$ se tranformant trivialement sous l'action de $\mathfrak{A}_{\tilde{G}}$, cf. [W] 7.3. Montrons que 
 
 (2) pour deux \'el\'ements distincts $(\Omega_{1},h^*_{1})$ et $(\Omega_{2},h^*_{2})\in E^{\tilde{G}}_{ell}$, on a l'\'egalit\'e $(\tilde{\pi}_{\Omega_{1},h^*_{1}},\tilde{\pi}_{\Omega_{2},h^*_{2}})_{ell}=0$;
 
 (3) quand $(\Omega,h^*)$ parcourt $E^{\tilde{G}}_{ell}$, les produits $(\tilde{\pi}_{\Omega,h^*},\tilde{\pi}_{\Omega,h^*})$ ne prennent qu'un nombre fini de valeurs r\'eelles strictement positives.
 
 Preuve de (2). Si $\Omega_{1}$ est distinct de $\Omega_{2}$, les caract\`eres sont de supports disjoints. Si $\Omega_{1}=\Omega_{2}$, le produit elliptique est, \`a une constante pr\`es, l'int\'egrale sur un unique  $\tilde{T}({\mathbb R})/(1-\theta)(T({\mathbb R}))$ du produit d'un certain module, du caract\`ere de $\tilde{\pi}_{\Omega_{2},h^*_{2}}$ et du conjugu\'e du  caract\`ere de $\tilde{\pi}_{\Omega_{1},h^*_{1}}$. D'apr\`es le r\'esultat \'evoqu\'e ci-dessus, c'est une combinaison lin\'eaire finie d'int\'egrales sur $exp(\mathfrak{t}^{\theta}({\mathbb R}))$ de fonctions du type $exp(X)\mapsto e^{<w_{2}(h^*_{2})-w_{1}(h^*_{1}),X>}$. L'int\'egrale d'une telle fonction n'est non nulle que si $w_{2}(h^*_{2})=w_{1}(h^*_{1})$. Une telle \'egalit\'e (o\`u ne peuvent intervenir que certains \'el\'ements du groupe de Weyl) entra\^{\i}ne que $h^*_{1}=h^*_{2}$ (pour $h^*_{1}$ et $h^*_{2}$ dans notre ensemble $H^*(\Omega_{1})$, cf. remarque ci-dessus). 
 
 Preuve de (3). Le m\^eme raisonnement montre que, si le caract\`ere de $\tilde{\pi}_{\Omega,h^*}$ en un point $exp(X)\gamma$ s'\'ecrit comme produit d'un module explicite et de
 $\sum_{w}c_{w}e^{<w(h^*),X>}$, alors
  $$(\tilde{\pi}_{\Omega,h^*},\tilde{\pi}_{\Omega,h^*})=m\sum_{w}\vert c_{w}\vert ^2,$$
  o\`u $m$ est une constante ne d\'ependant que des mesures. Or les constantes $c_{w}$ intervenant en [R] 15.17 sont de modules $1$. La somme ci-dessus est donc $m$ fois le nombre d'\'el\'ements de l'ensemble de sommation. Ce dernier \'etant un sous-ensemble de $W$, l'assertion s'ensuit. 
  \bigskip
  
  Renard d\'emontre une propri\'et\'e suppl\'ementaire de ses distributions. Le th\'eor\`eme [R] 15.18 signifie que
  
  (4) pour tout $(\Omega,h^*)\in E^{\tilde{G}}_{ell}$, $\tilde{\pi}_{\Omega,h^*}$ est supertemp\'er\'ee. 
  
  On renvoie \`a [Moe] pour cette notion. Or, dans [Moe], Moeglin d\'emontre qu'un \'el\'ement de $D_{temp}(\tilde{G}({\mathbb R}))$ qui est supertemp\'er\'e est combinaison lin\'eaire de caract\`eres de repr\'esentations elliptiques.   En joignant ce dernier r\'esultat \`a (1), on obtient que $\tilde{\pi}_{\Omega,h^*}\in D_{ell,0}(\tilde{G}({\mathbb R}))$ pour tout $(\Omega,h^*)\in E^{\tilde{G}}_{ell}$. 
  
 Pour tout $(\Omega,h^*)\in E^{\tilde{G}}_{ell}$, on peut donc exprimer $\tilde{\pi}_{\Omega,h^*}$ comme combinaison lin\'eaire des $\tilde{\pi}_{\boldsymbol{\tau}}$ pour $\boldsymbol{\tau}\in \underline{{\cal E}}_{ell,0}(\tilde{G})$ (rappelons que $\underline{{\cal E}}_{ell,0}(\tilde{G})$ est un ensemble de repr\'esentants fix\'e en 1.4). Plus pr\'ecis\'ement, pour toute $W$-orbite dans $\mathfrak{h}^*$, notons $E^{\tilde{G}}_{ell,\mu}$, resp. $\underline{{\cal E}}_{ell,0,\mu}(\tilde{G})$, le sous-ensemble des $(\Omega,h^*)\in E^{\tilde{G}}_{ell}$, resp. des $\boldsymbol{\tau}\in\underline{{\cal E}}_{ell,0}(\tilde{G})$, tels que le caract\`ere infinit\'esimal de $\tilde{\pi}_{\Omega,h^*}$, resp. $\tilde{\pi}_{\boldsymbol{\tau}}$, soit de param\`etre $\mu$. Comme on l'a dit, $E^{\tilde{G}}_{ell}$ est r\'eunion des $E^{\tilde{G}}_{ell,\mu}$ et, de m\^eme, $\underline{{\cal E}}_{ell,0}(\tilde{G})$ est r\'eunion des $\underline{{\cal E}}_{ell,0,\mu}(\tilde{G})$. L'\'ecriture d'un \'el\'ement de $D_{ell,0}(\tilde{G}({\mathbb R}))$ dans la base $\underline{{\cal E}}_{ell,0}(\tilde{G})$ est compatible avec l'action de $\mathfrak{Z}(G)$. Donc, pour tout $\mu$ et tout $(\Omega,h^*)\in E^{\tilde{G}}_{ell,\mu}$, $\tilde{\pi}_{\Omega,h^*}$ est combinaison lin\'eaire des 
 $\tilde{\pi}_{\boldsymbol{\tau}}$ pour $\boldsymbol{\tau}\in \underline{{\cal E}}_{ell,0,\mu}(\tilde{G})$. Notons que, d'apr\`es sa construction, l'ensemble $\underline{{\cal E}}_{ell,0,\mu}(\tilde{G})$ est fini et son nombre d'\'el\'ements est born\'e ind\'ependamment de $\mu$. Montrons que
 
 (5) la famille $(\tilde{\pi}_{\Omega,h^*})_{ (\Omega,h^*)\in E^{\tilde{G}}_{ell}}$ est une base de $D_{ell,0}(\tilde{G}({\mathbb R}))$.
 
Il est clair par (2) et (3) qu'elle est libre. On peut fixer $\mu$ et prouver que l'espace engendr\'e par
$(\tilde{\pi}_{\Omega,h^*})_{ (\Omega,h^*)\in E^{\tilde{G}}_{ell,\mu}}$ contient $\underline{{\cal E}}_{ell,0,\mu}(\tilde{G})$. Notons $PW^{\infty}_{ell,\mu}(\tilde{G})$ le sous-espace de $PW^{\infty}_{ell}(\tilde{G})$ form\'e des familles  $(\varphi_{\boldsymbol{\tau}})_{\boldsymbol{\tau}\in \underline{{\cal E}}_{ell,0}(\tilde{G})}$ telles que $\varphi_{\boldsymbol{\tau}}=0$ si $\mu(\boldsymbol{\tau})\not=\mu$. Notons $PW'_{ell,\mu}(\tilde{G})$ le sous-espace de $PW'_{ell}(\tilde{G})$ form\'e des familles $(\varphi_{\Omega,h^*})_{(\Omega,h^*)\in E^{\tilde{G}}_{ell}}$ telles que $\varphi_{\Omega,h^*}=0$ si $(\Omega,h^*)\not\in E^{\tilde{G}}_{ell,\mu}$. Il est clair que l'inclusion de l'espace engendr\'e par les $\tilde{\pi}_{\Omega,h^*}$ pour $(\Omega,h^*)\in E^{\tilde{G}}_{ell,\mu}$ dans celui engendr\'e par les $\tilde{\pi}_{\boldsymbol{\tau}}$ pour $\boldsymbol{\tau}\in \underline{{\cal E}}_{ell,0,\mu}(\tilde{G})$ induit une projection $PW^{\infty}_{ell,\mu}(\tilde{G})\to PW'_{ell,\mu}(\tilde{G})$. L'inclusion pr\'ec\'edente est surjective si et seulement si cette projection est injective. Supposons que ce ne soit pas le cas. En utilisant le th\'eor\`eme de Paley-Wiener pour les fonctions $K$-finies ([DM] repris en [W] 6.2), on peut construire une fonction $K$-finie $f\in C_{c}^{\infty}(\tilde{G}({\mathbb R}))$ dont l'image dans $PW^{\infty}(\tilde{G})$ est non nulle mais appartient au noyau de la projection. La premi\`ere propri\'et\'e implique que l'image de $f$ dans $I(\tilde{G}({\mathbb R}))$ est non nulle. La seconde entra\^{\i}ne que son image dans $PW'(\tilde{G})$ est nulle. D'apr\`es le th\'eor\`eme de Renard, l'image de $f$ dans $I(\tilde{G}({\mathbb R}))$ est nulle. Cette contradiction prouve (5). $\square$.

On a donc deux matrices de changement de base $(a_{(\Omega,h^*),\boldsymbol{\tau}})_{(\Omega,h^*)\in E^{\tilde{G}}_{ell},\boldsymbol{\tau}\in \underline{{\cal E}}_{ell,0}(\tilde{G})}$ et $(b_{\boldsymbol{\tau},(\Omega,h^*)})_{\boldsymbol{\tau}\in \underline{{\cal E}}_{ell,0}(\tilde{G}),(\Omega,h^*)\in E^{\tilde{G}}_{ell}}$, inverses l'une de l'autre, qui font passer de la base $(\tilde{\pi}_{\Omega,h^*})_{ (\Omega,h^*)\in E^{\tilde{G}}_{ell}}$ de $D_{ell,0}(\tilde{G}({\mathbb R}))$ \`a la base $(\tilde{\pi}_{\boldsymbol{\tau}})_{\boldsymbol{\tau}\in \underline{{\cal E}}_{ell,0}(\tilde{G})}$.  Comme on l'a dit ci-dessus, ces matrices induisent pour tout $\mu$ un isomorphisme de $PW^{\infty}_{ell,\mu}(\tilde{G})$ sur $PW'_{ell,\mu}(\tilde{G})$. Pour prouver que ces isomorphismes se globalisent en un isomorphisme de $PW^{\infty}_{ell}(\tilde{G})$ sur $PW'_{ell}(\tilde{G})$, on voit qu'il suffit de prouver que les matrices de changement de base v\'erifient les deux propri\'et\'es suivantes:

(6) il existe un entier $N$ tel que, pour tout $(\Omega,h^*)\in E^{\tilde{G}}_{ell}$, l'ensemble des $\boldsymbol{\tau}\in  \underline{{\cal E}}_{ell,0}(\tilde{G})$ tels que $a_{(\Omega,h^*),\boldsymbol{\tau}}\not=0$ a au plus $N$ \'el\'ements et, pour tout $\boldsymbol{\tau}\in  \underline{{\cal E}}_{ell,0}(\tilde{G})$, l'ensemble des $(\Omega,h^*)\in E^{\tilde{G}}_{ell}$ tels que $b_{\boldsymbol{\tau},(\Omega,h^*)}\not=0$ a au plus $N$ \'el\'ements;

(7) il existe une constante $C>0$ telle que, pour tous $(\Omega,h^*)\in E^{\tilde{G}}_{ell}$ et $\boldsymbol{\tau}\in  \underline{{\cal E}}_{ell,0}(\tilde{G})$, on a $\vert a_{(\Omega,h^*),\boldsymbol{\tau}}\vert \leq C$ et $\vert b_{\boldsymbol{\tau},(\Omega,h^*)}\vert \leq C$.

Les matrices se d\'ecomposent en blocs selon les param\`etres $\mu$. L'assertion (6) r\'esulte alors de ce que l'on a d\'ej\`a dit: le nombre d'\'el\'ements de  $\underline{{\cal E}}_{ell,0,\mu}(\tilde{G})$ est born\'e ind\'ependamment de $\mu$. La base $(\tilde{\pi}_{\boldsymbol{\tau}})_{\boldsymbol{\tau}\in \underline{{\cal E}}_{ell,0}(\tilde{G})}$ v\'erifie des propri\'et\'es analogues \`a (2) et (3), \`a savoir:

(8) pour deux \'el\'ements distincts $\boldsymbol{\tau}_{1}$, $\boldsymbol{\tau}_{2}$, on a $(\tilde{\pi}_{\boldsymbol{\tau}_{1}},\tilde{\pi}_{\boldsymbol{\tau}_{2}})_{ell}=0$;

(9) quand $\boldsymbol{\tau}$ parcourt $\underline{{\cal E}}_{ell,0}(\tilde{G})$, les produits $(\tilde{\pi}_{\boldsymbol{\tau}},\tilde{\pi}_{\boldsymbol{\tau}})_{ell}$ ne prennent qu'un nombre fini de valeurs r\'eelles strictement positives.

Cela r\'esulte de [W] th\'eor\`eme 7.3. Il  r\'esulte de (2) et (8) que les coefficients de nos matrices sont des rapports
$$\frac{(\tilde{\pi}_{\boldsymbol{\tau}},\tilde{\pi}_{\Omega,h^*})_{ell}}{(\tilde{\pi}_{\boldsymbol{\tau}},\tilde{\pi}_{\boldsymbol{\tau}})_{ell}}$$
ou
$$\frac{(\tilde{\pi}_{\Omega,h^*},\tilde{\pi}_{\boldsymbol{\tau}})_{ell}}{(\tilde{\pi}_{\Omega,h^*},\tilde{\pi}_{\Omega,h^*})_{ell}}.$$
Ils sont born\'es par 
$$\frac{(\tilde{\pi}_{\Omega,h^*},\tilde{\pi}_{\Omega,h^*})_{ell}^{1/2}}{(\tilde{\pi}_{\boldsymbol{\tau}},\tilde{\pi}_{\boldsymbol{\tau}})_{ell}^{1/2}}$$
ou par l'inverse de ce rapport. Ces deux rapports sont born\'es d'apr\`es (3) et (9). 

En conclusion, le changement de base identifie $PW'_{ell}(\tilde{G})$ et $PW^{\infty}_{ell}(\tilde{G})$. Il en r\'esulte une identification de $PW'(\tilde{G})$ avec $PW^{\infty}(\tilde{G})$. Il est clair que l'application du th\'eor\`eme 1.4 est la compos\'ee de l'application de Renard de $C_{c}^{\infty}(\tilde{G}({\mathbb R}))$ dans $PW'(\tilde{G})$ avec cet isomorphisme entre les deux espaces de Paley-Wiener. Le th\'eor\`eme 1.4 r\'esulte ainsi de celui de Renard. 

\bigskip

\subsection{Extension au cas $\omega\not=1$}
La m\'ethode est la m\^eme qu'en [W] 6.3. On suppose d'abord qu'il existe un caract\`ere $\mu$ de $G({\mathbb R})$ tel que $\omega=\mu\circ(1-\theta)$. Notons ${\bf 1}$ le caract\`ere trivial de $G({\mathbb R})$. A toute $\omega$-repr\'esentation $\tilde{\pi}$ de $\tilde{G}({\mathbb R})$, on associe la ${\bf 1}$-repr\'esentation $\tilde{\pi}_{1}$ d\'efinie par $\tilde{\pi}_{1}(g\gamma_{0})=\mu(g)\tilde{\pi}(g\gamma_{0})$ pour tout $g\in G({\mathbb R})$. C'est une bijection. Construisons un espace ${\cal F}(\tilde{G},\omega)$ comme on a construit $PW^{\infty}(\tilde{G},\omega)$, en oubliant les conditions de convergence. La bijection pr\'ec\'edente  induit un isomorphisme de ${\cal F}(\tilde{G},\omega)$ sur ${\cal F}(\tilde{G},{\bf 1})$, qui se restreint en un hom\'eomorphisme de $PW^{\infty}(\tilde{G},\omega)\to PW^{\infty}(\tilde{G},{\bf 1})$. On d\'efinit d'autre part une application
$$\begin{array}{ccc}C_{c}^{\infty}(\tilde{G}({\mathbb R}))&\to&C_{c}^{\infty}(\tilde{G}({\mathbb R}))\\ f&\mapsto& f_{1}\\ \end{array}$$
par $f_{1}(g\gamma_{0})=\mu(g)^{-1}f(g\gamma_{0})$. Notons $pw_{\tilde{G},\omega}:C_{c}^{\infty}(\tilde{G}({\mathbb R}))\to {\cal F}(\tilde{G},\omega)$ l'application du th\'eor\`eme 1.4 pour le caract\`ere $\omega$. Le th\'eor\`eme affirme que son image est $PW^{\infty}(\tilde{G},\omega)$ et que $pw_{\tilde{G},\omega}$ se factorise en un hom\'eomorphisme de $I(\tilde{G}({\mathbb R}),\omega)$ sur cette image. On v\'erifie que $pw_{\tilde{G},\omega}$ est la compos\'ee de l'application pr\'ec\'edente, de $pw_{\tilde{G},{\bf 1}}$ et de l'isomorphisme inverse ${\cal F}(\tilde{G},{\bf 1})\to {\cal F}(\tilde{G},\omega)$. Il est alors facile de d\'eduire les propri\'et\'es requises de l'application $pw_{\tilde{G},\omega}$ de celles d\'ej\`a connues de l'application $pw_{\tilde{G},{\bf 1}}$.

Dans le cas g\'en\'eral, on  peut \'evidemment supposer $\omega$ trivial sur $Z(G)^{\theta}({\mathbb R})$ (sinon le th\'eor\`eme affirme que $\{0\}$ est hom\'eomorphe \`a $\{0\}$). On introduit un couple $(G',\tilde{G}')$ comme en [W] 2.4. On a une suite exacte
$$1\to C\to G'\stackrel{p}{\to}G\to 1$$
o\`u $C$ est un tore central. On a une application $\tilde{p}:\tilde{G}'\to \tilde{G}$ compatible avec $p$. L'application $G'({\mathbb R})\to G({\mathbb R})$ est surjective. Enfin, il existe un caract\`ere $\mu'$ de $G'({\mathbb R})$ tel que $\omega\circ p=\mu'\circ(1-\theta')$ (o\`u $\theta'$ est l'analogue de $\theta$ pour $(G',\tilde{G}')$).  Notons que $\theta'$ n'est pas trivial sur $C$ en g\'en\'eral. Une $\omega$-repr\'esentation de $\tilde{G}({\mathbb R})$ s'identifie \`a une $\omega\circ p$-repr\'esentation $\tilde{\pi}'$ de $\tilde{G}'({\mathbb R})$ telle que la repr\'esentation sous-jacente $\pi'$  ait un caract\`ere central trivial sur $C({\mathbb R})$. Cette identification induit une application
${\cal F}(\tilde{G}',\omega\circ p)\to {\cal F}(\tilde{G},\omega)$. On voit que cette application se restreint en une application continue  $PW^{\infty}(\tilde{G}',\omega\circ p)\to PW^{\infty}(\tilde{G},\omega)$. On peut construire une section de ces applications de la fa\c{c}on suivante.   Fixons une fonction $\varphi_{C}$ sur ${\cal A}_{C,{\mathbb C}}^{\theta',*}$ qui est de Paley-Wiener et telle que $\varphi_{C}(0)=1$. Soit $\boldsymbol{\varphi}$ un \'el\'ement de ${\cal F}(\tilde{G},\omega)$. C'est donc une collection de fonctions $\varphi_{\tilde{L},\boldsymbol{\tau}}$, o\`u $\tilde{L}$ parcourt ${\cal L}(\tilde{M}_{0})$ et $\boldsymbol{\tau}$ parcourt $\underline{{\cal E}}_{ell,0}(\tilde{L},\omega)$, cette collection \'etant soumise \`a une condition d'invariance par $W(\tilde{M}_0)$. Soit $\tilde{L}\in {\cal L}(\tilde{M}_{0})$. Soit $\tilde{L}'$ son image r\'eciproque dans $\tilde{G}'$. On peut supposer que $\underline{{\cal E}}_{ell,0}(\tilde{L},\omega)$ s'identifie \`a un sous-ensemble de $\underline{{\cal E}}_{ell,0}(\tilde{L}',\omega\circ p)$. On peut identifier ${\cal A}_{\tilde{L}'}^*$ \`a ${\cal A}_{\tilde{L}}^*\oplus {\cal A}_{C}^{\theta',*}$. Pour $\boldsymbol{\tau}'\in \underline{{\cal E}}_{ell,0}(\tilde{L}',\omega\circ p)-\underline{{\cal E}}_{ell,0}(\tilde{L},\omega)$, on pose $\varphi'_{\tilde{L}',\boldsymbol{\tau}'}=0$. Pour $\boldsymbol{\tau}\in \underline{{\cal E}}_{ell,0}(\tilde{L},\omega)$, on d\'efinit une fonction $\varphi'_{\tilde{L}',\boldsymbol{\tau}}$ sur ${\cal A}_{\tilde{L}',{\mathbb C}}^*$ par
$$\varphi'_{\tilde{L}',\boldsymbol{\tau}}(\lambda_{C}+\lambda_{\tilde{L}})=\varphi_{C}(\lambda_{C})\varphi_{\tilde{L},\boldsymbol{\tau}}(\lambda_{\tilde{L}})$$
pour tous $\lambda_{C}\in {\cal A}_{C,{\mathbb C}}^{\theta',*}$ et $\lambda_{\tilde{L}}\in {\cal A}_{\tilde{L},{\mathbb C}}^*$.  La collection $\boldsymbol{\varphi}'$ de ces fonctions $\varphi'_{\tilde{L}',\boldsymbol{\tau}'}$ appartient \`a ${\cal F}(\tilde{G}',\omega\circ p)$. Il est clair que l'application $\boldsymbol{\varphi}\mapsto \boldsymbol{\varphi}'$ est une section de l'application ${\cal F}(\tilde{G}',\omega\circ p)\to {\cal F}(\tilde{G},\omega)$ ci-dessus et qu'elle envoie continuement $PW^{\infty}(\tilde{G},\omega)$ dans $PW^{\infty}(\tilde{G}',\omega\circ p)$. D'autre part, on d\'efinit une application
$$\begin{array}{ccc}C_{c}^{\infty}(\tilde{G}'({\mathbb R}))&\to&C_{c}^{\infty}(\tilde{G}({\mathbb R}))\\ f'&\mapsto&f\\ \end{array}$$
par
$$f(\gamma)=\int_{C({\mathbb R})}f(c\gamma')\,dc,$$
o\`u $\gamma'$ est un rel\`evement quelconque de $\gamma$ dans $\tilde{G}'({\mathbb R})$. Cette application est continue et admet clairement une section continue.

Le diagramme suivant est commutatif
$$\begin{array}{ccc}C_{c}^{\infty}(\tilde{G}'({\mathbb R}))&\stackrel{f'\mapsto f}{\to}&C_{c}^{\infty}(\tilde{G}({\mathbb R}))\\ \qquad \downarrow pw_{\tilde{G}',\omega\circ p}&&\qquad \downarrow pw_{\tilde{G},\omega}\\ {\cal F}(\tilde{G}',\omega\circ p)&\to&{\cal F}(\tilde{G},\omega)\\ \end{array}$$
On conna\^{\i}t d\'ej\`a le th\'eor\`eme pour $(G',\tilde{G}')$ et pour le caract\`ere $\omega\circ p$, puisque celui-ci est de la forme $\mu'\circ(1-\theta')$. Parce que l'application horizontale du haut est surjective, que l'application $pw_{\tilde{G}',\omega\circ p}$ a pour image $PW^{\infty}(\tilde{G}',\omega\circ
p)$ et que l'application horizontale du bas envoie ce dernier espace sur $PW^{\infty}(\tilde{G},\omega)$, on voit que $pw_{\tilde{G},\omega}$ envoie surjectivement $C_{c}^{\infty}(\tilde{G}({\mathbb R}))$ sur $PW^{\infty}(\tilde{G},\omega)$.  Parce que l'application horizontale du haut admet une section continue, la continuit\'e de l'application $pw_{\tilde{G},\omega}:C_{c}^{\infty}(\tilde{G}({\mathbb R}))\to PW^{\infty}(\tilde{G},\omega)$ r\'esulte de celle de l'application similaire $pw_{\tilde{G}',\omega\circ p}$ et de celle de l'application $PW^{\infty}(\tilde{G}',\omega\circ p)\to PW^{\infty}(\tilde{G},\omega)$. Comme on l'a d\'ej\`a dit, l'application $pw_{\tilde{G},\omega}$ se factorise par $I(\tilde{G}({\mathbb R}),\omega)$. Parce que l'application $C_{c}^{\infty}(\tilde{G}({\mathbb R}))\to I(\tilde{G}({\mathbb R}),\omega)$ est ouverte, notre application continue $pw_{\tilde{G},\omega}$ se factorise en une application continue
$$(1)\qquad I(\tilde{G}({\mathbb R}),\omega)\to PW^{\infty}(\tilde{G},\omega)$$
qui est encore surjective.  On sait par ailleurs qu'elle est injective ([W] th\'eor\`eme 5.5), donc bijective. Il reste \`a prouver que l'application r\'eciproque est continue.
On voit ais\'ement que le diagramme ci-dessus se factorise en un diagramme commutatif
 $$\begin{array}{ccc}C_{c}^{\infty}(\tilde{G}'({\mathbb R}))&\stackrel{f'\mapsto f}{\to}&C_{c}^{\infty}(\tilde{G}({\mathbb R}))\\ \downarrow&&\downarrow\\ I(\tilde{G}'({\mathbb R}),\omega\circ p)&\to&I(\tilde{G}({\mathbb R}),\omega)\\ \downarrow&&\downarrow\\ PW^{\infty}(\tilde{G}',\omega\circ p)&\to&PW^{\infty}(\tilde{G},\omega)\\ \end{array}$$
 L'application horizontale du milieu est continue car l'application du haut l'est, l'application verticale du haut \`a droite est continue et l'application verticale du haut \`a gauche est ouverte. L'application verticale du bas \`a gauche est un hom\'eomorphisme. Puisque l'application horizontale du bas admet une section continue, il en est de m\^eme de l'application verticale du bas \`a droite. C'est-\`a-dire que l'inverse de l'application (1) est continue. Cela ach\`eve la preuve. 
 
 \bigskip
 
 \section{Stabilit\'e}
 
 \bigskip
 
 \subsection{Quelques consid\'erations formelles}
 Soit ${\bf G}'=(G',{\cal G}',\tilde{s})$ une donn\'ee endoscopique relevante de $(G,\tilde{G},{\bf a})$. On doit fixer une application $H_{\tilde{G}'}:\tilde{G}'({\mathbb R})\to {\cal A}_{\tilde{G}'}={\cal A}_{G'}$ comme en 1.1. Supposons que ${\bf G}'$ est elliptique. On a alors un isomorphisme $\xi:{\cal A}_{\tilde{G}}\simeq{\cal A}_{\tilde{G}'}$ et
 
 (1) il existe une unique application $H_{\tilde{G}'}:\tilde{G}'({\mathbb R})\to {\cal A}_{\tilde{G}'}$ telle que 
 
 (i) $H_{\tilde{G}'}(y\delta)=H_{G'}(y)+H_{\tilde{G}'}(\delta)$ pour tous $y\in G'({\mathbb R})$ et $\delta\in \tilde{G}'({\mathbb R})$;
 
 (ii)  on a $H_{\tilde{G}'}(\delta)=\xi(H_{\tilde{G}}(\gamma))$ pour tout couple $(\delta,\gamma)\in \tilde{G}'({\mathbb R})\times \tilde{G}({\mathbb R})$  d'\'el\'ements semi-simples qui se correspondent.
 
  Preuve. Fixons des \'el\'ements semi-simples $\delta^{\natural}\in \tilde{G}'({\mathbb R})$ et $\gamma^{\natural}\in \tilde{G}({\mathbb R})$ qui se correspondent.  On d\'efinit $H_{\tilde{G}'}$ comme l'unique application v\'erifiant (i) et telle que $H_{\tilde{G}'}(\delta^{\natural})=\xi(H_{\tilde{G}}(\gamma^{\natural}))$. 
 On doit montrer que cette application  v\'erifie (ii).  Fixons un diagramme $(\delta^{\natural},B^{_{'}\natural},T^{_{'}\natural},B^{\natural},T^{\natural},\gamma^{\natural})$. Notons $\xi^{\natural}:T^{\natural}\to T^{_{'}\natural}$ l'homomorphisme qui s'en d\'eduit. Compl\'etons nos paires de Borel en des paires de Borel \'epingl\'ees ${\cal E}^{\natural}$ et ${\cal E}^{_{'}\natural}$. Ecrivons $\gamma^{\natural}=t^{\natural}e^{\natural}$, avec $t^{\natural}\in T^{\natural}$ et $e^{\natural}\in Z(\tilde{G},{\cal E}^{\natural})$, notons $e^{_{'}\natural}$ l'image naturelle de $e^{\natural}$ dans $Z(\tilde{G}',{\cal E}^{_{'}\natural})$, \'ecrivons $\delta^{\natural}=t^{_{'}\natural}e^{_{'}\natural}$ avec $t^{_{'}\natural}\in T^{_{'}\natural}$. On a alors $\xi(t^{\natural})=t^{_{'}\natural}$. Consid\'erons un autre diagramme quelconque $(\delta,B',T',B,T,\gamma)$. Fixons des \'el\'ements $x\in G$ et $y\in G'$ tels que $ad_{x}(B^{\natural},T^{\natural})=(B,T)$ et $ad_{y}(B^{_{'}\natural},T^{_{'}\natural})=(B',T')$. Posons $e=ad_{x}(e^{\natural})$ et $e'=ad_{y}(e^{_{'}\natural})$. On \'ecrit $\gamma=te$ et $\delta=t'e'$, avec $t\in T$ et $t'\in T'$. On a $\xi(t)=t'$, o\`u $\xi:T\to T'$ est d\'eduit du diagramme. On a $\gamma=ad_{x}(ue^{\natural})$ et $\delta=ad_{y}(u'e^{_{'}\natural})$, o\`u $u=ad_{x^{-1}}(t)$ et $u'=ad_{y^{-1}}(t')$. Puisque $ad_{y^{-1}}\circ\xi\circ ad_{x}=\xi^{\natural}$, on a $\xi^{\natural}(u)=u'$. Ecrivons $\gamma=g\gamma^{\natural}$ et $\delta=g'\delta^{\natural}$ avec $g\in G({\mathbb R})$ et $g'\in G'({\mathbb R})$. 
Pour prouver (ii), on doit prouver que $\xi(H_{\tilde{G}}(g))=H_{\tilde{G}'}(g')$.   Pour cela, il suffit de prouver que tout caract\`ere $\chi'\in X^{*}(G')^{\Gamma_{{\mathbb R}}}$ prend la m\^eme valeur sur les deux termes. Rappelons que $X^{*}(G')^{\Gamma_{{\mathbb R}}}\simeq X^{*}(G)^{\Gamma_{{\mathbb R}},\theta}$. Pr\'ecis\'ement, cet isomorphisme associe \`a $\chi'\in X^{*}(G')^{\Gamma_{{\mathbb R}}}$ l'unique \'el\'ement $\chi\in X^{*}(G)^{\Gamma_{{\mathbb R}},\hat{\theta}}$ tel que $\chi'\circ \xi^{\natural}=\chi$ sur $T^{\natural}$. Soit $\chi'\in X^{*}(G')^{\Gamma_{{\mathbb R}}}$ et $\chi$ l'\'el\'ement associ\'e. On calcule $g=xu\theta(x)^{-1}(t^{\natural})^{-1}$ (o\`u $\theta=ad_{e^{\natural}}$) et $g'=yu'y^{-1}(t^{_{'}\natural})^{-1}$. D'o\`u
$$\chi'(H_{\tilde{G}'}(g'))=log(\vert \chi'(g')\vert _{{\mathbb R}})=log(\vert \chi'(u'(t^{_{'}\natural})^{-1})\vert _{{\mathbb R}}),$$
$$\chi'\circ\xi(H_{\tilde{G}}(g))=\chi(H_{\tilde{G}}(g))=log(\vert \chi(u(t^{\natural})^{-1})\vert _{{\mathbb R}}).$$
Puisque $\xi^{\natural}(u(t^{\natural})^{-1})=u'(t^{_{'}\natural})^{-1}$, ces deux expressions sont \'egales, ce qui prouve (1). $\square$

 Dans le cas o\`u ${\bf G}'$ est elliptique, on identifie ${\cal A}_{\tilde{G}'}$ \`a ${\cal A}_{\tilde{G}}$ par $\xi$ et on choisit pour $H_{\tilde{G}'}$ l'application d\'efinie par (1).  
 
  On a expliqu\'e en [I] 2.5 comment d\'efinir des espaces $C_{c}^{\infty}({\bf G}')$, $I({\bf G}')$, $SI({\bf G}')$. Cette construction s'adapte aux espaces de distributions d\'efinis en 1.2. Pr\'ecis\'ement,  on a d\'efini en [I] 2.1 et 2.5 la notion de donn\'ees auxiliaires $G'_{1}$, $\tilde{G}'_{1}$, $C_{1}$, $\hat{\xi}_{1}$, $\Delta_{1}$. Il est facile de voir que l'on peut trouver de telles donn\'ees telles que le caract\`ere associ\'e $\lambda_{1}$  soit unitaire. On supposera toujours qu'il en est ainsi. Consid\'erons de telles donn\'ees. Notons $D_{spec,\lambda_{1}}(\tilde{G}'_{1}({\mathbb R}))$ le sous-espace de $D_{spec}(\tilde{G}'_{1}({\mathbb R}))$ engendr\'e par les caract\`eres de repr\'esentations $(\pi_{1},\tilde{\pi}_{1})$ de $\tilde{G}'_{1}({\mathbb R})$ (pour le caract\`ere trivial de $G'_{1}({\mathbb R})$) telles que le caract\`ere central de $\pi_{1}$ co\"{\i}ncide avec $\lambda_{1}$ sur $C_{1}({\mathbb R})$. Consid\'erons une telle repr\'esentation $(\pi_{1},\tilde{\pi}_{1})$ et d'autres donn\'ees auxiliaires $G'_{2}$,..., $\Delta_{2}$. On d\'efinit une repr\'esentation $(\pi_{2},\tilde{\pi}_{2})$ de $\tilde{G}'_{2}({\mathbb R})$, agissant dans le m\^eme espace que $(\pi_{1},\tilde{\pi}_{1})$, par les formules
 $$\pi_{2}(x_{2})=\lambda_{12}(x_{1},x_{2})^{-1}\pi_{1}(x_{1}),\,\, \tilde{\pi}_{2}(\delta_{2})=\tilde{\lambda}_{12}(\delta_{1},\delta_{2})^{-1}\tilde{\pi}_{1}(\delta_{1})$$
 pour $x_{2}\in G'_{2}({\mathbb R})$ et $\delta_{2}\in \tilde{G}'_{2}({\mathbb R})$, o\`u $x_{1}$ est un \'el\'ement quelconque de $G'_{1}({\mathbb R})$ qui a m\^eme projection que $x_{2}$ dans $G'({\mathbb R})$, o\`u $\delta_{1}$ est un \'el\'ement quelconque de $\tilde{G}'_{1}({\mathbb R})$ qui a m\^eme projection que $\delta_{2}$ dans $\tilde{G}'({\mathbb R})$ et o\`u $\lambda_{12}$ et $\tilde{\lambda}_{12}$ sont les fonctions d\'efinies en [I] 2.5.  Alors  le caract\`ere central de $\pi_{2}$ co\"{\i}ncide avec $\lambda_{2}$ sur $C_{2}({\mathbb R})$. L'application qui, au caract\`ere de $(\pi_{1},\tilde{\pi}_{1})$, associe celui de $(\pi_{2},\tilde{\pi}_{2})$, se prolonge en un isomorphisme de $D_{spec,\lambda_{1}}(\tilde{G}'_{1}({\mathbb R}))$  sur $D_{spec,\lambda_{2}}(\tilde{G}'_{2}({\mathbb R}))$. En recollant par ces isomorphismes canoniques les espaces $D_{spec,\lambda_{1}}(\tilde{G}'_{1}({\mathbb R}))$  associ\'es \`a toutes les donn\'ees auxiliaires possibles, on obtient un espace que l'on notera $D_{spec}({\bf G}')$. On d\'efinit de m\^eme le sous-espace $D_{temp}({\bf G}')$. On v\'erifie que l'application ci-dessus qui, \`a $(\pi_{1},\tilde{\pi}_{1})$, associe $(\pi_{2},\tilde{\pi}_{2})$, envoie une repr\'esentation elliptique sur une repr\'esentation elliptique. Il lui est associ\'e une bijection de ${\cal E}_{ell,\lambda_{1}}(\tilde{G}'_{1} )$ sur ${\cal E}_{ell,\lambda_{2}}(\tilde{G}'_{2})$, les indices $\lambda_{1}$ et $\lambda_{2}$ signifiant que l'on se restreint aux \'el\'ements dont le caract\`ere central se restreint \`a $C_{1}({\mathbb R})$, resp. $C_{2}({\mathbb R})$, en le caract\`ere $\lambda_{1}$, resp. $\lambda_{2}$. On en d\'eduit comme ci-dessus par recollement un ensemble ${\cal E}_{ell}({\bf G}')$ et un espace $D_{ell}({\bf G}')$. 
 
  {\bf Remarque.} Le corps de base est ici ${\mathbb R}$, mais les constructions ci-dessus valent aussi sur un corps de base local non-archim\'edien.

  Consid\'erons de nouveau des donn\'ees auxiliaires $G'_{1}$, ..., $\Delta_{1}$. Introduisons les espaces $\mathfrak{h}'$, $\mathfrak{h}'_{1}$ et $\mathfrak{h}_{C_{1}}$ analogues de $\mathfrak{h}$ pour les groupes $G'$, $G'_{1}$ et $C_{1}$. On a une suite exacte
  $$0\to \mathfrak{h}^{_{'}*}\to \mathfrak{h}_{1}^{_{'}*}\to \mathfrak{h}_{C_{1}}^*\to 0$$
  Au caract\`ere $\lambda_{1}$ est associ\'e un param\`etre $\mu(\lambda_{1})\in \mathfrak{h}_{C_{1}}^*$. Soit $(\pi_{1},\tilde{\pi}_{1})$ une repr\'esentation irr\'eductible de $\tilde{G}'_{1}({\mathbb R})$ dont le caract\`ere central co\"{\i}ncide avec $\lambda_{1}$ sur $C_{1}({\mathbb R})$. Alors son param\`etre $\mu(\tilde{\pi}_{1})$ appartient \`a l'ensemble $\mathfrak{h}_{1,\lambda_{1}}^{_{'}*}$ des \'el\'ements de $\mathfrak{h}_{1}^{_{'}*}$ qui se projettent sur $\mu(\lambda_{1})$. Cet ensemble est un espace affine sous $\mathfrak{h}^{_{'}*}$. Quand on change de donn\'ees auxiliaires, ces espaces affines se recollent. En effet, avec les notations ci-dessus, $\lambda_{12}$ est un caract\`ere du produit fibr\'e de $G'_{1}$ et $G'_{2}$ au-dessus de $G'$. Il lui est associ\'e un param\`etre 
  $$\mu(\lambda_{12})\in \mathfrak{h}_{12}^{_{'}*}=(\mathfrak{h}_{1}^{_{'}*}\times \mathfrak{h}_{2}^{_{'}*})/diag_{-}(\mathfrak{h}^{_{'}*})$$
  (on note $diag_{-}$ le plongement antidiagonal). La projection de $\mu(\lambda_{12})$ dans $\mathfrak{h}_{C_{1}}\times \mathfrak{h}_{C_{2}}$ est $(\mu(\lambda_{1}),-\mu(\lambda_{2}))$. Soit alors $\mu_{1}\in \mathfrak{h}_{1,\lambda_{1}}^{_{'}*}$. Il existe un unique \'el\'ement $\mu_{2}\in \mathfrak{h}_{2,\lambda_{2}}^{_{'}*}$ tel que $(\mu_{1},-\mu_{2})$ se projette sur $\mu(\lambda_{12})$ dans $\mathfrak{h}_{12}^{_{'}*}$. Le recollement associe $\mu_{2}$ \`a $\mu_{1}$. Par recollement, on obtient un espace affine sous $\mathfrak{h}^{_{'}*}$ que nous noterons $\mathfrak{h}^{{\bf G}',*}$.  On ne peut pas recoller les alg\`ebres $\mathfrak{Z}(G'_{1})$, mais on peut recoller leurs quotients $\mathfrak{Z}(G'_{1})/I(\lambda_{1})$, o\`u $I(\lambda_{1})$ est l'id\'eal bilat\`ere engendr\'e par les $X-<X,\mu(\lambda_{1})>$ pour $X\in \mathfrak{h}_{C_{1}}$. On obtient une alg\`ebre not\'ee $\mathfrak{Z}({\bf G}')$, qui s'identifie \`a l'alg\`ebre des polyn\^omes sur $\mathfrak{h}^{{\bf G}',*}$ invariants par $W^{G'}$.  Cette alg\`ebre agit naturellement sur  $C_{c}^{\infty}({\bf G}')$.

  On a d\'efini en 1.2 l'espace affine $\tilde{\mu}(\omega)+\mathfrak{h}^{\theta,*}$. Il s'identifie \`a $ \mathfrak{h}^{{\bf G}',*}$ par la construction suivante. Notons que $\mathfrak{h}'$ s'identifie naturellement \`a $\mathfrak{h}^{\theta}$ (l'identification d\'epend de choix de paires de Borel mais changer ces choix ne modifie l'identification que par l'action d'un \'el\'ement de $W^{\theta}$, ce qui nous importe peu).  Fixons des donn\'ees auxiliaires $G'_{1}$,...,$\Delta_{1}$. On a d\'efini en [I] 2.8 un \'el\'ement $b$ que l'on peut consid\'erer comme un \'el\'ement de
   $$(1) \qquad (\mathfrak{h}_{1}^{_{'}*}\times \mathfrak{h}^*)/diag_{-}(\mathfrak{h}^{_{'}*}).$$
    On v\'erifie que sa projection dans $\mathfrak{h}_{C_{1}}^*\times (1-\theta)(\mathfrak{h}^*)$ n'est autre que $(-\mu(\lambda_{1}),\tilde{\mu}(\omega))$. Pour $\mu\in \tilde{\mu}(\omega)+\mathfrak{h}^{\theta,*}$, il existe un unique $\mu_{1}\in \mathfrak{h}_{1,\lambda_{1}}^{_{'}*}$ tel que $(-\mu_{1},\mu)$ se projette sur $b$ dans l'espace (1). L'application $\mu\mapsto \mu_{1}$ fournit l'isomorphisme cherch\'e. Modulo cet isomorphisme, on a un homomorphisme 
    $$\begin{array}{ccc}\mathfrak{Z}(G)&\to& \mathfrak{Z}({\bf G}')\\ z&\mapsto &z^{{\bf G}'}\\ \end{array}$$ Un \'el\'ement $z$ de $\mathfrak{Z}(G)$ d\'efinit par restriction un polyn\^ome invariant par $W^{\theta}$ sur $\tilde{\mu}(\omega)+\mathfrak{h}^{\theta,*}$, qui s'identifie \`a un polyn\^ome $z^{{\bf G}'}$ invariant par $W^{G'}\subset W^{\theta}$ sur $ \mathfrak{h}^{{\bf G}',*}$. Le corollaire [I] 2.8 se reformule de la fa\c{c}on suivante.
    
    \ass{Lemme}{Soient ${\bf f}\in C_{c}^{\infty}(\tilde{G}({\mathbb R}))\otimes Mes(G({\mathbb R}))$, ${\bf f}'\in C_{c}^{\infty}({\bf G}')\otimes Mes(G'({\mathbb R}))$ et $z\in \mathfrak{Z}(G)$. Supposons que ${\bf f}'$ soit un transfert de ${\bf f}$. Alors $z^{{\bf G}'}{\bf f}'$ est un transfert de $z{\bf f}$.}
    
   Fixons des donn\'ees auxiliaires $G'_{1}$,...,$\Delta_{1}$. On a introduit ci-dessus un \'el\'ement $b$ de l'espace (1).   
   Son oppos\'e $-b$ se projette sur un \'el\'ement de $(\mathfrak{a}_{G'_{1}}^*\times \mathfrak{a}_{\tilde{G}}^*)/diag _{-}(\mathfrak{a}_{\tilde{G}}^*)\simeq \mathfrak{a}_{G'_{1}}^*$. Cet \'el\'ement d\'efinit un caract\`ere $\lambda_{\mathfrak{A}_{G'_{1}}}$ de $\mathfrak{A}_{G'_{1}}$, qui a   m\^eme restriction \`a $\mathfrak{A}_{C_{1}}$ que $\lambda_{1}$. Pour simplifier, fixons des mesures de Haar sur $G({\mathbb R})$ et $G'_{1}({\mathbb R})$.   Soit $a\in \mathfrak{A}_{\tilde{G}}$ et $a_{1}\in \mathfrak{A}_{G'_{1}}$ ayant m\^eme projection dans $\mathfrak{A}_{G'}$. Soient $f\in C_{c}^{\infty}(\tilde{G}({\mathbb R}))$ et $f_{1}\in C_{c,\lambda_{1}}^{\infty}(\tilde{G}'_{1}({\mathbb R}))$. Notons $f_{1}^{a_{1}}$, resp. $f^{a}$, les fonctions d\'efinies par $f_{1}^{a_{1}}(\delta_{1})=f_{1}(a_{1}\delta_{1})$, resp. $f^{a}(\gamma)=f(a\gamma)$. On a 
   
   (2) si $f_{1}$ est un transfert de $f$, alors $\lambda_{\mathfrak{A}_{G'_{1}}}(a_{1})f_{1}^{a_{1}}$ est un transfert de $f^{a}$. 
   
   Notons ${\cal E}_{ell,\lambda_{\mathfrak{A}_{G'_{1}}},\lambda_{1}}(\tilde{G}'_{1})$ le sous-ensemble des $\boldsymbol{\tau}\in {\cal E}_{ell,\lambda_{1}}(\tilde{G}'_{1})$ tels que la restriction \`a $\mathfrak{A}_{G'_{1}}$ du caract\`ere central de $\boldsymbol{\tau}$  soit \'egale \`a $\lambda_{\mathfrak{A}_{G'_{1}}}$. On note de m\^eme $D_{ell,\lambda_{\mathfrak{A}_{G'_{1}}},\lambda_{1}}(\tilde{G}'_{1})$ le sous-espace de $D_{ell,\lambda_{1}}(\tilde{G}'_{1})$ engendr\'e par les caract\`eres des $\tilde{\pi}_{\boldsymbol{\tau}}$ pour $\boldsymbol{\tau}\in {\cal E}_{ell,\lambda_{\mathfrak{A}_{G'_{1}}},\lambda_{1}}(\tilde{G}'_{1})$. Quand on fait varier les donn\'ees auxiliaires, on v\'erifie que ces objets se recollent en des objets que l'on note ${\cal E}_{ell,0}({\bf G}')$ et $D_{ell,0}({\bf G}')$. Il est utile de remarquer que l'on peut choisir des donn\'ees auxiliaires de sorte que $\lambda_{\mathfrak{A}_{G'_{1}}}$ soit trivial. En effet, via la projection naturelle $G'_{1}({\mathbb R})\to \mathfrak{A}_{G'_{1}}$, ce caract\`ere s'\'etend en un caract\`ere de $G'_{1}({\mathbb R})$. Celui-ci d\'etermine un cocycle $\zeta:W_{{\mathbb R}}\to Z(\hat{G}'_{1})$. On d\'efinit un nouveau plongement $\hat{\xi}'_{1}:{\cal G}'\to {^LG}'_{1}$ par $\hat{\xi}'_{1}(g,w)=\zeta(w)^{-1}\hat{\xi}_{1}(g,w)$ pour tout $(g,w)\in {\cal G}'$. On voit que, si l'on remplace $\hat{\xi}_{1}$ par ce nouveau plongement, $\lambda_{\mathfrak{A}_{G'_{1}}}$ devient trivial. 
  \bigskip
  
  \subsection{Les espaces $I_{cusp}(\tilde{G}({\mathbb R}))$ et $SI_{cusp}(\tilde{G}({\mathbb R}))$}
 {\bf Pour simplifier les notations,   on fixe dor\'enavant des mesures de Haar  sur tous les groupes rencontr\'es}.   Le th\'eor\`eme 1.4 entra\^{\i}ne que $I_{cusp}(\tilde{G}({\mathbb R}),\omega)$ est isomorphe \`a l'espace $PW^{\infty}_{ell}(\tilde{G},\omega)$. Rappelons que l'on note $C_{c}^{\infty}(\tilde{G}({\mathbb R}),K)$ l'espace des \'el\'ements de $C_{c}^{\infty}(\tilde{G}({\mathbb R}))$ qui sont $K$-finis \`a droite et \`a gauche.  Notons $C_{cusp}^{\infty}(\tilde{G}({\mathbb R}),K)$ le sous-espace des fonctions cuspidales et $I_{cusp}(\tilde{G}({\mathbb R}),\omega,K)$ son image dans $I_{cusp}(\tilde{G}({\mathbb R}),\omega)$. Le th\'eor\`eme de Delorme et Mezo (repris en [W] 6.2) entra\^{\i}ne que $I_{cusp}(\tilde{G}({\mathbb R}),\omega,K)$ s'identifie au sous-espace $PW_{ell}(\tilde{G},\omega)$ des familles $(\varphi_{\boldsymbol{\tau}})_{\boldsymbol{\tau}\in \underline{{\cal E}}_{ell,0}(\tilde{G},\omega)}$ telles que $\varphi_{\boldsymbol{\tau}}=0$ pour presque tout $\boldsymbol{\tau}\in \underline{{\cal E}}_{ell,0}(\tilde{G},\omega)$. 
 
 On peut modifier les d\'efinitions des espaces de Paley-Wiener de la fa\c{c}on suivante. On a d\'efini l'espace $D_{ell,0}(\tilde{G}({\mathbb R}),\omega)$ engendr\'e par les caract\`eres des repr\'esentations $\tilde{\pi}_{\boldsymbol{\tau}}$ pour $\boldsymbol{\tau}\in {\cal E}_{ell,0}(\tilde{G},\omega)$. Il est muni du produit elliptique, qui est hermitien et d\'efini positif. Pour une $W$-orbite $\mu$ dans $\mathfrak{h}^*$, on a d\'efini le sous-espace $D_{ell,0,\mu}(\tilde{G}({\mathbb R}),\omega)$ engendr\'e par les caract\`eres des $\tilde{\pi}_{\boldsymbol{\tau}}$ comme ci-dessus telles que $\mu(\boldsymbol{\tau})=\mu$. Ce sont des espaces de dimension finie uniform\'ement born\'ee. Consid\'erons une base $B$ de $D_{ell,0}(\tilde{G}({\mathbb R}),\omega)$ qui est r\'eunion de bases $B_{\mu}$ des sous-espaces  $D_{ell,0,\mu}(\tilde{G}({\mathbb R}),\omega)$, qui est orthogonale pour le produit hermitien et qui v\'erifie la condition
 
 (1) quand $\tilde{\pi}$ d\'ecrit $B$, les produits $(\tilde{\pi},\tilde{\pi})_{ell}$ ne prennent qu'un nombre fini de valeurs.
 
 On construit un espace de Paley-Wiener comme en 1.3, associ\'e \`a l'ensemble $B$ muni de la fonction $d$ qui vaut $\vert \mu\vert $ sur chaque $B_{\mu}$, les espaces $V$ \'etant tous \'egaux \`a ${\cal A}_{\tilde{G}}$. La m\^eme preuve qu'en 1.5 montre que cet espace s'identifie \`a $PW^{\infty}_{ell}(\tilde{G},\omega)$. On peut donc consid\'erer un \'el\'ement de cet espace comme une collection de fonctions $(\varphi_{\tilde{\pi}})_{\tilde{\pi}\in B}$. 
 
Tout ceci s'adapte si l'on remplace les espaces de fonctions sur $\tilde{G}({\mathbb R})$ par des espaces de fonctions invariantes par $\mathfrak{A}_{\tilde{G}}$. Par exemple, l'espace $I_{cusp}(\tilde{G}({\mathbb R})/\mathfrak{A}_{\tilde{G}},\omega)$ s'identifie \`a l'espace $PW^{\infty}_{ell,0}(\tilde{G},\omega)$ des familles  $(f_{\tilde{\pi}})_{\tilde{\pi}\in B}\in {\mathbb C}^B$ tels que, pour tout entier $N$, il existe  $C_{N}>0$ de sorte que
$$\vert f_{\tilde{\pi}}\vert \leq C_{N}(1+\vert \mu(\tilde{\pi})\vert )^{-N}$$
pour tout $\tilde{\pi}\in B$. Le sous-espace $I_{cusp}(\tilde{G}({\mathbb R})/\mathfrak{A}_{\tilde{G}},\omega,K)$ s'identifie au sous-espace $PW_{ell,0}(\tilde{G},\omega)$ des familles presque toutes nulles.  Notons que $\mathfrak{Z}(G)$ agit naturellement sur ces espaces.  En particulier, pour $(f_{\tilde{\pi}})_{\tilde{\pi}\in B}$ comme ci-dessus et pour $z\in \mathfrak{Z}(G)$, on a
$$z(f_{\tilde{\pi}})_{\tilde{\pi}\in B}=(z(\mu(\tilde{\pi}))f_{\tilde{\pi}})_{\tilde{\pi}\in B}.$$
On voit que le sous-espace $PW_{ell,0}(\tilde{G},\omega)$ co\"{\i}ncide avec celui des \'el\'ements $\mathfrak{Z}(G)$-finis de $PW^{\infty}_{ell,0}(\tilde{G},\omega)$. 

La formule des traces locale munit $I_{cusp}(\tilde{G}({\mathbb R})/\mathfrak{A}_{\tilde{G}},\omega)$ d'un produit hermitien d\'efini positif. Pour $f,f'\in I_{cusp}(\tilde{G}({\mathbb R})/\mathfrak{A}_{\tilde{G}},\omega)$, on a simplement
$$(f,f')=\sum_{\tilde{T}}\int_{\tilde{T}({\mathbb R})/\mathfrak{A}_{\tilde{G}}(1-\theta)(T({\mathbb R}))}\overline{I^{\tilde{G}}(\gamma,\omega,f)}I^{\tilde{G}}(\gamma,\omega,f')\,d\gamma,$$
o\`u l'on somme sur les classes de conjugaison par $G({\mathbb R})$ de tores tordus maximaux elliptiques et o\`u les mesures sont  convenablement normalis\'ees. 
 La th\'eorie des pseudo-coefficients nous dit qu'il existe un isomorphisme antilin\'eaire  
$$\begin{array}{ccc}D_{ell,0}(\tilde{G},\omega)&\to&  I_{cusp}(\tilde{G}({\mathbb R})/\mathfrak{A}_{\tilde{G}},\omega,K)\\ \tilde{\pi}&\mapsto &f[\tilde{\pi}]\\ \end{array}$$
tel que, pour tout $f\in I_{cusp}(\tilde{G}({\mathbb R})/\mathfrak{A}_{\tilde{G}},\omega)$ et tout $\tilde{\pi}\in D_{ell,0}(\tilde{G},\omega)$, on ait l'\'egalit\'e $I^{\tilde{G}}(\tilde{\pi},\omega,f)=(\bar{f}[\tilde{\pi}],f)$. C'est une isom\'etrie en ce sens que $(\tilde{\pi}_{1},\tilde{\pi}_{2})_{ell}=(f[\tilde{\pi}_{2}],f[\tilde{\pi}_{1}])$. 

Notons $\iota$ l'antiautomorphisme antilin\'eaire de  l'alg\`ebre enveloppante de $G({\mathbb R})$ qui prolonge antilin\'eairement l'application $X\mapsto -X$ de l'alg\`ebre de Lie. Il se restreint en un  automorphisme antilin\'eaire de $\mathfrak{Z}(G)$. En consid\'erant $\mathfrak{Z}(G)$ comme l'alg\`ebre des polyn\^omes sur $\mathfrak{h}^*$ invariants par $W$, on a $(\iota(z))(\lambda)=\overline{z(-\bar{\lambda})}$ pour tout $z\in \mathfrak{Z}(G)$ et tout $\lambda\in \mathfrak{h}^*$.  L'isom\'etrie ci-dessus v\'erifie la relation $f[z\tilde{\pi}]=\iota(z)(f[\tilde{\pi}])$. 
 Pour toute $W$-orbite  $\mu$ dans $\mathfrak{h}^*$, elle identifie $D_{ell,0,\mu}(\tilde{G},\omega)$ avec le sous-espace $I_{cusp,-\bar{\mu}}(\tilde{G}({\mathbb R})/\mathfrak{A}_{\tilde{G}},\omega)$ des \'el\'ements $f\in I_{cusp}(\tilde{G}({\mathbb R})/\mathfrak{A}_{\tilde{G}},\omega)$ tels que $zf=z(-\bar{\mu})f$ pour tout $z\in \mathfrak{Z}(G)$. 

{\bf On suppose pour la suite de la section que $(G,\tilde{G},\omega)$ est quasi-d\'eploy\'e et \`a torsion int\'erieure.}

On note $I_{cusp}^{inst}(\tilde{G}({\mathbb R}))$ le noyau de l'application naturelle $I_{cusp}(\tilde{G}({\mathbb R}))\to SI_{cusp}(\tilde{G}({\mathbb R}))$. On a introduit en [I] 4.14 le sous-espace $I_{cusp}^{st}(\tilde{G}({\mathbb R}))$ des \'el\'ements dont les int\'egrales orbitales sont constantes sur toute classe de conjugaison stable form\'ee d'\'el\'ements elliptiques fortement r\'eguliers. On a montr\'e que, par l'application naturelle ci-dessus, il s'envoyait bijectivement sur $SI_{cusp}(\tilde{G}({\mathbb R}))$. Il en r\'esulte que l'on a la d\'ecomposition
$$I_{cusp}(\tilde{G}({\mathbb R}))=I^{st}_{cusp}(\tilde{G}({\mathbb R}))\oplus I^{inst}_{cusp}(\tilde{G}({\mathbb R})).$$
Il est clair que chacun des sous-espaces est invariant par l'action de $\mathfrak{Z}(G)$. Des d\'efinitions et propri\'et\'es analogues valent pour les fonctions invariantes par $\mathfrak{A}_{\tilde{G}}$: on a
$$I_{cusp}(\tilde{G}({\mathbb R})/\mathfrak{A}_{\tilde{G}})=I^{st}_{cusp}(\tilde{G}({\mathbb R})/\mathfrak{A}_{\tilde{G}})\oplus I^{inst}_{cusp}(\tilde{G}({\mathbb R})/\mathfrak{A}_{\tilde{G}}).$$
Il r\'esulte de la d\'efinition du produit hermitien ci-dessus que la d\'ecomposition est orthogonale. La projection sur chacun des sous-espaces d'un \'el\'ement $\mathfrak{Z}(G)$-fini l'est aussi. On a donc une d\'ecomposition similaire
$$I_{cusp}(\tilde{G}({\mathbb R})/\mathfrak{A}_{\tilde{G}},K)=I^{st}_{cusp}(\tilde{G}({\mathbb R})/\mathfrak{A}_{\tilde{G}},K)\oplus I^{inst}_{cusp}(\tilde{G}({\mathbb R})/\mathfrak{A}_{\tilde{G}},K)$$
o\`u les sous-espaces sont les intersections des pr\'ec\'edents avec $I_{cusp}(\tilde{G}({\mathbb R})/\mathfrak{A}_{\tilde{G}},K)$.  Via l'isomorphisme $\tilde{\pi}\mapsto f[\tilde{\pi}]$, on obtient une d\'ecomposition
$$D_{ell,0}(\tilde{G})=D_{ell,0}^{st}(\tilde{G})\oplus D_{ell,0}^{inst}(\tilde{G}).$$
Pour chaque $W$-orbite $\mu$ dans $\mathfrak{h}^*$, elle se raffine en une d\'ecomposition
$$D_{ell,0,\mu}(\tilde{G})=D_{ell,0,\mu}^{st}(\tilde{G})\oplus D_{ell,0,\mu}^{inst}(\tilde{G}).$$
Pour chaque $\mu$, fixons des bases orthonorm\'ees $B^{st}_{\mu}$ et $B^{inst}_{\mu}$ de chacun des sous-espaces ci-dessus.  Notons $B^{st}$ la r\'eunion des $B^{st}_{\mu}$ et $B^{inst}$ celle des $B^{inst}_{\mu}$. Notons enfin $B=B^{st}\cup B^{inst}$. R\'ealisons comme plus haut notre espace $PW^{\infty}_{ell}(\tilde{G})$ \`a l'aide de cette base $B$. Notons $PW^{\infty,st}_{ell}(\tilde{G})$, resp. $PW^{\infty,inst}_{ell}(\tilde{G})$ le sous-espace des $(\varphi_{\tilde{\pi}})_{\tilde{\pi}\in B}$ tels que $\varphi_{\tilde{\pi}}=0$ si $\tilde{\pi}\in B^{inst}$, resp. $\tilde{\pi}\in B^{st}$. Rappelons enfin que les espaces $I_{cusp}(\tilde{G}({\mathbb R}))$ et $SI_{cusp}(\tilde{G}({\mathbb R}))$ sont munis de topologies ([I] 5.3, qui reprenait Bouaziz et Renard).

\ass{Lemme}{Via l'isomorphisme de $I_{cusp}(\tilde{G}({\mathbb R}))$ sur $PW^{\infty}_{ell}(\tilde{G})$, les espaces $I_{cusp}^{st}(\tilde{G}({\mathbb R}))$ et $I_{cusp}^{inst}(\tilde{G}({\mathbb R}))$ s'identifient respectivement \`a $PW^{\infty,st}_{ell}(\tilde{G})$ et $PW^{\infty,inst}_{ell}(\tilde{G})$. La projection de $I_{cusp}(\tilde{G}({\mathbb R}))$ sur $PW^{\infty,st}_{ell}(\tilde{G})$ se quotiente en un hom\'eomorphisme de $SI_{cusp}(\tilde{G}({\mathbb R}))$ sur $PW^{\infty,st}_{ell}(\tilde{G})$.} 

Preuve. Notons $pw:I_{cusp}(\tilde{G}({\mathbb R}))\to PW^{\infty}_{ell}(\tilde{G})$ notre isomorphisme. Soit $\tilde{\pi}\in B$, consid\'erons un \'el\'ement $f\in I_{cusp}(\tilde{G}({\mathbb R}))$ tel que $pw(f)=(\varphi_{\tilde{\pi}'})_{\tilde{\pi}'\in B}$ v\'erifie $\varphi_{\tilde{\pi}'}=0$ pour $\tilde{\pi}'\not=\tilde{\pi}$. On va prouver que

(2) si $\tilde{\pi}\in B^{st}$, resp. $\tilde{\pi}\in B^{inst}$, alors $f\in I_{cusp}^{st}(\tilde{G}({\mathbb R}))$, resp. $f\in I_{cusp}^{inst}(\tilde{G}({\mathbb R}))$.

Soit $\phi$ la fonction sur $\mathfrak{A}_{\tilde{G}}$ (que l'on identifie \`a ${\cal A}_{\tilde{G}}$ via l'exponentielle) dont la transform\'ee de Fourier est la restriction de $\varphi_{\tilde{\pi}}$ \`a $i{\cal A}_{\tilde{G}}^*$. Elle est $C^{\infty}$ et \`a support compact.  On d\'efinit une fonction $f'$ sur $\tilde{G}({\mathbb R})$ par $f'(\gamma)=(\tilde{\pi},\tilde{\pi})_{ell}^{-1}\phi(H_{\tilde{G}}(\gamma))f[\tilde{\pi}](\gamma)$. C'est un \'el\'ement de $C_{c}^{\infty}(\tilde{G}({\mathbb R}))$. Pour tout $\gamma\in \tilde{G}({\mathbb R})$ fortement r\'egulier, on a l'\'egalit\'e
$$I^{\tilde{G}}(\gamma,f')=(\tilde{\pi},\tilde{\pi})_{ell}^{-1}\phi(H_{\tilde{G}}(\gamma))I^{\tilde{G}}(\gamma,f[\tilde{\pi}]).$$
Il en r\'esulte que ces int\'egrales orbitales sont nulles si $\gamma$ n'est pas elliptique et qu'elles ont les m\^emes propri\'et\'es de stabilit\'e que celles de $f[\tilde{\pi}]$. Donc l'image de $f'$ dans $I(\tilde{G}({\mathbb R}))$ appartient \`a $I_{cusp}^{st}(\tilde{G}({\mathbb R}))$ si $\tilde{\pi}\in B^{st}$, \`a $I_{cusp}^{inst}(\tilde{G}({\mathbb R}))$ si $\tilde{\pi}\in B^{inst}$.

Pour $\tilde{\pi}'\in B$ et $\lambda\in {\cal A}_{\tilde{G},{\mathbb C}}^*$, on a
$$\tilde{\pi}'_{\lambda}(f')=\int_{\tilde{G}({\mathbb R})}\tilde{\pi}'_{\lambda}(\gamma)f'(\gamma)\,d\gamma$$
$$=(\tilde{\pi},\tilde{\pi})_{ell}^{-1}\int_{\tilde{G}({\mathbb R})^1}\int_{{\cal A}_{\tilde{G}}}\tilde{\pi}'(\gamma^1)e^{<H,\lambda>}f[\tilde{\pi}](\gamma^1)\phi(H)\,dH\,d\gamma^1$$
$$=(\tilde{\pi},\tilde{\pi})_{ell}^{-1}\varphi_{\tilde{\pi}}(\lambda)\tilde{\pi}'(f[\tilde{\pi}]).$$
Il en r\'esulte que
$$I^{\tilde{G}}(\tilde{\pi}'_{\lambda},f')=\left\lbrace\begin{array}{cc}0,&\text{ si }\tilde{\pi}'\not=\tilde{\pi},\\ \varphi_{\tilde{\pi}}(\lambda),&\text{ si }\tilde{\pi}'=\tilde{\pi}.\\ \end{array}\right.$$
Autrement dit, l'image de $f'$ dans $PW^{\infty}_{ell}(\tilde{G})$ co\"{\i}ncide avec celle de $f$. Donc $f$ est l'image de $f'$ dans $I_{cusp}(\tilde{G}({\mathbb R}))$. L'assertion (2) en r\'esulte. 

Soit maintenant $f\in I_{cusp}(\tilde{G}({\mathbb R}))$ telle que $pw(f)\in PW^{\infty,st}_{ell}(\tilde{G})$. Ecrivons $B^{st}$ comme r\'eunion croissante de sous-ensembles finis $B_{i}$, pour $i\in {\mathbb N}$. Ecrivons $pw(f)=(\varphi_{\tilde{\pi}})_{\tilde{\pi}\in B}$. Pour tout $i$, soit $f_{i}\in I_{cusp}(\tilde{G}({\mathbb R}))$ l'\'el\'ement tel que, en posant $pw(f_{i})=(\varphi_{i,\tilde{\pi}})_{\tilde{\pi}\in B}$, on ait $\varphi_{i,\tilde{\pi}}=\varphi_{\tilde{\pi}}$ si $\tilde{\pi}\in B_{i}$, $\varphi_{i,\tilde{\pi}}=0$ sinon. D'apr\`es (2), $f_{i}$ appartient \`a $I_{cusp}^{st}(\tilde{G}({\mathbb R}))$. Parce que $pw$ est un hom\'eomorphisme, $f$ est la limite des $f_{i}$. Il est clair par d\'efinition que $I^{st}_{cusp}(\tilde{G}({\mathbb R})$ est un sous-espace ferm\'e. Il en r\'esulte que $f\in I_{cusp}^{st}(\tilde{G}({\mathbb R}))$. Donc $pw^{-1}(PW^{\infty,st}_{ell}(\tilde{G}))\subset I^{st}_{cusp}(\tilde{G}({\mathbb R}))$. On prouve de m\^eme que $pw^{-1}(PW^{\infty,inst}_{ell}(\tilde{G}))\subset I^{inst}_{cusp}(\tilde{G}({\mathbb R}))$. Puisque la somme de $PW^{\infty,st}_{ell}(\tilde{G})$ et de $PW^{\infty,inst}_{ell}(\tilde{G})$ est l'espace $PW^{\infty}_{ell}(\tilde{G})$ tout entier et que l'intersection de $ I^{st}_{cusp}(\tilde{G}({\mathbb R}))$ et $ I^{inst}_{cusp}(\tilde{G}({\mathbb R}))$ 
  est clairement r\'eduite \`a $0$, les inclusions pr\'ec\'edentes sont des \'egalit\'es. Cela prouve les premi\`eres assertions du lemme. La derni\`ere en r\'esulte imm\'ediatement. $\square$
  
  \bigskip
  
  \subsection{Un th\'eor\`eme de Paley-Wiener d\'ecrivant l'espace $SI(\tilde{G}({\mathbb R}))$}
  Consid\'erons l'espace 
  $$\oplus_{\tilde{L}\in {\cal L}(\tilde{M}_{0})}PW^{\infty}_{ell}(\tilde{L}).$$
  Par les constructions du paragraphe pr\'ec\'edent, on peut le d\'ecomposer en somme directe
  $$\left(\oplus_{\tilde{L}\in {\cal L}(\tilde{M}_{0})}PW^{\infty,st}_{ell}(\tilde{L})\right)\oplus\left(\oplus_{\tilde{L}\in {\cal L}(\tilde{M}_{0})}PW^{\infty,inst}_{ell}(\tilde{L})\right).$$
  Le groupe $W(\tilde{M}_0)$   agit sur l'espace total. On v\'erifie que cette action conserve les deux composantes ci-dessus. On peut donc d\'efinir $PW^{\infty,st}(\tilde{G})$ et $PW^{\infty,inst}(\tilde{G})$ comme les sous-espaces des invariants par $W(\tilde{M}_{0})$ dans chacune des composantes. Notons $\underline{pw}^{st}$ le compos\'e de l'isomorphisme $pw:I(\tilde{G}({\mathbb R}))\to PW^{\infty}(\tilde{G})$ et de la projection sur $PW^{\infty,st}(\tilde{G})$. Supposons prouv\'e que, pour tout $\tilde{L}\in {\cal L}(\tilde{M}_{0})$, tout $\tilde{\pi}\in D_{ell,0}^{st}(\tilde{L}({\mathbb R}))$ et tout $\lambda\in {\cal A}_{\tilde{L},{\mathbb C}}^*$, la distribution $f\mapsto I^{\tilde{L}}(\tilde{\pi}_{\lambda},f_{\tilde{L}})$ soit stable.  Alors $\underline{pw}^{st}$ se quotiente en un homomorphisme continu  $pw^{st}:SI(\tilde{G}({\mathbb R}))\to PW^{\infty,st}(\tilde{G})$.  Remarquons que l'hypoth\`ese que l'on vient de faire r\'esulte de l'hypoth\`ese plus simple que, pour tout $\tilde{\pi}\in D_{ell,0}^{st}(\tilde{G}({\mathbb R})) $, la distribution $f\mapsto I^{\tilde{G}}(\tilde{\pi},f)$ est stable. En effet,  la tensorisation par un \'el\'ement $\lambda\in {\cal A}_{\tilde{G},{\mathbb C}}^*$ respecte la stabilit\'e.  Comme toujours, on suppose par r\'ecurrence que les propri\'et\'es vraies pour $\tilde{G}$ le sont aussi pour les groupes tordus plus petits, donc pour les espaces de Levi. Donc,  pour tout $\tilde{L}\in {\cal L}(\tilde{M}_{0})$, tout $\tilde{\pi}\in D_{ell,0}^{st}(\tilde{L}({\mathbb R}))$ et tout $\lambda\in {\cal A}_{\tilde{L},{\mathbb C}}^*$, la distribution $f\mapsto I^{\tilde{L}}(\tilde{\pi}_{\lambda},f)$ sur $I(\tilde{L}({\mathbb R}))$ est stable. Mais l'induction conserve la stabilit\'e. Donc la distribution $f\mapsto I^{\tilde{L}}(\tilde{\pi}_{\lambda},f_{\tilde{L}})$ est stable.
  
  \ass{Th\'eor\`eme}{(i) Pour tout $\tilde{\pi}\in D_{ell,0}^{st}(\tilde{G}({\mathbb R})) $, la distribution $f\mapsto I^{\tilde{G}}(\tilde{\pi},f)$ est stable.
  
  (ii) L'application lin\'eaire $pw^{st}:SI(\tilde{G}({\mathbb R}))\to PW^{\infty,st}(\tilde{G})$ est un hom\'eomorphisme. }
  
  Dans le paragraphe suivant, on ram\`enera le th\'eor\`eme \`a une autre assertion qui sera prouv\'ee en 2.7.
  
  \bigskip
 
 \subsection{Un r\'esultat d'instabilit\'e}
 Notons
 $$sym^W:\oplus_{\tilde{L}\in {\cal L}(\tilde{M}_{0})}PW^{\infty}_{ell}(\tilde{L})\to PW^{\infty}(\tilde{G})$$
 l'application de sym\'etrisation $sym^W(\varphi)=\sum_{w\in W(\tilde{M}_{0})}w(\varphi)$. Soit $\tilde{L}\in {\cal L}(\tilde{M}_{0})$. On munit $D_{ell,0}(\tilde{L}({\mathbb R}))$ d'une base $B(\tilde{L})$ ayant les m\^emes propri\'et\'es que la base $B$ de 2.2. En particulier, $B(\tilde{L})$ est r\'eunion de $B^{st}(\tilde{L})$ et $B^{inst}(\tilde{L})$. On r\'ealise l'espace $PW^{\infty}_{ell}(\tilde{L})$ en utilisant cette base $B(\tilde{L})$. Fixons $\tilde{\pi}\in B(\tilde{L})$. On note $PW_{\tilde{\pi}}(\tilde{L})$ le sous-espace des $(\varphi_{\tilde{\pi}'})_{\tilde{\pi}'\in B(\tilde{L})}\in PW^{\infty}_{ell}(\tilde{L})$ tels que $\varphi_{\tilde{\pi}'}=0$ si $\tilde{\pi}'\not=\tilde{\pi}$. C'est un unique espace de fonctions de Paley-Wiener sur ${\cal A}_{\tilde{L},{\mathbb C}}^*$. 
 
 \ass{Proposition}{Pour tout $\tilde{L}\in {\cal L}(\tilde{M}_{0})$ et tout $\tilde{\pi}\in B^{inst}(\tilde{L})$, l'espace  $pw^{-1}\circ sym^W(PW_{\tilde{\pi}}(\tilde{L}))$ est inclus dans $I^{inst}(\tilde{G}({\mathbb R}))$.}
 
 Montrons que cette proposition entra\^{\i}ne le th\'eor\`eme. Elle entra\^{\i}ne
 
 (1) l'espace $pw^{-1}(PW^{\infty,inst}(\tilde{G}))$ est inclus dans $I^{inst}(\tilde{G}({\mathbb R}))$.
 
  En effet, le m\^eme raisonnement que dans la preuve du lemme 2.2 montre que tout \'el\'ement  $f\in pw^{-1}(PW^{\infty,inst}(\tilde{G}))$ est limite d'une suite d'\'el\'ements $f_{i}$ qui sont combinaisons lin\'eaires finies d'\'el\'ements de $pw^{-1}\circ sym^W(PW_{\tilde{\pi}}(\tilde{L})$, pour des couples $(\tilde{L},\tilde{\pi})$ v\'erifiant les hypoth\`eses de la proposition. Puisque $I^{inst}(\tilde{G}({\mathbb R}))$ est ferm\'e, cette proposition entra\^{\i}ne qu'un tel $f$ appartient \`a cet espace.
  
  On a mieux:
  
  (2) $pw^{-1}(PW^{\infty,inst}(\tilde{G}))$ est \'egal \`a $I^{inst}(\tilde{G}({\mathbb R}))$.

  Preuve. Soit $f\in I^{inst}(\tilde{G}({\mathbb R}))$. On peut la d\'ecomposer en $f=f'+f''$, o\`u $pw(f')\in PW^{\infty,st}(\tilde{G})$ et $pw(f'')\in PW^{\infty,inst}(\tilde{G})$. L'assertion (1) entra\^{\i}ne que $f''$ appartient \`a $ I^{inst}(\tilde{G}({\mathbb R}))$, donc $f'$ aussi. Il suffit de prouver que $f'=0$. En oubliant cela, on consid\`ere $f\in I^{inst}(\tilde{G}({\mathbb R}))$ tel que $pw(f)\in PW^{\infty,st}(\tilde{G}({\mathbb R}))$ et on veut prouver que $f=0$. Par r\'ecurrence, on peut admettre l'assertion (i) du th\'eor\`eme pour tout espace de Levi propre $\tilde{L}$ de $\tilde{G}$. Soit $\tilde{L}$ un tel espace de Levi, soit $\tilde{\pi}\in D^{st}_{ell,0}(\tilde{L}({\mathbb R}))$ et soit $\lambda\in {\cal A}_{\tilde{L},{\mathbb C}}^*$. Comme on l'a expliqu\'e dans le paragraphe pr\'ec\'edent, la distribution $f'\mapsto I^{\tilde{L}}(\tilde{\pi}_{\lambda},f'_{\tilde{L}})$ est stable. Elle annule donc $f$. Il en r\'esulte que $pw(f)$ n'a de composantes non nulles que dans les $PW_{\tilde{\pi}}(\tilde{G})$ pour $\tilde{\pi}\in B^{st}(\tilde{G})$. Autrement dit, $pw(f)$ appartient au sous-espace $PW^{\infty,st}_{ell}(\tilde{G}({\mathbb R}))$. A fortiori, $f$ est cuspidale. Mais alors, le lemme  2.2 entra\^{\i}ne que $f=0$. $\square$
 
 Soit $\tilde{\pi}\in D_{ell,0}^{st}(\tilde{G}({\mathbb R}))$. Par d\'efinition de l'application $pw$, la distribution $f\mapsto I^{\tilde{G}}(\tilde{\pi},f)$ annule l'espace   $pw^{-1}( PW^{\infty,inst}(\tilde{G}))$. D'apr\`es (2), elle annule $I^{inst}(\tilde{G}({\mathbb R}))$, ce qui signifie qu'elle est stable. Cela prouve l'assertion (i) du th\'eor\`eme. L'assertion (2) entra\^{\i}ne aussi que l'application $\underline{pw}^{st}$ a pour noyau $I^{inst}(\tilde{G}({\mathbb R}))$, donc que $pw^{st}$ est injective. Elle est surjective puisque $pw$ l'est. Enfin, $pw^{st}$ admet une section continue, \`a savoir l'application compos\'ee
 $$PW^{\infty,st}(\tilde{G}({\mathbb R}))\to PW^{\infty}(\tilde{G}({\mathbb R}))\stackrel{pw^{-1}}{\to} I(\tilde{G}({\mathbb R}))\to SI(\tilde{G}({\mathbb R})).$$
 Donc $pw^{st}$ est un hom\'eomorphisme. Cela prouve le th\'eor\`eme.
 
 Commen\c{c}ons la preuve de la proposition. Si $\tilde{L}=\tilde{G}$, l'assertion r\'esulte du lemme 2.2. On suppose donc $\tilde{L}$ propre et on raisonne par r\'ecurrence sur $a_{\tilde{L}}=dim({\cal A}_{\tilde{L}})$ en supposant la proposition d\'emontr\'ee pour tout  couple $(\tilde{L}',\tilde{\pi}')$ analogue \`a $(\tilde{L},\tilde{\pi})$ tel que $a_{\tilde{L}'}<a_{\tilde{L}}$. On note simplement ${\cal F}$ l'espace des fonctions de Paley-Wiener sur ${\cal A}_{\tilde{L},{\mathbb C}}^*$ et on identifie $PW_{\tilde{\pi}}(\tilde{L})$ \`a ${\cal F}$. Pour $\varphi\in {\cal F}$, on pose $f_{\varphi}=pw^{-1}\circ sym^W(\varphi)$. Soit $\tilde{M}\in {\cal L}(\tilde{M}_{0})$ un espace de Levi propre de $\tilde{G}$. En notant $pw^{\tilde{M}}$ l'analogue de $pw$ pour l'espace $\tilde{M}$, le terme $pw^{\tilde{M}}(f_{\varphi,\tilde{M}})$ se d\'eduit ais\'ement de $pw(f_{\varphi})$.  On voit qu'il n'a de composantes non nulles que dans les sous-espaces $PW^{\infty,inst}_{ell,0}(\tilde{L}')$ pour des $\tilde{L}'\subset \tilde{M}$ conjugu\'es \`a $\tilde{L}$. En appliquant le th\'eor\`eme par r\'ecurrence \`a $\tilde{M}$, on voit que $f_{\varphi,\tilde{M}}$ appartient \`a $I^{inst}(\tilde{M}({\mathbb R}))$. Cela \'etant vrai pour tout $\tilde{M}$ propre, les int\'egrales orbitales stables de $f_{\varphi}$ sont donc nulles sur tout \'el\'ement fortement r\'egulier non elliptique de $\tilde{G}({\mathbb R})$. Donc l'image de $f_{\varphi}$ dans $SI(\tilde{G}({\mathbb R}))$ appartient au sous-espace $SI_{cusp}(\tilde{G}({\mathbb R}))$. D'apr\`es le lemme 2.2, il existe un unique $f_{\varphi}^{st}\in I^{st}_{cusp}(\tilde{G}({\mathbb R}))$ qui a m\^eme image que $f_{\varphi}$ dans $SI(\tilde{G}({\mathbb R}))$. Il s'agit de prouver que $f_{\varphi}^{st}$ est nul. Toujours d'apr\`es le lemme 2.2, on peut fixer $\tilde{\sigma}_{0}\in B^{st}(\tilde{G})$, $\lambda_{0}\in i{\cal A}_{\tilde{G}}^*$, poser $\tilde{\sigma}=\tilde{\sigma_{0,\lambda_{0}}}$ et prouver que $I^{\tilde{G}}(\tilde{\sigma},f_{\varphi}^{st})=0$.  Pour $\varphi\in {\cal F}$, posons $\ell(\varphi) = I^{\tilde{G}}(\tilde{\sigma},f_{\varphi}^{st})$. Etudions l'application lin\'eaire $\ell$ sur ${\cal F}$. Rappelons que $\mathfrak{Z}(G)$ agit naturellement sur ${\cal F}$: pour $z\in \mathfrak{Z}(G)$, $\varphi\in {\cal F}$ et $\lambda\in {\cal A}_{\tilde{L},{\mathbb C}}^*$, on a $(z\varphi)(\lambda)=z(\mu(\tilde{\pi})+\lambda)\varphi(\lambda)$. Parce que $f_{\varphi}^{st}$ est uniquement d\'etermin\'e, l'application $\varphi\mapsto f_{\varphi}^{st}$ est \'equivariante pour les actions de $\mathfrak{Z}(G)$. Pour $z\in \mathfrak{Z}(G)$, on a donc 
 $$\ell(z\varphi)=I^{\tilde{G}}(\tilde{\sigma},z(f_{\varphi}^{st}))=I^{\tilde{G}}(z\tilde{\sigma},f_{\varphi}^{st})=z(\mu(\tilde{\sigma}))I^{\tilde{G}}(\tilde{\sigma},f_{\varphi}^{st})=z(\mu(\tilde{\sigma}))\ell(\varphi).$$
 Notons $J$ l'id\'eal des \'el\'ements $z\in \mathfrak{Z}(G)$ tels que $z(\mu(\tilde{\sigma}))=0$. Alors $\ell$ annule $J{\cal F}$. Les lemmes des deux paragraphes suivants nous permettrons de pr\'eciser cet ensemble $J{\cal F}$.
  
 \bigskip
 
 \subsection{Un lemme sur les fonctions de Paley-Wiener}
 
 On conserve les notations pr\'ec\'edentes. Rappelons que $\mu(\tilde{\pi})$ est une $W^L$-orbite dans $\mathfrak{h}^{\tilde{L},*}$ et que $\mu(\tilde{\sigma})$ est une $W$-orbite dans $\mathfrak{h}^*$. Si l'intersection $(\mu(\tilde{\pi})+{\cal A}_{\tilde{L},{\mathbb C}}^*)\cap \mu(\tilde{\sigma})$ est non vide, on note $(\lambda_{i})_{i=1,...,m}$ la famille finie d'\'el\'ements de ${\cal A}_{\tilde{L},{\mathbb C}}^*$ tels que cette intersection soit la r\'eunion des $\mu(\tilde{\pi})+\lambda_{i}$ pour $i=1,...,m$. 
 
 \ass{Lemme}{(i) Si $(\mu(\tilde{\pi})+{\cal A}_{\tilde{L},{\mathbb C}}^*)\cap \mu(\tilde{\sigma})=\emptyset$, on a l'\'egalit\'e $J{\cal F}={\cal F}$.
 
 (ii) Si cette intersection est non vide, il existe un entier $N\geq1$ tel que $J{\cal F}$ contienne  toute fonction $\varphi\in {\cal F}$ qui s'annule \`a l'ordre au moins $N$ en chaque point $\lambda_{i}$ pour $i=1,...,m$.}
 
 Preuve de (i). On va montrer qu'il existe $z\in J$ tel que sa restriction \`a $\mu(\tilde{\pi})+{\cal A}_{\tilde{L},{\mathbb C}}^*$  soit constante de valeur $1$. L'assertion en r\'esulte puisqu'alors $\varphi=z\varphi$ pour tout $\varphi\in {\cal F}$. Notons $Y$ la projection de $\mu(\tilde{\sigma})$ sur $\mathfrak{h}^{\tilde{L},*}$ L'hypoth\`ese signifie que $\mu(\tilde{\pi})\cap Y=\emptyset$. Pour tout $x\in \mu(\tilde{\pi})$, on peut alors trouver un polyn\^ome  $q_{x}$ de degr\'e $1$ sur $\mathfrak{h}^{\tilde{L},*}$ tel que $q_{x}(x)=0$ et $q_{x}(y)\not=0$ pour tout $y\in Y$. Posons $q=\prod_{x\in \mu(\tilde{\pi})}q_{x}$. Consid\'erons $q$ comme un polyn\^ome sur $\mathfrak{h}^*$ via la projection $\mathfrak{h}^*\to \mathfrak{h}^{\tilde{L},*}$. D\'efinissons un polyn\^ome $z_{0}$ par $z_{0}(\nu)=\prod_{w\in W}q(w\nu)$. Il est invariant par $W$ donc appartient \`a $\mathfrak{Z}(G)$. Il s'annule sur $\mu(\tilde{\pi})+{\cal A}_{\tilde{L},{\mathbb C}}^*$ car le polyn\^ome $q$ lui-m\^eme s'annule sur cet ensemble. Pour $\nu\in \mu(\tilde{\sigma})$ et $w\in W$, $w\nu$ appartient aussi \`a $\mu(\tilde{\sigma})$ et se projette sur $\mathfrak{h}^{\tilde{L},*}$ en un point $y\in Y$. On a alors $q(w\nu)=q(y)\not=0$. Donc $z_{0}(\nu)\not=0$, ce que l'on peut noter $z_{0}(\mu(\tilde{\sigma}))\not=0$ puisque $z_{0}$ est invariant par $W$. Posons $z= 1-z_{0}(\mu(\tilde{\sigma}))^{-1}z_{0}$. Cet \'el\'ement r\'epond \`a la question.
 
 Preuve de (ii). Notons $ {\cal P}$ l'espace des polyn\^omes sur $\mathfrak{h}^*$. On montre d'abord
 
 (1) il existe un entier $N\geq1$ tel que tout \'el\'ement de $ {\cal P}$ qui s'annule \`a l'ordre au moins $N$ en tout point de $\mu(\tilde{\sigma})$ appartient \`a $J{\cal P}$.
 
 Rappelons que $\mathfrak{Z}(G)$ est le sous-espace des invariants par $W$ dans ${\cal P}$. Comme on sait, on peut fixer un sous-ensemble $(h_{w})_{w\in W}$  de ${\cal P}$ tel que tout \'el\'ement $Q\in {\cal P}$ s'\'ecrive de fa\c{c}on unique $Q=\sum_{w\in W}h_{w}Q_{w}$, avec $Q_{w}\in \mathfrak{Z}(G)$. On peut de plus fixer un \'el\'ement non nul $D\in {\cal P}$ et une matrice $(D_{w,w'})_{w,w'\in W}$ d'\'el\'ements de ${\cal P}$ de sorte que, pour tout $w\in W$ et tout $\lambda\in \mathfrak{h}^*$, on ait l'\'egalit\'e
 $$D(\lambda)Q_{w}(\lambda)=\sum_{w'\in W}D_{w,w'}(\lambda)Q(w'\lambda).$$
 On renvoie pour cela \`a [BR] preuve de la proposition 3.3. Notons $N-1$ le maximum des ordres d'annulation de $D$ aux diff\'erents points de $\mu(\tilde{\sigma})$. Il r\'esulte de la formule ci-dessus que, si $Q$ s'annule \`a l'ordre $N$ en tout point de $\mu(\tilde{\sigma})$, alors chaque $Q_{w}$ s'annule sur $\mu(\tilde{\sigma})$. Puisque $Q_{w}\in \mathfrak{Z}(G)$, cela signifie que $Q_{w}\in J$. Mais alors $Q=\sum_{w\in W}h_{w}Q_{w}$ appartient \`a $J{\cal P}$. Cela prouve (1). 
 
 Notons ${\cal H}$ l'espace des fonctions de Paley-Wiener sur $\mathfrak{h}^*$. Soit $Y$ un ensemble fini d'\'el\'ements de $\mathfrak{h}^*$ et soit ${\bf n}=(n_{y})_{y\in Y}$ une famille d'entiers naturels. Notons ${\cal H}_{{\bf n}}$, resp. $ {\cal P}_{{\bf n}}$, l'espace des \'el\'ements de ${\cal H}$, resp. ${\cal P}$, qui s'annulent en tout point $y\in Y$ \`a l'ordre au moins $n_{y}$. On va montrer
 
 (2) on a l'\'egalit\'e ${\cal H}_{{\bf n}}={\cal P}_{{\bf n}}{\cal H}$.
 
 C'est trivial si tous les $n_{y}$ sont nuls. Supposons qu'il existe un $y$ pour lequel $n_{y}>0$, fixons-en un que l'on note $y_{0}$. Notons ${\bf n}'=(n'_{y})_{y\in Y}$, o\`u $n'_{y}=n_{y}$ si $y\not=y_{0}$ et $n'_{y_{0}}=n_{y_{0}}-1$. Notons ${\cal P}_{y_{0}}$ l'espace des polyn\^omes qui s'annulent en $y_{0}$. L'assertion (2) r\'esulte par r\'ecurrence de l'assertion
 
 (3) on a l'\'egalit\'e ${\cal H}_{{\bf n}}={\cal P}_{y_{0}}{\cal H}_{{\bf n}'}$.
 
On ne perd rien \`a supposer $y_{0}=0$. On peut fixer un syst\`eme de coordonn\'ees sur $\mathfrak{h}_{\mathbb R}^*$ de sorte que les coordonn\'ees $y_{1},...,y_{n}$ de $y$ soient toutes non nulles pour  $y\in Y$, $y\not=0$. On \'ecrit tout \'el\'ement de $\mathfrak{h}^*$ sous la forme $\nu=(\nu_{1},...,\nu_{n})$. Fixons une fonction de Paley-Wiener $h$ sur ${\mathbb C}$ telle que $h(0)=1$ et $h$ s'annule \`a l'ordre au moins  $n_{y}$ en tout $y_{i}$, pour $y\in Y-\{0\}$ et $i=1,...,n$. Soit $\varphi\in {\cal H}_{{\bf n}}$. Pour $ i=1,...,n$, introduisons la fonction $\phi_{i}$ sur $\mathfrak{h}^*$ d\'efinie par
$$\phi_{i}(\nu)=\varphi(0,...,0,\nu_{i},...,\nu_{n}) \prod_{j=1,...,i-1}h(\nu_{j}).$$
Elle est de Paley-Wiener. Elle s'annule \`a l'ordre au moins $n_{0}$ en $0$ car il en est ainsi du premier facteur.  Le deuxi\`eme facteur s'annule \`a l'ordre au moins $(i-1)n_{y}$ en tout point $y\in Y-\{0\}$. Si $i\geq2$, $\phi_{i}$ s'annule donc \`a l'ordre au moins $n_{y}$ en un tel point. Si $i=1$, on a simplement $\phi_{1}=\varphi$ et cette propri\'et\'e est encore v\'erifi\'ee par hypoth\`ese. Donc $\phi_{i}\in {\cal H}_{{\bf n}}$. Par construction, pour $i=1,...,n-1$, la fonction $\phi_{i}-\phi_{i+1}$ s'annule sur l'hyperplan $\nu_{i}=0$. En posant $\phi_{n+1}=0$, il en est de m\^eme de la fonction $\phi_{n}-\phi_{n+1}$ car $\varphi$ s'annule en $0$.  On sait qu'une fonction de Paley-Wiener qui s'annule sur un tel hyperplan $\nu_{i}=0$ est divisible dans ${\cal H}$ par $\nu_{i}$. Il existe donc pour tout $i=1,...,n$ un \'el\'ement $\varphi_{i}\in {\cal H}$ tel que $\nu_{i}\varphi_{i}(\nu)=\phi_{i}(\nu)-\phi_{i+1}(\nu)$. Il est clair que $\varphi_{i}$ appartient \`a ${\cal H}_{{\bf n}'}$. On a d'autre part l'\'egalit\'e
$$\sum_{i=1,...,n}\nu_{i}\varphi_{i}(\nu)=\phi_{1}(\nu)-\phi_{n+1}(\nu)=\varphi(\nu),$$ 
ce qui r\'ealise $\varphi$ comme un \'el\'ement de  ${\cal P}_{0}{\cal H}_{{\bf n}'}$. Cela prouve (3) et (2).
 
En appliquant (2) \`a l'ensemble $Y=\mu(\tilde{\sigma})$ et \`a la famille ${\bf n}=(n_{y})_{y\in Y}$ telle que $n_{y}$ soit  pour tout $y$ un entier v\'erifiant la relation (1), on obtient

(4) il existe un entier $N$ tel que toute fonction $\varphi\in {\cal H}$ qui s'annule \`a l'ordre au moins $N$ en tout point de $\mu(\tilde{\sigma})$ appartienne \`a $J{\cal H}$.

Fixons un tel entier $N$. Notons maintenant $Y$ l'ensemble des projections dans $\mathfrak{h}^{\tilde{L},*}$ des \'el\'ements de $\mu(\tilde{\sigma})$. Il contient $\mu(\tilde{\pi})$ par hypoth\`ese. Fixons une fonction de Paley-Wiener $h$ sur $\mathfrak{h}^{\tilde{L},*}$ qui vaut $1$ en tout point de $\mu(\tilde{\pi})$ et qui s'annule \`a l'ordre au moins $N$ en tout point de $Y-\mu(\tilde{\pi})$. Soit $\varphi\in {\cal F}$ s'annulant \`a l'ordre au moins $N$ en tout $\lambda_{i}$, $i=1,...,m$. D\'efinissons une fonction $\varphi'$ sur $\mathfrak{h}^*$ par $\varphi'(\lambda^{\tilde{L}}+\lambda_{\tilde{L}})=h(\lambda^{\tilde{L}})\varphi(\lambda_{\tilde{L}})$, pour $\lambda^{\tilde{L}}\in \mathfrak{h}^{\tilde{L},*}$ et $\lambda_{\tilde{L}}\in {\cal A}_{\tilde{L},{\mathbb C}}^*$. C'est un \'el\'ement de ${\cal H}$ qui s'annule \`a l'ordre au moins $N$ en tout point de $\mu(\tilde{\sigma})$. Appliquant (4), on peut l'\'ecrire
$$\varphi'=\sum_{j=1,...,k}z_{j}\varphi'_{j},$$
o\`u les $z_{j}$ appartiennent \`a $J$ et les $\varphi'_{j}$ appartiennent \`a ${\cal H}$. On fixe un \'el\'ement $ x_{0}\in \mu(\tilde{\pi})$ et on d\'efinit pour tout $j$ une fonction $\varphi_{j}$ sur ${\cal A}_{\tilde{L},{\mathbb C}}^*$ par $\varphi_{j}(\lambda)=\varphi'_{j}(x_{0}+\lambda)$. C'est un \'el\'ement de ${\cal F}$. Pour $\lambda\in {\cal A}_{\tilde{L},{\mathbb C}}^*$, on a l'\'egalit\'e
$$\varphi(\lambda)=\varphi'(x_{0}+\lambda)=\sum_{j=1,...,k}z_{j}(\mu(\tilde{\pi})+\lambda)\varphi_{j}(\lambda).$$
Autrement dit, $\varphi=\sum_{j=1,...,k}z_{j}\varphi_{j}$ et $\varphi$ appartient \`a $J{\cal F}$. Cela ach\`eve la preuve. $\square$

Puisque notre forme lin\'eaire $\ell$ annule $J{\cal F}$, il s'ensuit du (i) du lemme que $\ell=0$ si $(\mu(\tilde{\pi})+{\cal A}_{\tilde{L},{\mathbb C}}^*)\cap \mu(\tilde{\sigma})=\emptyset$. Puisqu'on veut justement prouver que $\ell$ est nulle, on a termin\'e dans ce cas. On suppose dans la suite que $(\mu(\tilde{\pi})+{\cal A}_{\tilde{L},{\mathbb C}}^*)\cap \mu(\tilde{\sigma})\not=\emptyset$. En cons\'equence du (ii) du lemme, on a

(5) pour $i=1,...,m$, il existe un op\'erateur diff\'erentiel $D_{i}$ sur ${\cal A}_{\tilde{L},{\mathbb C}}^*$ tel que
$$\ell(\varphi)=\sum_{i=1,...,m}(D_{i}\varphi)(\lambda_{i})$$
pour tout $\varphi\in {\cal F}$. 

\bigskip

\subsection{Fonctions $f_{\varphi}$ \`a support assez r\'egulier}
Soit $\varphi\in {\cal F}$. Introduisons la fonction $\phi$ sur ${\cal A}_{\tilde{L}}$ dont la transform\'ee de Fourier est la restriction de $\varphi$ \`a $i{\cal A}_{\tilde{L}}^*$.  On introduit la fonction $f_{\varphi}^{\tilde{L}}$ sur $\tilde{L}({\mathbb R})$ d\'efinie par $f_{\varphi}^{\tilde{L}}(\gamma)=\phi(H_{\tilde{L}}(\gamma))f[\tilde{\pi}](\gamma)$ pour tout $\gamma\in \tilde{L}({\mathbb R})$. On a vu dans la preuve du lemme 2.2 que cette fonction appartient \`a $C_{cusp}^{\infty}(\tilde{L}({\mathbb R}))$ et que, en notant encore $f_{\varphi}^{\tilde{L}}$ son image dans $I(\tilde{L}({\mathbb R}))$, $pw^{\tilde{L}}(f_{\varphi}^{\tilde{L}})$ est l'\'el\'ement de $PW^{\infty}_{ell}(\tilde{L})$ dont les composantes sont nulles pour $\tilde{\pi}'\in B(\tilde{L})-\{\tilde{\pi}\}$ et dont la composante sur $\tilde{\pi}$ est  $\varphi$. Introduisons l'espace
$$(1) \qquad \oplus_{\tilde{M}\in {\cal L}(\tilde{M}_{0}); a_{\tilde{M}}=a_{\tilde{L}}}I_{cusp}(\tilde{M}({\mathbb R})).$$
  Le groupe $W(\tilde{M}_{0})$ y agit naturellement  et on a montr\'e en [I] 4.2 que le sous-espace des invariants \'etait le  gradu\'e $Gr^{a_{\tilde{L}}}I(\tilde{G}({\mathbb R}))$ d'ordre $a_{\tilde{L}}$ d'une certaine filtration $({\cal F}^nI(\tilde{G}({\mathbb R})))_{n=a_{\tilde{M}_{0}},...,a_{\tilde{G}}}$ de $I(\tilde{G}({\mathbb R}))$ (la notation ${\cal F}$ pour cette filtration n'a rien \`a voir avec notre espace de Paley-Wiener).  On peut consid\'erer $f_{\varphi}^{\tilde{L}}$ comme un \'el\'ement de  l'espace (1). Posons
  $$f_{\varphi}^{a_{\tilde{L}}}=\sum_{w\in W(\tilde{M}_{0})}w(f_{\varphi}^{\tilde{L}}).$$
  On voit alors que $f_{\varphi}$ est un \'el\'ement du terme ${\cal F}^{a_{\tilde{L}}}I(\tilde{G}({\mathbb R}))$ de notre filtration et que son image dans $Gr^{a_{\tilde{L}}}I(\tilde{G}({\mathbb R}))$ est $f_{\varphi}^{a_{\tilde{L}}}$. Donc, pour un espace de Levi $\tilde{M}\in {\cal L}(\tilde{M}_{0})$ et pour un \'el\'ement  $\gamma\in \tilde{M}({\mathbb R})$ elliptique dans $\tilde{M}({\mathbb R})$ et fortement r\'egulier dans $\tilde{G}({\mathbb R})$, on a les \'egalit\'es suivantes
  
  (2) si $a_{\tilde{M}}> a_{\tilde{L}}$ ou si $a_{\tilde{M}}=a_{\tilde{L}}$ et $\tilde{M}$ n'est pas conjugu\'e \`a $\tilde{L}$, $I^{\tilde{G}}(\gamma,f_{\varphi})=0$;
  
  (3) si $\tilde{M}=\tilde{L}$,
   $$I^{\tilde{G}}(\gamma,f_{\varphi})=\sum_{w\in W(\tilde{L})}I^{\tilde{L}}(w(\gamma),f_{\varphi}^{\tilde{L}})=\sum_{w\in W(\tilde{L})}\phi(w(H_{\tilde{L}}(\gamma)))I^{\tilde{L}}(w\gamma,f[\tilde{\pi}]),$$
   o\`u $W(\tilde{L})=Norm_{\tilde{G}({\mathbb R})}(\tilde{L})/L({\mathbb R})$.
  
  Supposons que 
  
  (4) le support de $f_{\varphi}^{\tilde{L}}$ soit form\'e d'\'el\'ements $\gamma\in \tilde{L}({\mathbb R})$ qui sont $\tilde{G}$-\'equisinguliers, c'est-\`a-dire  tels que $L_{\gamma}=G_{\gamma}$. 
  
  Dans ce cas, on peut trouver un voisinage compact $U_{1}$ dans $\tilde{G}({\mathbb R})$ de ce support  et un voisinage compact $U_{2}$ de $U_{1}$ tels que tout \'el\'ement de $U_{2}$ soit conjugu\'e par $G({\mathbb R})$ \`a un \'el\'ement de $\tilde{L}({\mathbb R})$. On peut trouver une fonction $h$ sur $\tilde{G}({\mathbb R})$ qui est invariante par conjugaison par $G({\mathbb R})$, qui vaut $1$ sur $U_{1}$ et qui vaut $0$ sur tout \'el\'ement $\gamma\in \tilde{G}({\mathbb R})$ qui n'est pas conjugu\'e \`a un \'el\'ement de $U_{2}$. Posons $g_{\varphi}=hf_{\varphi}$. Cette fonction v\'erifie des propri\'et\'es analogues \`a (2) et (3). La propri\'et\'e de $h$ entra\^{\i}ne que l'on a aussi $I^{\tilde{G}}(\gamma,g_{\varphi})=0$ pour $\tilde{M}$ et $\gamma$ comme plus haut, si $a_{\tilde{M}}< a_{\tilde{L}}$. En d\'efinitive, pour un espace de Levi $\tilde{M}\in {\cal L}(\tilde{M}_{0})$ et pour un \'el\'ement  $\gamma\in \tilde{M}({\mathbb R})$ elliptique dans $\tilde{M}({\mathbb R})$ et fortement r\'egulier dans $\tilde{G}({\mathbb R})$, on a les \'egalit\'es suivantes
  
  (5) si $\tilde{M}$ n'est pas conjugu\'e \`a $\tilde{L}$, $I^{\tilde{G}}(\gamma,g_{\varphi})=0$;
 
 (6) si $\tilde{M}=\tilde{L}$, 
 $$I^{\tilde{G}}(\gamma,g_{\varphi})=\sum_{w\in W(\tilde{L})}\phi(w(H_{\tilde{L}}(\gamma)))I^{\tilde{L}}(w\gamma,f[\tilde{\pi}]).$$ 
 
  De ces formules et  de l'instabilit\'e de $\tilde{\pi}$ r\'esulte que $g_{\varphi}\in I^{inst}(\tilde{G}({\mathbb R}))$. De la comparaison de (2) et (3) d'une part, de (5) et (6) d'autre part, r\'esulte que, pour tout $n\geq a_{\tilde{L}}$, les composantes dans
  $$(7) \qquad (\oplus_{\tilde{M}\in {\cal L}(\tilde{M}_{0}); a_{\tilde{M}}=n}PW^{\infty}_{ell}(\tilde{M}))^{W(\tilde{M}_{0})}$$
  de $pw(g_{\varphi})$ et de $pw(f_{\varphi})$ sont \'egales. Elles sont donc nulles si $n> a_{\tilde{L}}$ et \'egales \`a $sym^W(\varphi)$ si $n=a_{\tilde{L}}$. Supposons $a_{\tilde{G}}< n<a_{\tilde{L}}$. On peut admettre le th\'eor\`eme 2.3 par r\'ecurrence pour les $\tilde{M}$ tels que $a_{\tilde{M}}=n$. C'est-\`a-dire que $\tilde{\pi}'$ est stable pour tout $\tilde{\pi}'\in B^{st}(\tilde{M})$. L'instabilit\'e de $g_{\varphi}$ implique alors que la composante de $pw(g_{\varphi})$ dans l'espace (7) appartient au sous-espace
  $$ (\oplus_{\tilde{M}\in {\cal L}(\tilde{M}_{0});a_{\tilde{M}}=n}PW^{\infty,inst}_{ell}(\tilde{M}))^{W(\tilde{M}_{0})}.$$
On a pos\'e une hypoth\`ese de r\'ecurrence au d\'ebut de la preuve de la proposition 2.4. Elle nous dit que l'image par $pw^{-1}$ de l'espace ci-dessus est inclus dans $I^{inst}(\tilde{G}({\mathbb R})$. Reste le cas o\`u $n=a_{\tilde{G}}$. L'espace (7) se r\'eduit \`a $PW^{\infty}_{ell}(\tilde{G})$, que l'on peut d\'ecomposer en $PW^{\infty,st}_{ell}(\tilde{G})\oplus PW^{\infty,inst}_{ell}(\tilde{G})$. On note $pw^{st}_{ell}(g_{\varphi})$ et $pw^{inst}_{ell}(g_{\varphi})$ les  projections de $pw(g_{\varphi})$ dans chacun de ces sous-espaces. Le lemme 2.2 nous dit que $pw^{st}_{ell}(g_{\varphi})$ appartient \`a $I^{st}_{cusp}(\tilde{G}({\mathbb R}))$ tandis que $pw^{inst}_{ell}(g_{\varphi})$ appartient \`a $I^{inst}_{cusp}(\tilde{G}({\mathbb R})$. Posons $g_{\varphi}^{st}=pw^{-1}(pw^{st}_{ell}(g_{\varphi}))$. Cela prouve que
$$g_{\varphi}\in f_{\varphi}+g_{\varphi}^{st}+I^{inst}(\tilde{G}({\mathbb R})).$$
Puisque $g_{\varphi}$ appartient elle-m\^eme \`a $I^{inst}(\tilde{G}({\mathbb R}))$, on obtient
$$f_{\varphi}+g_{\varphi}^{st}\in I^{inst}(\tilde{G}({\mathbb R})).$$
En se rappelant la d\'efinition de $f_{\varphi}^{st}$, on obtient $f_{\varphi}^{st}=-g_{\varphi}^{st}$. Alors 
$$\ell(\varphi)=I^{\tilde{G}}(\tilde{\sigma},f_{\varphi}^{st})=-I^{\tilde{G}}(\tilde{\sigma},g_{\varphi}^{st})=-I^{\tilde{G}}(\tilde{\sigma},g_{\varphi}),$$
la derni\`ere \'egalit\'e r\'esultant de la d\'efinition de $g_{\varphi}^{st}$. Comme on sait, le caract\`ere de $\tilde{\sigma}$ est donn\'e par une fonction localement int\'egrable $\Theta_{\tilde{\sigma}}$ sur $\tilde{G}({\mathbb R})$. Notons $\tilde{T}$ un tore tordu maximal elliptique de $\tilde{L}$ (l'existence de $\tilde{\pi}\in D_{ell,0}(\tilde{L}({\mathbb R}))$ implique l'existence d'un tel tore). Parce que l'on est dans une situation \`a torsion int\'erieure, on sait qu'il n'y en a qu'un, \`a conjugaison pr\`es  par $L({\mathbb R})$ et on peut supposer que tout \'el\'ement de $W(\tilde{L})$ conserve ce tore tordu. L'\'egalit\'e pr\'ec\'edente et les relations (5) et (6) entra\^{\i}nent l'\'egalit\'e
$$\ell(\varphi)=-c\sum_{w\in W(\tilde{L})}\int_{\tilde{T}({\mathbb R})}D^{\tilde{G}}(\gamma)^{1/2}\Theta_{\tilde{\sigma}}(\gamma)\phi(w(H_{\tilde{L}}(\gamma)))I^{\tilde{L}}(w\gamma,f[\tilde{\pi}])\,d\gamma,$$
la constante $c>0$ ne d\'ependant que des mesures.  Puisque $\Theta_{\tilde{\sigma}}$ est invariant par $W(\tilde{M}_{0})$, cette expression se simplifie d'ailleurs en
$$\ell(\varphi)=-c\vert W(\tilde{L})\vert \int_{\tilde{T}({\mathbb R})}D^{\tilde{G}}(\gamma)^{1/2}\Theta_{\tilde{\sigma}}(\gamma)\phi(H_{\tilde{L}}(\gamma))I^{\tilde{L}}(\gamma,f[\tilde{\pi}])\,d\gamma.$$
En se rappelant la d\'ecomposition $\tilde{L}({\mathbb R})={\cal A}_{\tilde{L}}\times \tilde{L}({\mathbb R})^1$ et en posant $\tilde{T}({\mathbb R})^1=\tilde{T}({\mathbb R})\cap \tilde{L}({\mathbb R})^1$, on obtient
$$(8) \qquad \ell(\varphi)=-c\vert W(\tilde{L})\vert \int_{\tilde{T}({\mathbb R})^1}\int_{{\cal A}_{\tilde{L}}}D^{\tilde{G}}(exp(H)\gamma^1)^{1/2}$$
$$\Theta_{\tilde{\sigma}}(exp(H)\gamma^1)\phi(H )I^{\tilde{L}}(\gamma^1,f[\tilde{\pi}])\,dH\,d\gamma^1.$$

\bigskip

\subsection{Utilisation de la propri\'et\'e: une repr\'esentation elliptique est supertemp\'er\'ee}
Fixons $\varphi\in {\cal F}$. Pour $X\in {\cal A}_{\tilde{L}}$, d\'efinissons la fonction $\varphi_{X}\in {\cal F}$ par $\varphi_{X}(\lambda)=e^{<X,\lambda>}\varphi(\lambda)$ pour tout $\lambda\in {\cal A}_{\tilde{L},{\mathbb C}}^*$. Si $\phi$, resp. $\phi_{X}$, est la fonction sur ${\cal A}_{\tilde{L}}$ dont la transform\'ee de Fourier est la restriction de $\varphi$, resp. $\varphi_{X}$, \`a $i{\cal A}_{\tilde{L}}^*$, on a $\phi_{X}(H)=\phi(H-X)$. Donc, $\varphi$ \'etant fix\'e, la fonction $f_{\varphi_{X}}$ v\'erifie l'hypoth\`ese (4) du paragraphe pr\'ec\'edent pour tout $X$ assez grand. Fixons un r\'eel $\epsilon>0$ assez petit et un r\'eel $C>0$. Consid\'erons le domaine ${\cal D}$ des $X\in {\cal A}_{\tilde{L}}$ qui v\'erifient les deux conditions

(1) $\vert X\vert >C$;

(2) pour tout $\alpha\in \Sigma(A_{\tilde{L}})$, $\vert \alpha(X)\vert > \epsilon \vert X\vert $.

On fixe $C$ assez grand pour que $f_{\varphi_{X}}$ v\'erifie l'hypoth\`ese (4) du paragraphe pr\'ec\'edent pour tout $X\in {\cal D}$. Soit $\tilde{P}\in {\cal P}(\tilde{L})$, notons $\Delta_{\tilde{P}}$ la base de $\Sigma(A_{\tilde{L}})$ associ\'ee \`a $\tilde{P}$ et ${\cal C}_{\tilde{P}}\subset {\cal A}_{\tilde{L}}$ la chambre positive associ\'ee \`a $\tilde{P}$.  Supposons $X\in {\cal D}\cap {\cal C}_{\tilde{P}}$. On peut appliquer la formule (8) du paragraphe pr\'ec\'edent, qui, par changement de variable $H\mapsto H+X$, nous donne
$$(3) \qquad \ell(\varphi_{X})=-c\vert W(\tilde{L})\vert \int_{\tilde{T}({\mathbb R})^1}\int_{{\cal A}_{\tilde{L}}}D^{\tilde{G}}(exp(H+X)\gamma^1)^{1/2}$$
$$\Theta_{\tilde{\sigma}}(exp(H+X)\gamma^1)\phi(H )I^{\tilde{L}}(\gamma^1,f[\tilde{\pi}])\,dH\,d\gamma^1.$$
Puisque $\phi$ est \`a support compact, on a $H+X\in {\cal C}_{\tilde{P}}$ si  $\phi(H)\not=0$, pourvu que $C$ soit  choisi assez grand. On sait que l'on peut  alors d\'evelopper $D^{\tilde{G}}(exp(H+X)\gamma^1)^{1/2}\Theta_{\tilde{\sigma}}(exp(H+X)\gamma^1)$ comme somme finie
$$\sum_{j=1,...,k}q_{j}(H+X)e^{<H+X,\nu_{j}>}\Theta_{j}(\gamma^1),$$
o\`u les $q_{j}$ sont des polyn\^omes, les $\nu_{j}$ sont des \'el\'ements de ${\cal A}_{\tilde{L},{\mathbb C}}^*$ et les $\Theta_{j}$ sont des fonctions born\'ees sur $\tilde{T}({\mathbb R})^1$, $C^{\infty}$ sur le sous-ensemble des \'el\'ements r\'eguliers. On se rappelle que $\tilde{\sigma}=\tilde{\sigma}_{0,\lambda_{0}}$, o\`u $\tilde{\lambda}_{0}\in i{\cal A}_{\tilde{G}}^*$ et $\tilde{\sigma}_{0}$ est une repr\'esentation elliptique. Moeglin a prouv\'e que $\tilde{\sigma}_{0}$ \'etait supertemp\'er\'ee ([Moe]). Il en r\'esulte que les $Re(\nu_{j})$  sont des combinaisons lin\'eaires \`a coefficients n\'egatifs ou nuls des \'el\'ements de $\Delta_{\tilde{P}}$, l'un au moins des coefficients \'etant strictement n\'egatif.   Il r\'esulte alors de la formule (3) que $lim_{X\in {\cal D}\cap {\cal C}_{\tilde{P}},\vert X\vert \to \infty} \ell(\varphi_{X})=0$. Cela \'etant vrai pour tout $\tilde{P}$, on a plus simplement $lim_{X\in {\cal D}, \vert X\vert \to \infty}\ell(\varphi_{X})=0$. Mais il r\'esulte de 2.5(5) et de la d\'efinition de $\varphi_{X}$ que $\ell(\varphi_{X})$ est un polyn\^ome exponentiel en $X$, c'est-\`a-dire que
$$\ell(\varphi_{X})=\sum_{i=1,...,m}q_{i}(X)e^{<X,\lambda_{i}>},$$
pour des \'el\'ements $\lambda_{i}$ de ${\cal A}_{\tilde{L},{\mathbb C}}^*$ et des polyn\^omes $q_{i}$. Notons que cette \'egalit\'e est cette fois vraie pour tout $X$. Le lemme ci-dessous (appliqu\'e \`a l'ouvert $U$ des \'el\'ements de ${\cal A}_{\tilde{L}}$ qui v\'erifient (2)) entra\^{\i}ne que ce polyn\^ome exponentiel est identiquement nul. En particulier, pour $X=0$, $\ell(\varphi)=0$. Cela ach\`eve la preuve de la proposition 2.4 et du th\'eor\`eme 2.3.

\ass{Lemme}{Soit $V$ un espace  euclidien (de dimension finie), soit $U$ un ouvert de $V$  invariant par multiplication par ${\mathbb R}^{\times}$ et soit $f$ un polyn\^ome exponentiel sur $V$. Supposons $lim_{X\in U, \vert X\vert \to \infty}f(X)=0$. Alors $f=0$. }

Preuve. On \'ecrit comme ci-dessus
$$f(X)=\sum_{i=1,...,m}q_{i}(X)e^{<X,\lambda_{i}>}$$
o\`u les $\lambda_{i}$ sont distincts et les $q_{i}$ sont non nuls.   Si $f$ n'est pas nul (i.e. $m\geq1$), on peut trouver un  \'el\'ement $X_{0}\not=0$ de $U$ tel que  les produits $<X_{0},\lambda_{i}>$ soient deux-\`a-deux distincts et   $q_{i}(X_{0})\not=0$ pour tout $i$. La restriction $f_{D}$ de $f$ \`a la droite $D$ passant par  $X_{0}$ est alors un polyn\^ome exponentiel non nul qui v\'erifie $lim_{\vert X\vert \to \infty}f_{D}(X)=0$. Cette construction nous ram\`ene au cas o\`u $V$ est une droite,  auquel cas l'assertion est un simple exercice que l'on laisse au lecteur. $\square$

\bigskip

\subsection{L'espace $D_{spec}^{st}(\tilde{G}({\mathbb R}))$}
On note $D_{spec}^{st}(\tilde{G}({\mathbb R}))$ le sous-espace des \'el\'ements de $D_{spec}(\tilde{G}({\mathbb R}))$ qui sont des distributions stables, c'est-\`a-dire qui annulent $I^{inst}(\tilde{G}({\mathbb R}))$. Le corollaire suivant r\'esulte imm\'ediatement du th\'eor\`eme 2.3 et des relations 1.2(2) et 2.4(2).

\ass{Corollaire}{On a l'\'egalit\'e
 $$D_{spec}^{st}(\tilde{G}({\mathbb R}))=\left(\oplus_{\tilde{L}\in {\cal L}(\tilde{M}_{0})}Ind_{\tilde{L}}^{\tilde{G}}(D_{ell,{\mathbb C}}^{st}(\tilde{L}({\mathbb R})))\right)^{W(\tilde{M}_0)}.$$}
\bigskip

\subsection{L'espace $SI(\tilde{G}({\mathbb R}),K)$}
On note $SI(\tilde{G}({\mathbb R}),K)$ l'image de $C_{c}^{\infty}(\tilde{G}({\mathbb R}),K)$ dans $SI(\tilde{G}({\mathbb R}),K)$. On fixe des bases $B(\tilde{L})$ comme en 2.4. On note $PW^{st}(\tilde{G})$ le sous-espace de $PW^{\infty,st}(\tilde{G})$ form\'e des familles $(\varphi_{\tilde{L},\tilde{\pi}})_{\tilde{L}\in {\cal L}(\tilde{M}_{0}),\tilde{\pi}\in B(\tilde{L})}$ telles que l'ensemble des $(\tilde{L},\tilde{\pi})$ pour lesquels $\varphi_{\tilde{L},\tilde{\pi}}\not=0$ est fini.

\ass{Corollaire}{L'application $pw^{st}$ se restreint en un isomorphisme de $SI(\tilde{G}({\mathbb R}),K)$ sur $PW^{st}(\tilde{G})$.}

Preuve. D'apr\`es le th\'eor\`eme de Delorme et Mezo, $C_{c}^{\infty}(\tilde{G}({\mathbb R}),K)$ s'envoie par $pw$ sur $PW(\tilde{G})$. La projection dans $PW^{\infty,st}(\tilde{G})$ de cet espace est \'egale \`a $PW^{st}(\tilde{G})$. Donc $pw^{st}$ se restreint en une surjection de $SI(\tilde{G}({\mathbb R}),K)$ sur $PW^{st}(\tilde{G})$. Cette restriction est \'evidemment injective puisque $pw^{st}$ l'est. $\square$

\bigskip

\section{Transfert}

\bigskip

\subsection{D\'efinition d'un transfert spectral elliptique }
On revient au cas o\`u $(G,\tilde{G},\omega)$ est quelconque. Comme on l'a dit en 2.2, on fixe des mesures de Haar sur tous les groupes qui apparaissent. Soit ${\bf G}'=(G',{\cal G}',\tilde{s})$ une donn\'ee endoscopique elliptique et relevante de $(G,\tilde{G},\omega)$. Shelstad a prouv\'e l'existence du transfert 
$$\begin{array}{ccc}I(\tilde{G}({\mathbb R}),\omega)&\to& SI({\bf G}')\\ f&\mapsto& f^{{\bf G}'}\\ \end{array}$$
cf. [S]. Il envoie $I_{cusp}(\tilde{G}({\mathbb R}),\omega)$ dans $SI_{cusp}({\bf G}')$. 
Soit $\tilde{\sigma}\in D_{ell,0}^{st}({\bf G}')$, cf. 2.1. Consid\'erons l'application lin\'eaire sur $I_{cusp}(\tilde{G}({\mathbb R}),\omega)$ d\'efinie par $f\mapsto S^{{\bf G}'}(\tilde{\sigma},f^{{\bf G}'})$.    Le lemme 2.1 et le fait que $\tilde{\sigma}$ est $\mathfrak{Z}({\bf G}')$-finie implique que cette forme lin\'eaire est elle-m\^eme $\mathfrak{Z}(G)$-finie. On a d\'efini l'espace $D_{ell,0}^{st}({\bf G}')$ de sorte que cette forme lin\'eaire se factorise en une forme lin\'eaire sur $I_{cusp}(\tilde{G}({\mathbb R})/\mathfrak{A}_{\tilde{G}},\omega)$, cf. 2.1(2). 
Il existe donc un \'el\'ement $\tilde{\tau}\in D_{ell,0}(\tilde{G}({\mathbb R}),\omega)$, n\'ecessairement unique, tel que
$$I^{\tilde{G}}(\tilde{\tau},f)=S^{{\bf G}'}(\tilde{\sigma},f^{{\bf G}'})$$
pour tout $f\in I_{cusp}(\tilde{G}({\mathbb R}),\omega)$. Plus g\'en\'eralement, pour $\tilde{\sigma}$ comme ci-dessus et pour $\lambda\in {\cal A}_{\tilde{G}',{\mathbb C}}^*\simeq {\cal A}_{\tilde{G},{\mathbb C}}^*$, on a 
$$I^{\tilde{G}}(\tilde{\tau}_{\lambda},f)=S^{{\bf G}'}(\tilde{\sigma}_{\lambda},f^{{\bf G}'})$$
pour tout $f\in I_{cusp}(\tilde{G}({\mathbb R}),\omega)$. En pr\'ecisant le raisonnement ci-dessus, on voit qu'avec les d\'efinitions de 2.1, on a l'\'egalit\'e $\tilde{\mu}(\tilde{\tau}_{\lambda})=\tilde{\mu}(\tilde{\sigma}_{\lambda})$. 

 On appelle transfert spectral elliptique l'application lin\'eaire 
 $$\begin{array}{ccc}D_{ell}^{st}({\bf G}')&\to& D_{ell}(\tilde{G}({\mathbb R}),\omega)\\ \tilde{\sigma}_{\lambda}&\mapsto &\tilde{\tau}_{\lambda}\\ \end{array}$$
 o\`u ici, $\lambda$ appartient \`a $i{\cal A}_{\tilde{G}}^*$. 

\bigskip

\subsection{Le th\'eor\`eme}  
Les donn\'ees sont les m\^emes que dans le paragraphe pr\'ec\'edent.
\ass{Th\'eor\`eme}{Soit $\tilde{\sigma}\in D_{ell}^{st}({\bf G}')$, notons $\tilde{\tau}\in D_{ell}(\tilde{G}({\mathbb R}),\omega)$ son transfert spectral elliptique. Alors $\tilde{\tau}$ est le transfert de $\tilde{\sigma}$, c'est-\`a-dire que l'on a l'\'egalit\'e
$$I^{\tilde{G}}(\tilde{\tau},f)=S^{{\bf G}'}(\tilde{\sigma},f^{{\bf G}'})$$
pour tout $f\in I(\tilde{G}({\mathbb R}),\omega)$.}

Preuve. Pour tout $\tilde{L}\in {\cal L}(\tilde{M}_{0})$, on  fixe une base orthonorm\'ee $B(\tilde{L})$ de $D_{ell,0}(\tilde{L},\omega)$  et on r\'ealise l'espace $PW^{\infty}_{ell}(\tilde{L},\omega)$ \`a l'aide de cette base.  Les deux membres de  l'\'egalit\'e de l'\'enonc\'e sont continus en $f$. Gr\^ace au th\'eor\`eme de Paley-Wiener, il suffit de fixer $\tilde{L}\in {\cal L}(\tilde{M}_{0})$ et $\tilde{\pi}\in  B(\tilde{L})$ et de prouver que cette \'egalit\'e est v\'erifi\'ee pour $f\in pw^{-1}\circ sym^W(PW _{\tilde{\pi}}(\tilde{L},\omega))$, o\`u on utilise les m\^emes notations qu'en 2.4. C'est vrai par d\'efinition  du transfert spectral elliptique si $\tilde{L}=\tilde{G}$. On suppose donc $\tilde{L}\not=\tilde{G}$ et on raisonne par r\'ecurrence sur $a_{\tilde{L}}$. Puisque $\tilde{\tau}$ est elliptique et $\tilde{L}\not=\tilde{G}$, on a $I^{\tilde{G}}(\tilde{\tau},f)=0$ pour tout $f\in pw^{-1}\circ sym^W(PW _{\tilde{\pi}}(\tilde{L},\omega))$. Il s'agit de prouver que l'on a aussi $S^{{\bf G}'}(\tilde{\sigma},f^{{\bf G}'})=0$. On  utilise les m\^emes notations qu'en 2.4: on note ${\cal F}=PW_{\tilde{\pi}}(\tilde{L},\omega)$ et $f_{\varphi}=pw^{-1}\circ sym^W(\varphi)$ pour $\varphi\in {\cal F}$. On d\'efinit l'application lin\'eaire $\ell$ sur ${\cal F}$ par $\ell(\varphi)=S^{{\bf G}'}(\tilde{\sigma},(f_{\varphi})^{{\bf G}'})$. On a d\'efini la $W^{G'}$-orbite $\tilde{\mu}(\tilde{\sigma})$ dans $\tilde{\mu}(\omega)+\mathfrak{h}^{\theta,*}$. On note simplement $\mu_{\tilde{\sigma}}$ la $W$-orbite qu'elle engendre dans $\mathfrak{h}^*$. Le lemme 2.1 entra\^{\i}ne que l'on a l'\'egalit\'e
$$\ell(z\varphi)=z(\mu_{\tilde{\sigma}})\ell(\varphi)$$
pour tout $\varphi\in {\cal F}$ et $z\in \mathfrak{Z}(G)$. Le lemme 2.5 a de nouveau  les cons\'equences suivantes:

(1) si $(\mu(\tilde{\pi})+{\cal A}_{\tilde{L},{\mathbb C}}^*)\cap \mu_{\tilde{\sigma}}=\emptyset$, alors $\ell=0$.

Dans ce cas, on a fini. Supposons au contraire que $(\mu(\tilde{\pi})+{\cal A}_{\tilde{L},{\mathbb C}}^*)\cap \mu_{\tilde{\sigma}}\not=\emptyset$. Alors

(2) il y a un nombre fini de points $\lambda_{i}\in {\cal A}_{\tilde{L},{\mathbb C}}^*$, pour $i=1,...,m$, et, pour tout $i$, un op\'erateur diff\'erentiel $D_{i}$ de sorte que
$$\ell(\varphi)=\sum_{i=1,...,m}(D_{i}\varphi)(\lambda_{i}).$$

Soit $\varphi\in {\cal F}$. On introduit comme en 2.6 une fonction $f_{\varphi}^{\tilde{L}}\in I_{cusp}(\tilde{L}({\mathbb R}),\omega)$. Supposons que le support de cette fonction v\'erifie la condition 2.6(4). On construit alors l'\'el\'ement $g_{\varphi}\in I(\tilde{G}({\mathbb R}),\omega)$ qui v\'erifie des \'egalit\'es analogues \`a 2.6(5) et (6). Pr\'ecis\'ement, pour un espace de Levi $\tilde{M}\in {\cal L}(\tilde{M}_{0})$ et un \'el\'ement $\gamma\in \tilde{M}({\mathbb R})$ elliptique dans $\tilde{M}({\mathbb R})$ et fortement r\'egulier dans $\tilde{G}({\mathbb R})$, on a

(3) si $\tilde{M}$ n'est pas conjugu\'e \`a $\tilde{L}$, $I^{\tilde{G}}(\gamma,\omega,g_{\varphi})=0$;

(4) si $\tilde{M}=\tilde{L}$, 
$$I^{\tilde{G}}(\gamma,\omega,g_{\varphi})=\sum_{w\in W(\tilde{L})}\phi(w(H_{\tilde{L}}(\gamma)))I^{\tilde{L}}(\gamma,\omega,w^{-1}(f[\tilde{\pi}])).$$

  L'image $g_{\varphi}$ est somme de $f_{\varphi}$ et de termes $g_{\varphi,n}$, pour $n=a_{\tilde{G}},...a_{\tilde{L}}-1$, o\`u $pw(g_{\varphi,n})$ appartient \`a  
$$(\oplus_{\tilde{M}\in {\cal L}(\tilde{M}_{0}); a_{\tilde{M}}=n}PW^{\infty}_{ell}(\tilde{M},\omega))^{W(\tilde{M}_{0})}$$
  Parce que l'application $f\mapsto S^{{\bf G}'}(\tilde{\sigma},f)$ est continue, l'hypoth\`ese de r\'ecurrence nous dit que $S^{{\bf G}'}(\tilde{\sigma},g_{\varphi,n}^{{\bf G}'})=0$ pour tout $n=a_{\tilde{G}}+1,...,a_{\tilde{L}}-1$. En posant simplement $g_{\varphi,ell}=g_{\varphi,a_{\tilde{G}}}$, on obtient
  $$(5) \qquad \ell(\varphi)=S^{{\bf G}'}(\tilde{\sigma},f_{\varphi}^{{\bf G}'})=S^{{\bf G}'}(\tilde{\sigma},g_{\varphi}^{{\bf G}'})-S^{{\bf G}'}(\tilde{\sigma},g_{\varphi,ell}^{{\bf G}'}).$$
  Puisque $g_{\varphi,ell}$ est cuspidale, la d\'efinition du transfert spectral elliptique nous permet de r\'ecrire 
  $$\ell(\varphi)= S^{{\bf G}'}(\tilde{\sigma},g_{\varphi}^{{\bf G}'})-I^{\tilde{G}}(\tilde{\tau},g_{\varphi,ell})=S^{{\bf G}'}(\tilde{\sigma},g_{\varphi}^{{\bf G}'})-I^{\tilde{G}}(\tilde{\tau},g_{\varphi}).$$
 On fixe $\varphi\in {\cal F}$.  Soit $X\in {\cal A}_{\tilde{L}}$. On d\'efinit la fonction $\varphi_{X}$ comme en 2.7. Pour un domaine ${\cal D}$ comme dans ce paragraphe, la fonction $\varphi_{X}$ v\'erifie la condition 2.6(4) pour $X\in {\cal D}$. Parce que $\tilde{\tau}$ est elliptique, donc supertemp\'er\'ee, le m\^eme raisonnement qu'en 2.7 montre que
 $$(6) \qquad lim_{X\in {\cal D}; \vert X\vert \to \infty}I^{\tilde{G}}(\tilde{\tau},g_{\varphi_{X}})=0.$$
 A l'aide des formules  (3) et (4), on calcule ais\'ement les int\'egrales orbitales stables de la  fonction $g_{\varphi_{X}}^{{\bf G}'}$. On fixe des donn\'ees auxiliaires $G'_{1}$,...,$\Delta_{1}$ pour ${\bf G}'$, en supposant pour simplifier que le caract\`ere $\lambda_{\mathfrak{A}_{G'_{1}}}$ de 2.1 est trivial. Soit $\tilde{L}'$ un espace de Levi de $\tilde{G}'$. Notons $\tilde{L}'_{1}$ son image r\'eciproque dans $\tilde{G}'_{1}$. Soit $\delta_{1}\in \tilde{L}_{1}'({\mathbb R})$ un \'el\'ement $\tilde{G}$-r\'egulier et elliptique dans $\tilde{L}'_{1}({\mathbb R})$. Il lui correspond un sous-ensemble de  $\tilde{G}({\mathbb R})$ qui est soit vide, soit une classe de conjugaison stable d'\'el\'ements fortement r\'eguliers. Si cet ensemble ne contient 
  aucun \'el\'ement elliptique  de $\tilde{L}({\mathbb R})$, il r\'esulte de (3) que
 $$S^{G'_{1}}_{\lambda_{1}}(\delta_{1},g_{\varphi_{X}}^{\tilde{G}'_{1}})=0.$$
 Supposons qu'il corresponde \`a $\delta_{1}$ un \'el\'ement $\gamma\in \tilde{L}({\mathbb R})$ qui est elliptique. Cela entra\^{\i}ne que les espaces ${\cal A}_{\tilde{L}'}$ et ${\cal A}_{\tilde{L}}$ sont isomorphes  et que $H_{\tilde{L}}(\gamma)=H_{\tilde{L}'}(\delta)$ (cf. 2.1). Fixons un ensemble de repr\'esentants $(\gamma_{j})_{j=1,...,k}$ des classes de conjugaison par $L({\mathbb R})$ dans la classe de conjugaison stable de $\gamma$ dans $\tilde{L}({\mathbb R})$. On sait que c'est aussi un ensemble de repr\'esentants des classes de conjugaison par $G({\mathbb R})$ dans la classe de conjugaison stable de $\gamma$ dans $\tilde{G}({\mathbb R})$. De plus, $H_{\tilde{L}}(\gamma_{j})=H_{\tilde{L}}(\gamma)$ pour tout $j$.  On a alors
 $$S^{G'_{1}}_{\lambda_{1}}(\delta_{1},g_{\varphi_{X}}^{\tilde{G}'_{1}})=\sum_{j=1,...,k}[Z_{G}(\gamma_{j};{\mathbb R}):G_{\gamma_{j}}({\mathbb R})]^{-1}\Delta_{1}(\delta_{1},\gamma_{j})$$
 $$\sum_{w\in W(\tilde{L})}\phi(w(H_{\tilde{L}}(\gamma_{j}))-X)I^{\tilde{L}}(\gamma_{j},\omega,w^{-1}(f[\tilde{\pi}])).$$
 Pour tout $w\in W(\tilde{L})$, notons $g_{w}$ le transfert dans $SI_{\lambda_{1}}(\tilde{L}'_{1}({\mathbb R})/\mathfrak{A}_{\tilde{L}'_{1}})$ de la fonction $w^{-1}(f[\tilde{\pi}])$. On obtient l'\'egalit\'e
 $$S^{G'_{1}}_{\lambda_{1}}(\delta_{1},g_{\varphi_{X}}^{\tilde{G}'_{1}})=\sum_{w\in W(\tilde{L})}\phi(w(H_{\tilde{L}'}(\delta))-X)S^{\tilde{G}'_{1}}_{\lambda_{1}}(\delta_{1},g_{w}).$$
 A l'aide de ces formules, on peut adapter la preuve de 2.7: le fait que $\tilde{\sigma}$ soit supertemp\'er\'e entra\^{\i}ne que
 $$lim_{X\in {\cal D};\vert X\vert \to \infty}S^{{\bf G}'}(\tilde{\sigma},g_{\varphi_{X}}^{{\bf G}'})=0.$$
 Gr\^ace \`a (5) et (6), cela entra\^{\i}ne
 $$lim_{X\in {\cal D};\vert X\vert \to \infty}\ell(\varphi_{X})=0.$$
 De nouveau,   le lemme 2.7 montre que cette relation est contradictoire avec (2), sauf si $\ell=0$. Cela d\'emontre cette nullit\'e et le th\'eor\`eme. 
 
 \bigskip
 
 \subsection{Le transfert spectral}
  Fixons des donn\'ees auxiliaires $G'_{1}$,...,$\Delta_{1}$ pour notre donn\'ee endoscopique ${\bf G}'$. On a la variante du corollaire 2.8, avec des notations \'evidentes:
  $$D_{spec,\lambda_{1}}^{st}(\tilde{G}'_{1}({\mathbb R}))=\left(\oplus_{\tilde{L}'\in {\cal L}(\tilde{M}'_{0})}Ind_{\tilde{L}'_{1}}^{\tilde{G}'_{1}}(D_{ell, \lambda_{1},{\mathbb C}}^{st}(\tilde{L}'_{1}({\mathbb R})))\right)^{W^{G'}(\tilde{M}'_{0})}.$$
  Quand on fait varier les donn\'ees auxiliaires, on a une petite difficult\'e. Soit $\tilde{L}'\in {\cal L}(\tilde{M}'_{0})$. Deux cas se pr\'esentent. Si $\tilde{L}'$ est relevant pour $\tilde{G}$, les espaces $D_{ell, \lambda_{1},{\mathbb C}}^{st}(\tilde{L}'_{1}({\mathbb R}))$ se recollent naturellement en un espace que l'on peut noter $D_{ell,{\mathbb C}}^{st}({\bf L}')$. Si $\tilde{L}'$ n'est pas relevant, il n'y a pas de recollement intrins\`eque (c'est-\`a-dire ne d\'ependant que de $\tilde{L}'$ et pas du groupe ambiant $\tilde{G}'$) entre ces espaces. Par contre, les espaces induits $Ind_{\tilde{L}'_{1}}^{\tilde{G}'_{1}}(D_{ell, \lambda_{1},{\mathbb C}}^{st}(\tilde{L}'_{1}({\mathbb R})))$ se recollent. Cela ne nous g\^ene pas trop car il est clair que les transferts \`a $\tilde{G}({\mathbb R})$ des \'el\'ements de ces espaces sont  nuls. Notons ${\cal L}^{\tilde{G}}(\tilde{M}'_{0})$ le sous-ensemble des \'el\'ements de ${\cal L}(\tilde{M}'_{0})$ qui sont relevants pour $\tilde{G}$. Notons $D_{spec}^{st,\tilde{G}-nul}({\bf G}')$ la somme des images dans $D_{spec}^{st}({\bf G}')$ des espaces $Ind_{\tilde{L}'_{1}}^{\tilde{G}'_{1}}(D_{ell, \lambda_{1},{\mathbb C}}^{st}(\tilde{L}'_{1}({\mathbb R})))$ pour les $\tilde{L}'\in {\cal L}(\tilde{M}'_{0})-{\cal L}^{\tilde{G}}(\tilde{M}'_{0})$. On a alors
  $$D_{spec}^{st}({\bf G}')=\left(\oplus_{\tilde{L}'\in {\cal L}^{\tilde{G}}(\tilde{M}'_{0})}Ind_{{\bf L}'}^{{\bf G}'}(D_{ell,{\mathbb C}}^{st}({\bf L}'))\right)^{W^{G'}(\tilde{M}'_{0})}\oplus D_{spec}^{st,\tilde{G}-nul}({\bf G}').$$
  On d\'efinit une application lin\'eaire, que par anticipation, on note
  $$transfert:D_{spec}^{st}({\bf G}')\to D_{spec}(\tilde{G}({\mathbb R}),\omega)$$
  de la fa\c{c}on suivante. Elle est nulle sur $D_{spec}^{st,\tilde{G}-nul}({\bf G}')$. Soit $\tilde{L}'\in {\cal L}^{\tilde{G}}(\tilde{M}'_{0})$. On peut identifier ${\bf L}'$ \`a une donn\'ee endoscopique de $(L,\tilde{L},{\bf a})$, o\`u $\tilde{L}$ est un \'el\'ement de ${\cal L}(\tilde{M}_{0})$. On a le transfert spectral elliptique, qui se prolonge en une application lin\'eaire $D_{ell,{\mathbb C}}^{st}({\bf L}')\to D_{ell,{\mathbb C}}(\tilde{L}({\mathbb R}),\omega)$. Par induction, on en d\'eduit une application 
 lin\'eaire
 $$Ind_{{\bf L}'}^{{\bf G}'}(D_{ell,{\mathbb C}}^{st}({\bf L}'))\to Ind_{\tilde{L}}^{\tilde{G}}(D_{ell,{\mathbb C}}(\tilde{L}({\mathbb R}),\omega)).$$
 Le dernier espace s'envoie naturellement dans  $D_{spec}(\tilde{G}({\mathbb R}),\omega)$. Alors l'application $transfert$ co\"{\i}ncide sur $Ind_{{\bf L}'}^{{\bf G}'}(D_{ell,{\mathbb C}}^{st}({\bf L}'))$ avec le compos\'e des deux applications pr\'ec\'edentes. Puisque le transfert commute \`a l'induction, le corollaire suivant r\'esulte imm\'ediatement du th\'eor\`eme 3.2.
 
 \ass{Corollaire}{L'application $transfert$ est bien le transfert endoscopique, c'est-\`a-dire que, pour tout $\tilde{\sigma}\in D_{spec}^{st}({\bf G}')$ et tout $f\in I(\tilde{G}({\mathbb R}),\omega)$, on a l'\'egalit\'e
 $$I^{\tilde{G}}(transfert(\tilde{\sigma}),f)=S^{{\bf G}'}(\tilde{\sigma},f^{{\bf G}'}).$$ }
 
 {\bf Remarque.} Si l'on r\'etablit plus canoniquement les espaces de mesures, l'application transfert devient une application lin\'eaire
$$transfert:D_{spec}^{st}({\bf G}')\otimes Mes(G'({\mathbb R}))^*\to D_{spec}(\tilde{G}({\mathbb R}),\omega)\otimes Mes(G({\mathbb R}))^*.$$

 \bigskip
 
 \subsection{Transfert $K$-fini}
  Soit $\Omega$ un ensemble fini de $K$-types, c'est-\`a-dire de repr\'esentations irr\'eductibles de $K$. Pour une fonction $f$ sur $\tilde{G}({\mathbb R})$ et pour $k\in K$, notons $\lambda_{k}(f)$, resp. $\rho_{k}(f)$, la fonction $\gamma\mapsto f(k^{-1}\gamma)$, resp. $\gamma\mapsto f(\gamma k)$. On dit que $f$ se transforme \`a gauche, resp. \`a droite selon $\Omega$ si la repr\'esentation de $K$ dans l'ensemble de fonctions $\{\lambda_{k}(f);k\in K\}$, resp. $\{\rho_{k}(f);k\in K\}$, se d\'ecompose en repr\'esentations irr\'eductibles appartenant \`a $\Omega$. On note $C_{c}^{\infty}(\tilde{G}({\mathbb R}),\Omega)$ l'espace des fonctions qui se transforment \`a droite et \`a gauche selon $\Omega$. On note $I(\tilde{G}({\mathbb R}),\omega,\Omega)$ son image dans $I(\tilde{G}({\mathbb R}),\omega)$. L'espace $C_{c}^{\infty}(\tilde{G}({\mathbb R}),K)$ est la r\'eunion des $C_{c}^{\infty}(\tilde{G}({\mathbb R}),\Omega)$ quand $\Omega$ parcourt tous les ensembles finis de $K$-types.  
 
 Soit $(\pi,\tilde{\pi})$ une $\omega$-repr\'esentation telle que $\pi$ soit irr\'eductible ou, plus g\'en\'eralement telle que $\pi$ soit de longueur finie et que toutes ses composantes irr\'eductibles aient un m\^eme param\`etre infinit\'esimal. Notons $\mu$ ce param\`etre, qui est une orbite dans $\mathfrak{h}^*$ pour l'action de $W$. On a vu en 1.2 que l'intersection $(\tilde{\mu}(\omega)+\mathfrak{h}^{\theta,*})\cap \mu$ \'etait une unique orbite sous l'action de $W^{\theta}$. Notons-la $\tilde{\mu}(\tilde{\pi})$.   Soient $\tilde{L}\in {\cal L}(\tilde{M}_{0})$ et $\tilde{\pi}\in {\cal E}_{ell,0}(\tilde{L},\omega)$. A $\tilde{\pi}$ est associ\'e une orbite  $\tilde{\mu}(\tilde{\pi})\subset \tilde{\mu}(\omega)+\mathfrak{h}^{\theta,*}$ pour l'action de $W^{L,\theta}$. Notons $\tilde{\mu}(\tilde{\pi})^G$ l'orbite pour l'action de $W^{\theta}$ engendr\'ee par  $\tilde{\mu}(\tilde{\pi})$. Appelons "type spectral" une orbite dans $\tilde{\mu}(\omega)+\mathfrak{h}^{\theta,*}$ pour l'action de $W^{\theta}$. R\'ealisons l'espace $PW^{\infty}(\tilde{G},\omega)$ en fixant des bases comme en 2.2. Pour $\tilde{L}\in {\cal L}(\tilde{M}_{0})$, on note $B(\tilde{L})$ la base fix\'ee de $D_{ell,0}(\tilde{L}({\mathbb R}),\omega)$. Pour un type spectral $\tilde{\mu}$, d\'efinissons le sous-espace $PW(\tilde{G},\omega,\tilde{\mu})$ des familles $(\varphi_{\tilde{L},\tilde{\pi}})_{\tilde{L}\in {\cal L}(\tilde{M}_{0}),\tilde{\pi}\in B(\tilde{L})}$  appartenant \`a $PW^{\infty}(\tilde{G},\omega)$ v\'erifiant la condition suivante:
 
 - pour $\tilde{L}\in {\cal L}(\tilde{M}_{0})$ et $\tilde{\pi}\in B(\tilde{L})$, on a $\varphi_{\tilde{L},\tilde{\pi}}=0$ si $\tilde{\mu}(\tilde{\pi})^G\not=\tilde{\mu}$.
 
 Il est clair que $PW(\tilde{G},\omega)$ est la somme directe des $PW(\tilde{G},\omega,\tilde{\mu})$ sur tous les types spectraux $\tilde{\mu}$ (et $PW^{\infty}(\tilde{G},\omega)$ est le compl\'et\'e de cette somme). Plus g\'en\'eralement, pour un ensemble fini $\Omega^{pw}$ de types spectraux, notons $PW(\tilde{G},\omega,\Omega^{pw})$ la somme des $PW(\tilde{G},\omega,\tilde{\mu})$ pour $\tilde{\mu}\in \Omega^{pw}$. 
 Le th\'eor\`eme de Delorme et Mezo affirme pr\'ecis\'ement que
  
  (1) pour tout ensemble fini $\Omega$ de $K$-types, il existe un ensemble fini $\Omega^{pw}$ de types spectraux tel que $pw_{\tilde{G},\omega}(I(\tilde{G}({\mathbb R}),\omega,\Omega))\subset PW(\tilde{G},\omega,\Omega^{pw})$;
  
  (2) pour tout ensemble fini $\Omega^{pw}$ de types spectraux, il existe un ensemble fini $\Omega$ de $K$-types tel que $PW(\tilde{G},\omega,\Omega^{pw})\subset pw_{\tilde{G},\omega}(I(\tilde{G}({\mathbb R}),\omega,\Omega)$. 
  
  Dans le cas o\`u $(G,\tilde{G},{\bf a})$ est quasi-d\'eploy\'e et \`a torsion int\'erieure, on a les variantes stables $SI(\tilde{G}({\mathbb R}),\Omega)$ et $PW^{st}(\tilde{G},\Omega^{pw})$. En reprenant la preuve du corollaire 2.9, on obtient des variantes de (1) et (2) pour ces variantes stables.

 Fixons des donn\'ees auxiliaires $G'_{1}$,...,$\Delta_{1}$ pour notre donn\'ee endoscopique ${\bf G}'$. On fixe un sous-groupe compact maximal $K'_{1}$ de $G'_{1}({\mathbb R})$. Les constructions ci-dessus s'adaptent \`a ces donn\'ees auxiliaires, en consid\'erant alors des fonctions et des repr\'esentations qui se transforment selon le caract\`ere $\lambda_{1}$ de $C_{1}({\mathbb R})$. Il convient peut-\^etre de parler alors de $\lambda_{1}$-type spectral
 
 \ass{Corollaire}{(i) Pour toute fonction $f\in C_{c}^{\infty}(\tilde{G}({\mathbb R}))$ qui est $K$-finie, le transfert $f^{{\bf G}'}$ est l'image dans $SI({\bf G}')$ d'une fonction $K_{1}'$-finie dans $C_{c,\lambda_{1}}^{\infty}(\tilde{G}'_{1}({\mathbb R}))$.
 
 (ii) Soit $\Omega$ un ensemble fini de $K$-types. Alors il existe un ensemble fini $\Omega'_{1}$ de $K'_{1}$-types tel que, pour tout  $f\in C_{c}^{\infty}(\tilde{G}({\mathbb R}),\Omega)$, le transfert $f^{{\bf G}'}$ est l'image dans $SI({\bf G}')$ d'un \'el\'ement de $C_{c,\lambda_{1}}^{\infty}(\tilde{G}'_{1}({\mathbb R}),\Omega'_{1})$.
 
 (iii) Soit $\Omega'_{1}$ un ensemble fini de $K'_{1}$-types. Alors il existe un ensemble fini $\Omega$ de $K$-types de sorte que la condition suivante soit v\'erifi\'ee. Soit $f_{0}\in C_{c}^{\infty}(\tilde{G}({\mathbb R}))$, supposons que son  transfert $f_{0}^{{\bf G}'}$ soit l'image dans $SI({\bf G}')$ d'un \'el\'ement de $C_{c,\lambda_{1}}^{\infty}(\tilde{G}'_{1}({\mathbb R}),\Omega'_{1})$. Alors il existe $f\in C_{c}^{\infty}(\tilde{G}({\mathbb R}),\Omega)$ tel que $f^{{\bf G}'}=f_{0}^{{\bf G}'}$.}
 
Preuve. Evidemment, (i) r\'esulte de (ii).  Le transfert s'identifie \`a une application lin\'eaire
$$transfert:PW^{\infty}(\tilde{G},\omega)\to PW_{\lambda_{1}}^{\infty,st}(\tilde{G}'_{1}).$$
Notons $Im$ son image. Gr\^ace aux propri\'et\'es (1) et (2) ci-dessus et \`a leurs analogues stables, les assertions (ii) et (iii) r\'esultent des propri\'et\'es suivantes

(3) soit $\Omega^{pw}$ un ensemble fini de types spectraux pour $\tilde{G}$; alors il existe un ensemble fini $\Omega^{_{'}pw}$ de $\lambda_{1}$-types spectraux pour $\tilde{G}'_{1}$ de sorte que
$$transfert(PW(\tilde{G},\omega,\Omega^{pw}))\subset PW_{\lambda_{1}}^{st}(\tilde{G}'_{1},\Omega^{_{'}pw});$$

(4) soit $\Omega^{_{'}pw}$ un ensemble fini de $\lambda_{1}$-types spectraux pour $\tilde{G}'_{1}$; alors il existe un ensemble fini $\Omega^{pw}$ de types spectraux  pour $\tilde{G}$ de sorte que
$$Im\cap PW_{\lambda_{1}}^{st}(\tilde{G}'_{1},\Omega^{_{'}pw})\subset transfert(PW(\tilde{G},\omega,\Omega^{pw})).$$

Un $\lambda_{1}$-type spectral $\tilde{\mu}'$ pour $\tilde{G}'_{1}$ est une orbite pour l'action de $W^{G'}$ dans 
$\mathfrak{h}^{{\bf G}',*}$. Cet ensemble est isomorphe \`a $\tilde{\mu}(\omega)+\mathfrak{h}^{\theta,*}$, l'isomorphisme \'etant compatible avec les actions de $W^{G'}$ et $W^{\theta}$ et le plongement de $W^{G'}$ dans $W^{\theta}$. On peut donc d\'efinir l'orbite pour  $W^{\theta}$ engendr\'ee par $\tilde{\mu}'$, que l'on note $(\tilde{\mu}')^G$. On obtient une application $q:\tilde{\mu}'\mapsto (\tilde{\mu}')^G$ de l'ensemble des $\lambda_{1}$-types spectraux pour $\tilde{G}'_{1}$ dans l'ensemble des types spectraux pour $\tilde{G}$. Cette application est \`a fibres finies. Montrons que

(5) pour tout type spectral $\tilde{\mu}$ pour $\tilde{G}$, on a l'inclusion 
$$transfert(PW(\tilde{G},\omega,\tilde{\mu})\subset PW_{\lambda_{1}}^{st}(\tilde{G}'_{1},q^{-1}(\tilde{\mu})).$$

  Soit $(\varphi_{\tilde{L},\tilde{\pi}})_{\tilde{L}\in {\cal L}(\tilde{M}_{0}),\tilde{\pi}\in B(\tilde{L})}$  un \'el\'ement de $PW^{\infty}(\tilde{G},\omega)$, notons  $(\varphi_{\tilde{L}'_{1},\tilde{\sigma}_{1}})_{\tilde{L}'\in {\cal L}(\tilde{M}'_{0}),\tilde{\sigma}_{1}\in B^{st}(\tilde{L}'_{1})}$ son transfert (avec une notation \'evidente). D'apr\`es le corollaire 3.3, on a la description suivante. Soit $\tilde{L}'\in {\cal L}(\tilde{M}'_{0})$. Si $\tilde{L}'$   ne correspond pas \`a un espace de Levi de $\tilde{G}$, alors $\varphi_{\tilde{L}'_{1},\tilde{\sigma}_{1}}=0$ pour tout $\tilde{\sigma}_{1}$. 
 Supposons que $\tilde{L}'$ soit associ\'e \`a un espace de Levi $\tilde{L}\in {\cal L}(\tilde{M}_{0})$. Pour $\tilde{\sigma}_{1}\in B^{st}(\tilde{L}'_{1})$, on peut \'ecrire $transfert(\tilde{\sigma}_{1})$ comme une combinaison lin\'eaire finie d'\'el\'ements de la base $B(\tilde{L})$:
$$ transfert(\tilde{\sigma}_{1})=\sum_{\tilde{\pi}\in B(\tilde{L})}c_{\tilde{\pi}}\tilde{\pi}.$$
 On a alors
$$ \varphi_{\tilde{L}'_{1},\tilde{\sigma}_{1}}=\sum_{\tilde{\pi}\in B(\tilde{L})}c_{\tilde{\pi}}\varphi_{\tilde{L},\tilde{\pi}}.$$
  Par compatibilit\'e du transfert avec les actions des centres de l'alg\`ebre enveloppante, on sait que les \'el\'ements $\tilde{\pi}$ qui apparaissent ont un caract\`ere central param\'etr\'e par la $W^L$-orbite dans $\mathfrak{h}^*$ engendr\'ee par la $W^{L'}$-orbite $\mu(\tilde{\sigma}_{1})$ dans $\mathfrak{h}^{{\bf G}',*}\simeq \tilde{\mu}(\omega)+\mathfrak{h}^{\theta,*}$. Il en r\'esulte que $\tilde{\mu}(\tilde{\pi})=(\tilde{\mu}(\tilde{\sigma}_{1}))^{G}$.  Supposons que la famille $(\varphi_{\tilde{L},\tilde{\pi}})_{\tilde{L}\in {\cal L}(\tilde{M}_{0}),\tilde{\pi}\in B(\tilde{L})}$ appartienne \`a $PW(\tilde{G},\omega,\tilde{\mu})$. Si $\tilde{\mu}(\tilde{\sigma}_{1})\not\in q^{-1}(\tilde{\mu})$,  le coefficient $c_{\tilde{\pi}}$ est nul ou $\varphi_{\tilde{L},\tilde{\pi}}$ est nul. Donc $\varphi_{\tilde{L}'_{1},\tilde{\sigma}_{1}}=0$. Cela prouve que la famille $(\varphi_{\tilde{L}'_{1},\tilde{\sigma}_{1}})_{\tilde{L}'\in {\cal L}(\tilde{M}'_{0}),\tilde{\sigma}_{1}\in B^{st}(\tilde{L}'_{1})}$ appartient \`a $PW_{\lambda_{1}}^{st}(\tilde{G}'_{1},q^{-1}(\tilde{\mu}))$. Cela prouve (5). 
  
   L'assertion (3) en r\'esulte imm\'ediatement: il suffit de prendre pour $\Omega^{_{'}pw}$ la r\'eunion des $q^{-1}(\tilde{\mu})$ pour $\mu\in \Omega^{pw}$. 
   
   Pour tout type spectral $\tilde{\mu}$ pour $\tilde{G}$, on a d\'efini $PW(\tilde{G},\omega,\tilde{\mu})$ comme un sous-espace de $PW^{\infty}(\tilde{G},\omega)$. On a un projecteur $p_{\tilde{\mu}}:PW^{\infty}(\tilde{G},\omega)\to PW(\tilde{G},\omega,\tilde{\mu})$. Il associe \`a une famille $(\varphi_{\tilde{L},\tilde{\pi}})_{\tilde{L}\in {\cal L}(\tilde{M}_{0}),\tilde{\pi}\in B(\tilde{L})}$ la famille $(\varphi'_{\tilde{L},\tilde{\pi}})_{\tilde{L}\in {\cal L}(\tilde{M}_{0}),\tilde{\pi}\in B(\tilde{L})}$, o\`u 
   $$\varphi'_{\tilde{L},\tilde{\pi}}=\left\lbrace\begin{array}{cc}\varphi_{\tilde{L},\tilde{\pi}},&\text{ si }\tilde{\mu}(\tilde{\pi})=\tilde{\mu};\\ 0,&\text{ sinon.}\\ \end{array}\right.$$
   Plus g\'en\'eralement, on a un projecteur $p_{\Omega^{pw}}:PW^{\infty}(\tilde{G},\omega)\to PW(\tilde{G},\omega,\Omega^{pw})$ pour tout ensemble fini $\Omega^{pw}$ de types spectraux. La preuve de (5) prouve plus g\'en\'eralement que, pour tout tel ensemble, on a
   $$transfert\circ p_{\Omega^{pw}}=p_{q^{-1}(\Omega^{pw})}\circ transfert.$$
    Soit $\Omega^{_{'}pw}$ un ensemble fini de $\lambda_{1}$-types spectraux pour $\tilde{G}'_{1}$ , posons $\Omega^{pw}=q(\Omega^{_{'}pw})$.   Soit $\boldsymbol{\varphi}'\in Im\cap PW_{\lambda_{1}}^{st}(\tilde{G}'_{1},\Omega^{_{'}pw})$. Introduisons un \'el\'ement $\boldsymbol{\varphi}\in PW^{\infty}(\tilde{G},\omega)$ tel que $transfert(\boldsymbol{\varphi})=\boldsymbol{\varphi}'$. Parce que $\boldsymbol{\varphi}'\in PW_{\lambda_{1}}^{st}(\tilde{G}'_{1},\Omega^{_{'}pw})$, on a $p_{q^{-1}(\Omega^{pw})}(\boldsymbol{\varphi}')=\boldsymbol{\varphi}'$. Donc
    $$\boldsymbol{\varphi}'=p_{q^{-1}(\Omega^{pw})}(\boldsymbol{\varphi}')=p_{q^{-1}(\Omega^{pw})}\circ transfert(\boldsymbol{\varphi})=transfert\circ p_{\Omega^{pw}}(\boldsymbol{\varphi}).$$
    Donc $\boldsymbol{\varphi}'$ est le transfert de l'\'el\'ement $p_{\Omega^{pw}}(\boldsymbol{\varphi})\in PW(\tilde{G},\omega,\Omega^{pw})$. Cela prouve (4) et le corollaire.
  $\square$

\bigskip

\subsection{Transfert $K$-fini, version g\'en\'erale}
Consid\'erons dans ce paragraphe un $K$-espace $K\tilde{G}$, cf. [I] 1.11. On fixe pour chaque composante connexe $K\tilde{G}_{p}$ un espace de Levi minimal $\tilde{M}_{p,0}$ et un sous-groupe compact maximal $K_{p}$ de $G_{p}({\mathbb R})$ en bonne position relativement \`a $\tilde{M}_{p,0}$. Les d\'efinitions du paragraphe pr\'ec\'edent se g\'en\'eralisent imm\'ediatement. Un ensemble fini de $K$-types est la r\'eunion disjointe d'ensembles finis de $K_{p}$-types. On a une application de transfert
$$I(K\tilde{G}({\mathbb R}),\omega)\to \oplus_{{\bf G}'\in {\cal E}(\tilde{G},{\bf a})}SI({\bf G}').$$
On renvoie \`a [I] 4.11 pour les notations (on a supprim\'e les espaces
 de mesures, les mesures \'etant fix\'ees depuis 2.2). Pour tout ${\bf G}'\in {\cal E}(\tilde{G},{\bf a})$, fixons des donn\'ees auxiliaires $G'_{1}$,...,$\Delta_{1}$ et un sous-groupe compact maximal $K^{G'_{1}}$ de $G'_{1}({\mathbb R})$.
 
 \ass{Corollaire}{(i) Soit $\Omega$ un ensemble fini de $K$-types. Alors il existe pour tout ${\bf G}'\in {\cal E}(\tilde{G},{\bf a})$ un ensemble fini de $K^{G'_{1}}$-types $\Omega^{G'_{1}}$ de sorte que, pour tout $f\in C_{c}^{\infty}(K\tilde{G}({\mathbb R}),\omega,\Omega)$ et tout ${\bf G}'\in {\cal E}(\tilde{G},{\bf a})$, le transfert $f^{{\bf G}'}$ soit l'image dans $SI({\bf G}')$ d'un \'el\'ement de $C_{c,\lambda_{1}}^{\infty}(\tilde{G}'_{1}({\mathbb R}),\Omega^{G'_{1}})$.   
  
  (ii) Pour tout ${\bf G}'\in {\cal E}(\tilde{G},{\bf a})$, soit $\Omega^{G'_{1}}$ un ensemble fini de $K^{G'_{1}}$-types. Alors il existe un ensemble fini $\Omega$ de $K$-types de sorte que la condition suivante soit v\'erifi\'ee. Soit $f_{0}\in C_{c}^{\infty}(K\tilde{G}({\mathbb R}))$, supposons que, pour tout ${\bf G}'\in {\cal E}(\tilde{G},{\bf a})$, le  transfert $f^{{\bf G}'}$ soit l'image dans $SI({\bf G}')$ d'un \'el\'ement de $C_{c,\lambda_{1}}^{\infty}(\tilde{G}'_{1}({\mathbb R}),\Omega^{G'_{1}})$. Alors il existe un \'el\'ement $f\in C_{c}^{\infty}(K\tilde{G}({\mathbb R}),\Omega)$ tel que $f^{{\bf G}'}=f_{0}^{{\bf G}'}$ pour tout ${\bf G}'\in {\cal E}(\tilde{G},{\bf a})$.}
  
  Preuve. Le (i) se d\'eduit imm\'ediatement du (ii) du corollaire pr\'ec\'edent. Le (ii) se prouve de la m\^eme fa\c{c}on que le (iii) de ce corollaire. On laisse les d\'etails au lecteur. $\square$
  
  \bigskip
  
  \subsection{Le cas du corps de base ${\mathbb C}$}
  Dans tout l'article, le corps de base \'etait ${\mathbb R}$. Rempla\c{c}ons-le maintenant par ${\mathbb C}$. Les m\^emes d\'efinitions et th\'eor\`emes restent valables. En effet, consid\'erons un triplet $(G,\tilde{G},{\bf a})$ d\'efini sur ${\mathbb C}$. Par restriction des scalaires, on en d\'eduit un triplet $(G_{{\mathbb R}},\tilde{G}_{{\mathbb R}},{\bf a}_{{\mathbb R}})$ sur ${\mathbb R}$. Il suffit alors d'appliquer les th\'eor\`emes \`a   ce triplet. En fait, le cas complexe est beaucoup plus simple. Il n'y a de fonctions cuspidales sur $\tilde{G}({\mathbb C})$ que si $G$ est un tore. Le groupe $G$ est forc\'ement d\'eploy\'e. Dans le cas o\`u $(G,\tilde{G},{\bf a})$ est  \`a torsion int\'erieure,  on a l'\'egalit\'e $SI(\tilde{G}({\mathbb C}))=I(\tilde{G}({\mathbb C}))$.

 \bigskip 
{\bf Bibliographie}

[A]  J. Arthur  {\it On local character relations}, Selecta Math. 2 (1996), p. 501-579

[BR] G. Barban\c{c}on, M. Ra\"{\i}s: {\it Sur le th\'eor\`eme de Hilbert diff\'erentiablle pour les groupes lin\'eaires finis (d'apr\`es E. Noether)}, Ann. Sc. ENS 16 (1983), p. 355- 373

[BT] A. Borel, J. Tits: {\it Groupes r\'eductifs},  Publ. Math. IHES 27 (1965), p. 55- 150

[DM]  P. Delorme, P. Mezo: {\it A twisted invariant Paley-Wiener theorem for real reductive groups}, Duke Math. J. 144 (2008), p.341-380

[Me] P. Mezo: {\it Spectral transfer in the twisted endoscopy of real groups}, pr\'epublication 2013

[Moe] C. Moeglin: {\it Repr\'esentations elliptiques; caract\'erisation et formules de transfert  de caract\`eres}, pr\'epublication 2013

[R] D. Renard:  {\it Int\'egrales orbitales tordues sur les groupes de Lie r\'eductifs r\'eels}, J. Funct. Analysis 145 (1997), p. 374- 454

[S] D. Shelstad: {\it On geometric transfer in real twisted endoscopy}, pr\'epublication 2011

[W] J.-L. Waldspurger: {\it La formule des traces locale tordue}, pr\'epublication 2012

[I] -------------------------: {\it Stabilisation de la formule des traces tordue I: endoscopie tordue sur un corps local}, pr\'epublication 2014

 \bigskip
 
 Institut de Math\'ematiques de Jussieu -CNRS
 
 2 place Jussieu 75005 Paris
 
 e-mail: waldspur@math.jussieu.fr

\end{document}